\documentclass[a4paper,11pt,reqno]{article}

\usepackage[utf8]{inputenc}
\usepackage[T1]{fontenc}
\usepackage{lmodern}
\usepackage[english]{babel}
\usepackage{microtype}
\usepackage{csquotes}

\usepackage[backend=biber,hyperref=true,style=alphabetic,giveninits=true]{biblatex}
\usepackage{biblatex-mr}
\usepackage{tikz}

\definecolor{midnightblue}{rgb}{0.1, 0.1, 0.44}
\usepackage[colorlinks,allcolors=midnightblue]{hyperref}
\AtBeginDocument{}
\AtBeginDocument{}
\AtBeginDocument{}

\usepackage[a4paper,hmargin=2.5cm,bottom=3cm,top=3cm,footskip=3\baselineskip,
	marginparwidth=2.4cm,marginparsep=0.05cm]{geometry}

\usepackage{fsmath}
\usepackage{bm}

\usepackage[normalem]{ulem}
\usepackage{enumerate}
\usepackage{footmisc}

\usepackage{array}
\usepackage{booktabs}
\newcommand{\otoprule}{\midrule[\heavyrulewidth]}

\newcommand{\loc}{\mathrm{loc}}
\DeclareMathOperator{\TV}{TV}

\newcommand{\dx}{\d x}
\newcommand{\dy}{\d y}
\newcommand{\dz}{\d z}
\newcommand{\dt}{\d t}
\newcommand{\ds}{\d s}

\newcommand{\II}{I\!I}
\newcommand{\III}{I\!I\!I}

\newcommand{\Intq}{\int_0^T \!\!\!\! \int_0^T \!\!\!\! \int_{\setR^n} \!\! \int_{\setR^n}}
\newcommand{\Intt}{\int_0^T \!\!\!\! \int_{\setR^n} \!\! \int_{\setR^n}}
\newcommand{\Intd}{\int_0^T \!\!\!\! \int_{\setR^n}}

\newcommand{\intt}{\int_{t_1}^{t_2} \!\!\!\! \int_{\setR^n} \!\! \int_{\setR^n}}
\newcommand{\intd}{\int_{t_1}^{t_2} \!\!\!\! \int_{\setR^n}}
\newcommand{\intrr}{\int_{\setR^n} \!\! \int_{\setR^n}}

\begin{document}

\title{Stability of quasi-entropy solutions \\ of non-local scalar conservation laws}
\author{
	Elio Marconi\thanks{Università degli Studi di Padova, \url{elio.marconi@unipd.it}}
	\and
	Emanuela Radici\thanks{Università degli Studi dell'Aquila, \url{emanuela.radici@univaq.it}}
	\and
	Federico Stra\thanks{Politecnico di Torino, \url{federico.stra@polito.it}}
}

\maketitle

\begin{abstract}
We prove the stability of entropy solutions of nonlinear conservation laws with respect to perturbations of the initial datum, the space-time dependent flux and the entropy inequalities.

Such a general stability theorem is motivated by the study of problems in which the flux $P[u](t,x,u)$ depends possibly non-locally on the solution itself. For these problems we show the conditional existence and uniqueness of entropy solutions.

Moreover, the relaxation of the entropy inequality allows to treat approximate solutions arising from various numerical schemes. This can be used to derive the rate of convergence of the recent particle method introduced in \cite{RadiciStra} to solve a one-dimensional model of traffic with congestion, as well as recover already known rates for some other approximation methods.
\end{abstract}

\tableofcontents

\vspace*{2cm}

\section{Introduction}

The study of conservation laws $\partial_t u + \div_x\bigl(P(t,x,u)\bigr) = 0$ has been initiated by \cite{Kruzkov}, who introduced the notion of entropy solutions as a selecting criterion among the more general distributional solutions. This notion is encoded in the distributional inequalities
\[
\partial_t\abs{u-c} + \div_x\bigl[\sign(u-c)\bigl(P(t,x,u)-P(t,x,c)\bigr)\bigr]
	+ \sign(u-c) \div_x\!P(t,x,c)
\leq 0, \qquad \forall c\in\setR.
\]

By means of his celebrated \emph{doubling of variables} technique, in \cite{Kruzkov} the author shows the $L^1_\loc$ stability of $L^\infty$ entropy solutions with respect to perturbations of the initial datum.

The stability has been extended to allow the treatment of approximate solutions, for instance arising from numerical methods. This notion of quasi-entropy solutions amounts to introducing error terms in the right hand side of the entropy inequalities.

In \cite{Kuznetsov} this flexibility is exploited to derive convergence rates for several numerical methods in the particular case of a flux $P(u)$ depending only on the density.
Following this research direction, \cite{BouchutPerthame1998} codified a more general notion of quasi-entropy solutions where the error terms in the right hand side of the entropy inequalities are derivatives with respect to $t$ and $x$ of measures. Their result is still restricted to the case of a flux of the form $P(u)$, which is a severe limitation in many applications.

An improvement in this direction, but limited to exact entropy solutions, is given by \cite[Theorem 1.3]{KarlsenRisebro}, where the authors show the stability of entropy solutions with respect to two distinct fluxes, which are taken to be of product form $P(x,u)=f(u)k(x)$. Dealing with two distinct fluxes forces to consider $BV$ solutions instead of merely $L^\infty$, as will be pointed out when we discuss the strategy of our proof.

\begin{table}[h!]
\caption{Summary of previous stability theorems.}
\centering
\begin{tabular}[t]{lcclc}
\toprule
 & no. of fluxes & type of flux & type of solutions & regularity in space \\\otoprule
\cite{Kruzkov} & 1 & $P(t,x,u)$ & entropy & $L^\infty$ \\\midrule
\cite{Kuznetsov} & 1 & $P(u)$ & quasi-entropy, numerical schemes & $L^\infty$ \\\midrule
\cite{BouchutPerthame1998} & 1 & $P(u)$ & quasi-entropy, more general & $L^\infty$ \\\midrule
\cite{KarlsenRisebro} & 2 & $f(u)k(x)$ & entropy & $BV$ \\\midrule
Ours & 2 & $P(t,x,u)$ & quasi-entropy & $BV$ \\\bottomrule
\end{tabular}
\end{table}

\paragraph{Goal of the article.}
Merging the various lines of improvement, we want to obtain a stability result for quasi-entropy solutions with distinct fluxes depending also on $t$ and $x$.
The notion of $\mu$-quasi-entropy solutions is a modification of the weak formulation of the entropy inequality in which the right hand side is allowed to be an error of the form
\[
-\int_0^T \left(
	\int_{\setR^n}\abs{\phi(t,x)}\d\mu_{0,t}(x)
	+ \int_{\setR^n}\abs{\nabla_x\phi(t,x)}\d\mu_{1,t}(x)
	\right)\dt
\]
when tested against $\phi\in C^\infty_c\bigl((0,T)\times\setR;[0,\infty)\bigr)$ for some non-negative measures $\mu_{0,t},\mu_{1,t}$, instead of being zero.

The structure of the stability theorem that we obtain can be synthesized in the following manner. This statement is of course imprecise and we refer the reader to \autoref{sec:stability} for the correct definitions and formulations.

\begin{theorem}[Informal version of \autoref{thm:stability}]\label{thm:informal-stability}
Given two fluxes $P,Q$ satisfying suitable regularity assumptions, let $u$ and $v$ be a $\mu$-quasi-entropy solution and $\nu$-quasi-entropy solution of the conservation laws
\begin{align*}
\partial_t u + \div_x\bigl(P(t,x,u)\bigr) &= 0 & &\text{and} &
\partial_t v + \div_x\bigl(Q(t,x,v)\bigr) &= 0
\end{align*}
respectively,
with the bounds
\begin{align*}
u(t,\plchldr),v(t,\plchldr) &\leq R(t), &
\norm{u(t,\plchldr)}_{BV},\norm{v(t,\plchldr)}_{BV} & \leq B(t),
\end{align*}
for some increasing functions $R,B:[0,T)\to[1,\infty)$.

Then
\[
\begin{split}
\int_{\setR^n} \abs{u(t,x)-v(t,x)}\dx
&\leq \int_{\setR^n} \abs{u(0,x)-v(0,x)}\dx
	+ \int_0^t C_{P,Q}(s) \int_{\setR^n} \abs{u(s,x)-v(s,x)}\dx\ds
	\spliteq
	+ B(t) \int_0^t \norm{P(s)-Q(s)}_* \ds
	+ M_0 + c_nM_1
	+ B(t)C_{P,Q}(t)\sqrt{M_1}
\end{split}
\]
where
\begin{align*}
M_0 &= \int_0^t (\mu_{0,s}+\nu_{0,s})(\setR^n) \ds, &
M_1 &= \int_0^t (\mu_{1,s}+\nu_{1,s})(\setR^n) \ds,
\end{align*}
\[
\norm{P(s)-Q(s)}_*
= \norm{\div_x(P-Q)(s,x,u)}_{L^1_xL^\infty_u}
	+ \norm{\partial_u(P-Q)(s,x,u)}_{L^\infty_xL^\infty_u},
\]
and $C_{P,Q}$ is a function of the norms of some derivatives of $P$ and $Q$ which is bounded under the assumptions.
\end{theorem}

The major simplification we introduced in this informal statement is that we integrate over the whole space, whereas in \autoref{thm:stability} the estimate is localized with a suitable weight function.

The need for this article is motivated by some limitations in the available literature. The possibility to work at the same time with distinct fluxes depending on $(t,x,u)$ and quasi-entropy solutions has several primary benefits.


Firstly, it allows to study problems with fluxes $P[u](t,x,u)$ which depend non-locally on the solution $u$ itself, for which it is unavoidable to require the space-time dependence and to consider two distinct fluxes. In this context, the term $\norm{P-Q}_*$ sometimes can be estimated with $\norm{u-v}_{L^1_x}$ allowing to close a Grönwall type inequality (\autoref{prop:gronwall}). Moreover, in the situations where there is an estimate of the form $\norm{P-Q}_*\leq\norm{u-v}_1$, we show in \autoref{sec:existence} and \autoref{sec:uniqueness} the conditional existence and uniqueness of entropy solutions respectively.


Secondly, quasi-entropy solutions arise naturally in the study of numerical methods and the stability theorem can be used to derive their rates of convergence (\autoref{sec:rates}).
In particular, in \autoref{sec:particles} we obtain for the first time the rate of convergence of a deterministic particle method presented in \cite{RadiciStra} (see also \cite{FagioliTse}) to solve a scalar conservation law in one dimension inspired by a model of traffic with congestion (\cite{DiFrancescoRosini,DiFrancescoFagioliRadici,DiFrancescoStivaletta}). For this application we really need both improvements (two space-time fluxes, quasi-entropy solutions) with respect to the present literature: indeed, the stability result of \cite{KarlsenRisebro} only treats entropy solutions whereas the discrete approximations produced by the particle method are only shown to be quasi-entropy solutions (\autoref{prop:particles-quasi-entropy}). The approximation error with $N$ particles that we obtain is of order $N^{-1/2}$, which is shown to be sharp in \autoref{sec:particles}. The adopted technique also provides an independent way to prove the existence of entropy-solutions which bypasses the compactness argument used in \cite{RadiciStra}.

In addition we recover the known rates of convergence of the vanishing viscosity method (\autoref{sec:viscosity}), extending its validity to the case of a flux $P(t,x,u)$, and of the front tracking method (\autoref{sec:front-tracking}).

%

\paragraph{Strategy of the proof of the stability theorem.}
Our approach is inspired by \cite{Kruzkov,KarlsenRisebro}, with some modifications.

From \cite{Kruzkov} we adopt the general framework of doubling of variables.
However, dealing with two distinct fluxes has some important consequences: some algebraic symmetries exploited in \cite{Kruzkov} do not hold anymore and this causes the appearance of additional mixed terms which have to be estimated, for instance
\[
\div_x\bigl[\sign(u-v)\bigl(P(t,x,u)-P(t,x,v)-Q(t,y,u)+Q(t,y,v)\bigr)\bigr]
\]
where $u=u(t,x)$ and $v=v(t,y)$.

The way we deal with this increased complexity shares a closer resemblance with the scheme presented in \cite[Theorem 1.3]{KarlsenRisebro}, which also treats two fluxes, albeit of product form; in particular we decompose the entropy inequality into similar terms and to estimate the one described above we apply the chain rule.
In order to do this, we need to require $u$ and $v$ to be $BV$ instead of merely $L^\infty$ and at the beginning of the proof we also need to regularize the two fluxes and the absolute value and sign functions.

The fact that we consider quasi-entropy solutions instead of exact entropy solutions excludes the possibility of performing the full dedoubling of the space variables $\abs{x-y}\to0$.
As in \cite{Kuznetsov,BouchutPerthame1998}, we collapse the space variables $\abs{x-y}\lesssim\beta$ at an optimal scale $\beta$ which is determined by the balance between the mass of the error terms and the modulus of continuity in $L^1$ of the translations of the solutions.


\subsection{Future perspectives}

\paragraph{More general definition of quasi-entropy solutions.}
In \autoref{def:mu-quasi-entropy-solution} of quasi-entropy solutions we impose that the right hand side is estimated with measures whose disintegration with respect to time is of the form $\d\mu(x)\dt$. A natural generalization is to extend it to arbitrary measures $\d\mu(t,x)$ whose projection on time is not necessarily absolutely continuous with respect to the Lebesgue measure.
This generalization goes in the direction of \cite{BouchutPerthame1998}.
The usefulness of this extensions comes from the possibility of studying the convergence rate of numerical schemes which are discrete/discontinuous processes in time, for instance layering/smoothing, finite difference/volume/elements such as Godunov, or higher order ones.
With this more flexible notion one can no longer expect quasi-solutions to be in $C\bigl([0,T);L^1_\loc(\setR^n)\bigr)$, so this makes it more difficult to perform the dedoubling of the time variables.

\paragraph{Conservation laws with diffusion and source terms.}
An interesting and useful generalization of the present result that we intend to pursue consists in considering conservation laws of the form
\[
\partial_t u + \div_x\bigl(P(t,x,u)\bigr) = \Delta A(u) + f(t,x,u).
\]

This research direction would improve upon the work of \cite{VolpertParabolic}, which treats a single flux and entropy solutions, and \cite{KarlsenRisebro}, which treats two time-independent fluxes of product type and entropy solutions.
Merging our approach with theirs we plan to extend the stability theorem to two general fluxes and $BV$ quasi-entropy solutions.

This is for example relevant for determining the rate of convergence of deterministic particle schemes introduced in \cite{FagioliRadiciDiffusion,DaneriRadiciRunaLinear,DaneriRadiciRunaNonlinear} for solving non-local conservation laws where the flux is a convolution with the solutions itself and there is a linear/nonlinear mobility to model congestion.

\paragraph{Non-conditional existence for non-local problems.}
In \autoref{sec:existence} we prove the conditional existence of entropy solutions provided that one can construct an approximating sequence of quasi-entropy solutions. It is natural to ask whether, under some general assumptions on the fluxes, such a sequence can be constructed.

This approach could be especially beneficial to show the existence of entropy solutions to non-local problems, for which to the best of our knowledge the results are rather sparse and specific to some particular equations (for instance traffic, pedestrian, chemotaxis models).

\subsection*{Acknowledgments}

The three authors wish to thank the EPFL, where part of the work has been carried out while being affiliated to the Mathematics department there, and in particular the SNF grant 182565.

In addition, E.M.\ acknowledges the support received from the European Union's Horizon 2020 research and innovation program under the Marie Sk{\l}odowska-Curie grant 101025032.

\section{Stability}\label{sec:stability}

The aim of this section is to state and prove a precise formulation of the stability result that was informally presented in \autoref{thm:informal-stability}. We begin by establishing the regularity assumptions required for the fluxes.

\begin{assumptions}[Regularity of the flux]\label{as:flux}
We require $P:[0,T)\times\setR^n\times[0,\infty)\to\setR^n$ to be a flux satisfying the following conditions:
\begin{enumerate}[({A}1)]
\item $(t,x)\mapsto P(t,x,0) \in L^1_\loc\bigl([0,T)\times\setR^n\bigr)$;
\item $(t,x)\mapsto \div_xP(t,x,0) \in L^1_\loc\bigl([0,T)\times\setR^n\bigr)$;
\item $u\mapsto P(t,x,u) \in \Lip_\loc(\setR)$ locally uniformly for $x\in\setR^n$ with dependence $L^1_\loc\bigl([0,T)\bigr)$ in time, i.e. for every $R>0$ there is a function $C_R\in L^1_\loc\bigl([0,T);[0,\infty)\bigr)$ such that
\[
\abs{P(t,x,u)-P(t,x,v)} \leq C_R(t)\abs{u-v} \qquad \forall u,v\in[0,R],\ \forall(t,x)\in[0,T)\times B_R(0);
\]
\item $u\mapsto \div_xP(t,x,u) \in \Lip_\loc(\setR)$ locally uniformly for $x\in\setR^n$ with dependence $L^1_\loc\bigl([0,T)\bigr)$ in time, i.e. for every $R>0$ there is a function $C_R\in L^1_\loc\bigl([0,T);[0,\infty)\bigr)$ such that
\[
\abs{\div_xP(t,x,u)-\div_xP(t,x,v)} \leq C_R(t)\abs{u-v} \qquad \forall u,v\in[0,R],\ \forall(t,x)\in[0,T)\times B_R(0);
\]
\item $(t,x)\mapsto\div_x P(t,x,u)$ and $(t,x)\mapsto\partial_u P(t,x,u) \in L^1_\loc\bigl([0,T);W^{1,\infty}_\loc(\setR^n)\bigr)$ locally uniformly for $u\in[0,\infty)$, i.e.\ in particular for every $R>0$
\[
\norm{\nabla_x\div_x P(t,x,u)}_{L^\infty\bigl(B_R(0)\times[0,R]\bigr)},\ 
\norm{\nabla_x\partial_u P(t,x,u)}_{L^\infty\bigl(B_R(0)\times[0,R]\bigr)}
\in L^1_\loc\bigl([0,T)\bigr).
\]
\end{enumerate}
\end{assumptions}

Notice in particular that under the previous assumptions $P$ turns out to be a Carathéodory function, therefore $P(t,x,u(t,x))$ is a well defined measurable function of $(t,x)$ when $u(t,x)$ which is itself measurable. Moreover, $\abs{P(t,x,u(t,x))} \leq \abs{P(t,x,0)}+C_R(t)\abs{u(t,x)}$ for some function $C_R(t)>0$ locally in $x$, hence $P(t,x,u(t,x))$ is $L^1_\loc\bigl([0,T)\times\setR^n\bigr)$ as soon as so is $u$. The same considerations apply to $(\div_x P)(t,x,u(t,x))$. This observation ensures that the formulations \eqref{eq:def-entropy-ineq-hyp} and \eqref{eq:def-quasi-entropy-ineq-hyp} of the following definitions make sense.

The assumption (A5) is unnecessary to ensure the meaningfulness of the entropy inequality and is required only to estimate the error arising in the stability of quasi-entropy solutions. When working with entropy solutions (for instance in \autoref{rmk:entropic-statement}, \autoref{thm:uniqueness-nonlocal-compact} and \autoref{thm:uniqueness-nonlocal}) this assumption can be omitted.

For the reader's convenience, we recall the classical definition of \emph{entropy solution} in the sense of \cite{Kruzkov}.

\begin{definition}[Entropy solution]\label{def:entropy-solution}
Let $P:[0,T)\times\setR^n\times[0,\infty)\to\setR^n$ a flux satisfying (A1)\nobreakdash--(A4) of \autoref{as:flux}.
We say that a non-negative function $u\in C\bigl([0,T);L^1_\loc(\setR^n)\bigr) \cap L^\infty_\loc\bigl([0,T);BV_\loc(\setR^n)\bigr)$ is an entropy solution of the scalar conservation law
\[
\partial_t u + \div_x(P(t,x,u)) = 0
\]
if the following entropy inequality
\begin{equation}\label{eq:def-entropy-ineq-hyp}
\begin{split}
\Intd & \bigl\{
	\abs{u-c}\partial_t\phi + \sign(u-c)\bigl[
		\bigl(P(t,x,u)-P(t,x,c)\bigr)\cdot\nabla_x\phi - \div_xP(t,x,c)\phi
	\bigr] \bigr\} \dx \dt
\geq 0
\end{split}
\end{equation}
holds for every constant $c\in[0,\infty)$ and non-negative test function $\phi\in C^\infty_c\bigl((0,T)\times\setR;[0,\infty)\bigr)$.\footnote{\label{nt:c1}Observe that it is equivalent to require the inequality for all $\phi\in C^1_c\bigl((0,T)\times\setR;[0,\infty)\bigr)$.}
\end{definition}

As anticipated in the introduction, the purpose of the article is to extend the stability beyond entropy solutions. We therefore introduce a notion of \emph{quasi-entropy solution} that will be suitable for our needs.

\begin{definition}[$(\mu_0,\mu_1)$-quasi-entropy solution]\label{def:mu-quasi-entropy-solution}
Let $\mu_{0,t},\mu_{1,t}\in L^1_\loc\bigl([0,T);\Meas_+(\setR^n)\bigr)$ be locally-finite non-negative Borel measures and $P:[0,T)\times\setR^n\times[0,\infty)\to\setR$ a flux satisfying (A1)\nobreakdash--(A4) of \autoref{as:flux}.
We say that a non-negative function $u\in C\bigl([0,T);L^1_\loc(\setR^n)\bigr) \cap L^\infty_\loc\bigl([0,T);BV_\loc(\setR^n)\bigr)$ is a $(\mu_0,\mu_1)$-quasi-entropy solutions of the scalar conservation law
\[
\partial_t u + \div_x(P(t,x,u)) = 0
\]
if the following entropy inequality
\begin{equation}\label{eq:def-quasi-entropy-ineq-hyp}
\begin{split}
\Intd & \bigl\{
	\abs{u-c}\partial_t\phi + \sign(u-c)\bigl[
		\bigl(P(t,x,u)-P(t,x,c)\bigr)\cdot\nabla_x\phi - \div_xP(t,x,c)\phi
	\bigr] \bigr\} \dx \dt \\
&\geq -	\Intd \abs{\phi(t,x)} \d\mu_{0,t}(x)\dt
	- \Intd \abs{\nabla_x\phi(t,x)} \d\mu_{1,t}(x)\dt
\end{split}
\end{equation}
holds for every constant $c\in[0,\infty)$ and non-negative test function $\phi\in C^\infty_c\bigl((0,T)\times\setR^n;[0,\infty)\bigr)$.\footref{nt:c1}
\end{definition}

Notice that in comparison to \autoref{def:entropy-solution}, the weaker notion of \autoref{def:mu-quasi-entropy-solution} allows for a controlled violation of \eqref{eq:def-entropy-ineq-hyp}. The measures that control this right hand side will then play a crucial role in the stability estimates.

%
%
%

For the sake of keeping the terms appearing in both the statement of the main theorem and its proof shorter and more readable, we introduce some notation.

First of all, to avoid possible confusion when differentiating composite functions, we denote the partial derivatives of a function $F:[0,T)\times\setR^n\times[0,\infty)\to\setR$ with respect to its three arguments as $\partial_1F$, $\nabla_2F$ and $\partial_3F$ respectively. Analogously, this convention extends to $\div_2F$ as well.

Moreover, because of the local nature of \autoref{thm:stability} that we are going to state, we will need to refer to various norms of the fluxes computed locally, as opposed to the whole $[0,T)\times\setR^n\times[0,\infty)$.
In the statement there will be two functions responsible for the localization: a time dependent weight $\Theta$ localizing in space and a function $R$ bounding the density on the support of $\Theta$.
The notation that we introduce is implicit with respect to $\Theta$ and $R$, and has to be understood in the context where these functions are fixed.

\begin{notation}\label{rmk:notation-norms}
Given a functions $R:[0,T)\to[0,\infty)$ and a function $\Theta\in C^1\bigl([0,T)\times\setR^n;[0,\infty)\bigr)$ compactly supported in space for every time, we introduce the following notation to estimate some norms of functions $F(t,x,u)$ only in the domain $(x,u)\in\bigl(\supp\Theta(t,\plchldr)\bigr)_1\times[0,R(t)]$, where with $E_r$ we denote the enlargement of radius $r$ in space of the set $E$:
\begin{gather*}
\Lip_2(F(t)) = \sup\set*{\frac{\abs{F(t,x_2,u)-F(t,x_1,u)}}{\abs{x_2-x_1}}}
	{x_1,x_2\in\bigl(\supp\Theta(t,\plchldr)\bigr)_1,\ u\in[0,R(t)]}, \\
\Lip_3(F(t)) = \sup\set*{\frac{\abs{F(t,x,u_2)-F(t,x,u_1)}}{\abs{u_2-u_1}}}
	{x\in\bigl(\supp\Theta(t,\plchldr)\bigr)_1,\ u_1,u_2\in[0,R(t)]}.
\end{gather*}
\end{notation}

\begin{remark}
With the notation introduced above, observe that (A1)\nobreakdash--(A4) of \autoref{as:flux} ensures that $\Lip_3(P(t))$ and $\Lip_3(\div_2 P(t))$ are functions belonging to $L^1_\loc([0,T))$, whereas (A5) ensures that $\Lip_2(\div_2 P(t))$ and $\Lip_2(\partial_3P(t))$ are $L^1_\loc([0,T))$ too.
\end{remark}

\subsection{Main results}

In this section we present our main results, which deal with non-negative quasi-entropy solutions belonging to the space
\begin{equation}\label{eq:cand-sol}
\mathscr{S}\bigl([0,T),\setR^n\bigr)
= C\bigl([0,T); L^1_\loc(\setR^n)\bigr)
	\cap L^\infty_\loc\bigl([0,T); L^\infty_\loc(\setR^n)\bigr)
	\cap L^\infty_\loc\bigl([0,T); BV_\loc(\setR^n)\bigr).
\end{equation}

\begin{theorem}[Stability]\label{thm:stability}
Let $P,Q:[0,T)\times\setR^n\times[0,\infty)\to\setR^n$ be two fluxes satisfying \autoref{as:flux}.
Let $u,v\in\mathscr{S}\bigl([0,T),\setR^n\bigr)$ be non-negative quasi-entropy solutions of
\begin{align*}
\partial_t u + \div_x(P(t,x,u)) &= 0 & &\text{and} &
\partial_t v + \div_x(Q(t,x,v)) &= 0
\end{align*}
in the sense of \autoref{def:mu-quasi-entropy-solution}, i.e.\ there are measures $\mu_{0,t},\mu_{1,t},\nu_{0,t},\nu_{1,t}\in L^1_\loc\bigl([0,T);\Meas_+(\setR^n)\bigr)$ such that the quasi-entropy inequalities
\begin{subequations}\label{eq:entropy-ineq-hyp}
\begin{equation}
\begin{split}
\Intd & \bigl\{
	\abs{u-c}\partial_t\phi + \sign(u-c)\bigl[
		\bigl(P(t,x,u)-P(t,x,c)\bigr)\cdot\nabla_x\phi - \div_xP(t,x,c)\phi
	\bigr] \bigr\} \dx \dt \\
&\geq -	\Intd \abs{\phi(t,x)} \d\mu_{0,t}(x)\dt - \Intd \abs{\nabla_x\phi(t,x)} \d\mu_{1,t}(x)\dt,
\end{split}
\end{equation}
\vspace{-.25cm}
\begin{equation}
\begin{split}
\Intd & \bigl\{
	\abs{v-c}\partial_t\phi + \sign(v-c)\bigl[
		\bigl(Q(t,x,v)-Q(t,x,c)\bigr)\cdot\nabla_x\phi - \div_xQ(t,x,c)\phi
	\bigr] \bigr\} \dx \dt \\
&\geq -	\Intd \abs{\phi(t,x)} \d\nu_{0,t}(x)\dt - \Intd \abs{\nabla_x\phi(t,x)} \d\nu_{1,t}(x)\dt
\end{split}
\end{equation}
\end{subequations}
hold for every constant $c\in[0,\infty)$ and non-negative test function $\phi\in C^1_c\bigl((0,T)\times\setR;[0,\infty)\bigr)$.

Let $\Theta\in C^1\bigl([0,T)\times\setR^n;[0,\infty)\bigr)$ be a fixed weight function compactly supported in space for every time and let $R,B:[0,T)\to[0,\infty)$ be two increasing functions satisfying the following properties: the estimates
\begin{align*}
\norm{u(t,\plchldr)}_{L^\infty(\Omega_t)} &\leq R(t),
	& \TV_{\Omega_t}\bigl(u(t,\plchldr)\bigr) &\leq B(t), \\
\norm{v(t,\plchldr)}_{L^\infty(\Omega_t)} &\leq R(t),
	& \TV_{\Omega_t}\bigl(v(t,\plchldr)\bigr) &\leq B(t),
\end{align*}
hold for every $t\in[0,T)$ where $\Omega_t=\bigl(\supp\Theta(t,\plchldr)\bigr)_1\subset\setR^n$, and 
\begin{equation}\label{eq:Theta-slope}
\partial_t\Theta(t,x)
\leq-\Lip_3\bigl(Q(t)\bigr)\abs{\nabla_x\Theta(t,x)},
\qquad \forall (t,x)\in[0,T)\times\setR^n.
\end{equation}

%

Then for every $0\leq t_1<t_2<T$ the following inequality holds:
\begin{equation}\label{eq:MRS}
\begin{split}
&\left[\int_{\setR^n} \abs{u(t,x)-v(t,x)} \Theta(t,x) \dx \right]_{t_1}^{t_2} \\
&\leq \int_{t_1}^{t_2} [3\Lip_3(\div_2P(t))+\Lip_3(\div_2Q(t))]
		\int_{\setR^n} \abs{u(t,x)-v(t,x)}\Theta(t,x) \dx\dt
	\spliteq
	+\intd \norm{\div_2(P-Q)(t,x,\plchldr)}_{L^\infty([0,R(t)])}
		\Theta(t,x) \dx\dt
	\spliteq
	+2B(t_2) \int_{t_1}^{t_2}
		\norm{\Lip_3\bigl((P-Q)(t,\plchldr)\bigr) \Theta(t,\plchldr)}_\infty \dt
	\spliteq
	+\norm{\Theta}_\infty \int_{t_1}^{t_2}
		(\mu_{0,t}+\nu_{0,t})(\Omega_t) \dt
	+\frac12\norm{\nabla_x\Theta}_\infty \int_{t_1}^{t_2}
		(\mu_{1,t}+\nu_{1,t})(\Omega_t) \dt
	\spliteq
	+C(t_1,t_2) \min\oleft\{ M(t_1,t_2)^{1/2}, 1 \right\}
		+ c_n \norm{\Theta}_\infty
		\max\oleft\{ M(t_1,t_2)^{1/2}, M(t_1,t_2) \right\},
\end{split}
\end{equation}
where
\begin{subequations}\label{eq:M-C}
\begin{equation}\label{eq:M}
M(t_1,t_2) = \int_{t_1}^{t_2} (\mu_{1,t}+\nu_{1,t})(\Omega_t) \dt,
\end{equation}
\begin{equation}\label{eq:C}
\begin{split}
C(t_1,t_2) &=
	\norm{\Theta}_\infty B(t_2) \left( 2 + \int_{t_1}^{t_2}
		[3\Lip_3(\div_2P(t))+\Lip_3(\div_2Q(t))
			+2\Lip_2(\partial_3Q(t))] \dt \right)
	\spliteq
	+\norm{\Theta}_{L^\infty_t L^1_x} \int_{t_1}^{t_2} \Lip_2(\div_2P(t)) \dt
	+\norm{\Theta}_{L^\infty_t\Lip_x} B(t_2)
		\int_{t_1}^{t_2} \Lip_3\bigl((P-Q)(t)\bigr) \dt
	\spliteq
	+\frac12\norm{\Theta}_{L^\infty_t\Lip_x}
		\int_{t_1}^{t_2} \!\!\!\! \int_{\Omega_t}
			\norm{\div_2(P-Q)(t,x,\plchldr)}_{L^\infty([0,R(t)])} \dx\dt,
\end{split}
\end{equation}
\end{subequations}
and $c_n>0$ is a dimensional constant.
\end{theorem}

\begin{remark}[Compact solutions]\label{rmk:compact-support}
When the functions $u$ and $v$ are compactly supported in space uniformly in time, we can take any function $\Theta\in C^1\bigl([0,\infty)\times\setR^n;[0,\infty)\bigr)$ such that $\Theta=1$ in $\bigl(\supp(u)\cup\supp(v)\bigr)_1$, without requiring \eqref{eq:Theta-slope}. Indeed, inside the proof of \autoref{thm:stability} the slope condition \eqref{eq:Theta-slope} is used only in the last inequality of \eqref{eq:time-survivors}, which remains valid for such a $\Theta$ because $[\partial_1\Theta+\Lip_3(Q^\alpha(t))\abs{\nabla_2\Theta}]\eta_\eps(u-v)\omega_\beta$ is identically zero.
\end{remark}

\begin{remark}[Entropy solutions]\label{rmk:entropic-statement}
The result of \autoref{thm:stability} can be simplified significantly under the assumption that the functions $u$ and $v$ are exact entropy solutions of the scalar conservations laws instead of merely quasi-entropy solutions. Indeed in such a case all four measures $\mu_{0,t},\mu_{1,t},\nu_{0,t},\nu_{1,t}$ vanish for every $t$, and the analogue of \eqref{eq:MRS} that we obtain is
\[
\begin{split}
&\left[\int_{\setR^n} \abs{u(t,x)-v(t,x)} \Theta(t,x) \dx \right]_{t_1}^{t_2} \\
&\leq \int_{t_1}^{t_2} [3\Lip_3(\div_2P(t))+\Lip_3(\div_2Q(t))]
		\int_{\setR^n} \abs{u(t,x)-v(t,x)}\Theta(t,x) \dx\dt
	\spliteq
	+\intd \norm{\div_2(P-Q)(t,x,\plchldr)}_{L^\infty([0,R(t)])}
		\Theta(t,x) \dx\dt
	\spliteq
	+B(t_2) \int_{t_1}^{t_2}
		\norm{\Lip_3\bigl((P-Q)(t,\plchldr)\bigr) \Theta(t,\plchldr)}_\infty \dt.
\end{split}
\]
Interestingly, in this case we can omit the assumption (A5) of \autoref{as:flux}. Please refer to \autoref{rmk:entropic-statement-proof} for details on how to obtain this modified statement.
This result extends \cite[Theorem~1.3]{KarlsenRisebro}, which is limited to fluxes of product form $P(t,x,u)=f(u)k(x)$.


If in addition $P=Q$, then the estimate is independent of the total variation bound $B(t)$. This is consistent with the result by \cite{Kruzkov} which holds for $L^\infty$ entropy solutions.
\end{remark}

When the terms of the form $P-Q$ in \eqref{eq:MRS} can be estimated by some integral of $\abs{u-v}$, the whole inequality assumes a form suitable for the application of Grönwall theorem. More precisely, we can state the following proposition.

\begin{proposition}[Grönwall estimate]\label{prop:gronwall}
Let $u,v$, $P,Q$, $\mu_{0,t},\mu_{1,t},\nu_{0,t},\nu_{1,t}$, $\Theta,\Omega_t$, $R,B$ as in \autoref{thm:stability}. Assume in addition that the fluxes $P$ and $Q$ are close together in the following sense: there is a non-negative function $h\in L^1_\loc([0,T))$ such that for a.e.\ $t\in[0,T)$
\begin{subequations}\label{eq:gronwall-assumptions}
\begin{equation}
\int_{\setR^n} \norm{\div_2(P-Q)(t,x,\plchldr)}_{L^\infty([0,R(t)])} \Theta(t,x) \dx
\leq h(t)\int_{\setR^n} \abs{u(t,x)-v(t,x)} \Theta(t,x) \dx,
\end{equation}
\begin{equation}
\norm{\Lip_3\bigl((P-Q)(t,\plchldr)\bigr) \Theta(t,\plchldr)}_{L^\infty(\setR^n)}
\leq h(t) \int_{\setR^n} \abs{u(t,x)-v(t,x)} \Theta(t,x) \dx.
\end{equation}
\end{subequations}
Then for $0\leq t_1<t_2<T$ we have the estimate
\begin{equation}\label{eq:gronwall}
\begin{split}
\int_{\setR^n} & \abs{u(t_2,x)-v(t_2,x)} \Theta(t_2,x) \dx \\
&\leq \left( \int_{\setR^n} \abs{u(t_1,x)-v(t_1,x)} \Theta(t_1,x) \dx
	+ \Phi(t_1,t_2) \right) \exp\oleft(\int_{t_1}^{t_2} f(t)\dt\right),
\end{split}
\end{equation}
where
\begin{align*}
f(t) &= [3\Lip_3(\div_2P(t))+\Lip_3(\div_2Q(t))] + [1+2B(T)]h(t), \\
\Phi(t_1,t_2) &= \norm{\Theta}_\infty \int_{t_1}^{t_2}
		(\mu_{0,t}+\nu_{0,t})(\Omega_t) \dt
	+\frac12\norm{\nabla_x\Theta}_\infty \int_{t_1}^{t_2}
		(\mu_{1,t}+\nu_{1,t})(\Omega_t) \dt
	\spliteq
	+C(t_1,t_2) \min\left\{ M(t_1,t_2)^{1/2}, 1 \right\}
		+c_n \norm{\Theta}_\infty
		\max\left\{ M(t_1,t_2)^{1/2}, M(t_1,t_2) \right\},
\end{align*}
and $M(t_1,t_2)$ and $C(t_1,t_2)$ are the same as in \autoref{thm:stability}.
\end{proposition}

This proposition will be used in \autoref{sec:existence} and \autoref{sec:uniqueness} to show the conditional existence and the uniqueness of entropy solutions for the conservation law with non-local flux $P[u](t,x,u)$. Moreover, in \autoref{sec:particles} we show an application to a conservation law where the assumptions \eqref{eq:gronwall-assumptions} are satisfied.

\subsection{Proofs}

\begin{proof}[Proof of \autoref{thm:stability}]
For convenience of the reader we split the proof in several distinct steps indicated by paragraphs.

We prove the statement for $t_1>0$. The case $t_1=0$ is recovered by continuity.

\paragraph{Regularization of the fluxes.}
The first step consists in regularizing by convolution the fluxes with respect to the space and density variables. This forces us to consider the fluxes evaluated at negative densities too. Therefore, for convenience we extend the fluxes to $P,Q:[0,T)\times\setR^n\times\setR\to\setR^n$ by setting $P(t,x,u)=P(t,x,-u)$ and $Q(t,x,u)=Q(t,x,-u)$ for $u<0$. The extended fluxes enjoy the same regularity assumptions as the original ones.

Given $\alpha\in(0,1]$, we define the regularized fluxes $P^\alpha,Q^\alpha:[0,T)\times\setR^n\times\setR\to\setR^n$ obtained by convolution in space and with respect to the density
\begin{align*}
P^\alpha(t,x,u) = \int_{\setR^n}\int_\setR P(t,x',u')\omega_\alpha(x-x')\rho_\alpha(u-u')\d u' \dx' , \\
Q^\alpha(t,x,u) = \int_{\setR^n}\int_\setR Q(t,x',u')\omega_\alpha(x-x')\rho_\alpha(u-u')\d u' \dx',
\end{align*}
where $\omega$ and $\rho$ are symmetric $C^\infty_c$ mollifiers in $\setR^n$ and $\setR$ respectively supported in the corresponding unit balls, and $\omega_\alpha(x)=\alpha^{-n}\omega(\alpha^{-1}x)$ and $\rho_\alpha(u)=\alpha^{-1}\rho(\alpha^{-1}u)$ are the rescalings that preserve the $L^1$ norm. For every $t$ and almost every $x$, the regularized fluxes enjoy the regularity estimates
\begin{align*}
\norm{P^\alpha(t,x,\plchldr)-P(t,x,\plchldr)}_{L^\infty([0,R(t)])}
	&\leq \alpha \norm{\partial_3P(t,\plchldr,\plchldr)}_{L^\infty(B(x,\alpha)\times[-\alpha,R(t)+\alpha])} , \\
\norm{\div_2P^\alpha(t,x,\plchldr)-\div_2P(t,x,\plchldr)}_{L^\infty([0,R(t)])}
	&\leq \alpha \norm{\partial_3\div_2P(t,\plchldr,\plchldr)}_{L^\infty(B(x,\alpha)\times[-\alpha,R(t)+\alpha])} .
\end{align*}
Adding and subtracting in the quasi-entropy inequalities \eqref{eq:entropy-ineq-hyp} the corresponding regularized terms with $P^\alpha$ and $Q^\alpha$, we obtain
\begin{subequations}\label{eq:entropy-ineq-reg}
\begin{equation}\label{eq:entropy-ineq-reg-u}
\begin{split}
\Intd & \bigl\{
	\abs{u-c}\partial_t\phi + \sign(u-c)\bigl[
		\bigl(P^\alpha(t,x,u)-P^\alpha(t,x,c)\bigr)\cdot\nabla_x\phi
		- \div_xP^\alpha(t,x,c)\phi \bigr] \bigr\} \dx \dt \\
	&\geq -	\Intd \abs{\phi(t,x)} \d\mu_{0,t}(x)\dt
	- \Intd \abs{\nabla_x\phi(t,x)} \d\mu_{1,t}(x)\dt
	-\alpha e(P,\phi),
\end{split}
\end{equation}
\vspace{-.75cm}
\begin{equation}\label{eq:entropy-ineq-reg-v}
\begin{split}
\Intd & \bigl\{
	\abs{v-c}\partial_t\phi + \sign(v-c)\bigl[
		\bigl(Q^\alpha(t,x,v)-Q^\alpha(t,x,c)\bigr)\cdot\nabla_x\phi
		- \div_xQ^\alpha(t,x,c)\phi \bigr] \bigr\} \dx \dt \\
	&\geq -	\Intd \abs{\phi(t,x)} \d\nu_{0,t}(x)\dt
	- \Intd \abs{\nabla_x\phi(t,x)} \d\nu_{1,t}(x)\dt
	-\alpha e(Q,\phi) ,
\end{split}
\end{equation}
\end{subequations}
where $e(P,\phi)$ is defined as the right hand side of the following inequality used to estimate the error introduced by the convolution:
\begin{equation}\label{eq:e-P-phi}
\begin{split}
& \Intd \biggl\{
	2\norm{\partial_3P(t,\plchldr,\plchldr)}_{L^\infty(B(x,\alpha)\times[-\alpha,R(t)+\alpha])} \abs{\nabla_x\phi}
	\spliteq\qquad\qquad
	+\norm{\partial_3\div_2P(t,\plchldr,\plchldr)}_{L^\infty(B(x,\alpha)\times[-\alpha,R(t)+\alpha])} \abs{\phi}
	\biggr\} \dx\dt \\
&\leq \int_0^T \left\{
	2\Lip_3\bigl(P(t)\bigr) \int_{\setR^n} \abs{\nabla_x\phi} \dx
	+\Lip_3\bigl(\div_2P(t)\bigr) \int_{\setR^n} \abs{\phi} \dx \right\} \dt
\eqqcolon e(P,\phi).
\end{split}
\end{equation}

\paragraph{Regularization of the absolute value.}\label{par:reg-abs}
We now introduce a second convolution in order to regularize the absolute value function. Given the mollifier $\rho_\eps$ we define the regularized absolute value
\[
\eta_\eps(u) = \int_\setR \abs{s} \rho_\eps(u-s) \d s
	-\int_\setR \abs{s} \rho_\eps(s) \d s
\]
and its translation $\eta_{\eps,v}(u)=\eta_\eps(u-v)$.

The goal of this section is to show that convolving \eqref{eq:entropy-ineq-reg-u} with $\frac12\eta_{\eps,v}''(c)=\rho_\eps(c-v)$ has the effect of replacing the functions $\eta_{0,c}(u)=\abs{u-c}$, $\sign(u-c)\bigl(P^\alpha(t,x,u)-P^\alpha(t,x,c)\bigr)$ and $\sign(u-c)\div_xP^\alpha(t,x,c)$ with $\eta_{\eps,v}(u)$,
\[
P^\alpha_\eps(t,x,u,v) = \int_v^u \partial_3P^\alpha(t,x,s)\eta_{\eps,v}'(s) \d s,
\quad\text{and}\quad
\int_v^u \div_xP^\alpha(t,x,s)\eta_{\eps,v}''(s)\d s
\]
respectively. Similarly, convolving \eqref{eq:entropy-ineq-reg-v} with $\frac12\eta_{\eps,u}''(c)=\rho_\eps(c-u)$ has the effect of replacing the functions $\eta_{0,c}(v)=\abs{v-c}$, $\sign(v-c)\bigl(Q^\alpha(t,x,v)-Q^\alpha(t,x,c)\bigr)$ and $\sign(v-c)\div_xQ^\alpha(t,x,c)$ with $\eta_{\eps,u}(v)$,
\[
Q^\alpha_\eps(t,x,u,v) = \int_u^v \partial_3Q^\alpha(t,x,s)\eta_{\eps,u}'(s) \d s,
\quad\text{and}\quad
\int_u^v \div_xQ^\alpha(t,x,s)\eta_{\eps,u}''(s)\d s
\]
respectively. Notice that in the definitions of $P^\alpha_\eps$ and $Q^\alpha_\eps$ the role of $u$ and $v$ is interchanged, although we write them in the same order as arguments.

With regard to \eqref{eq:entropy-ineq-reg-u}, with the mentioned convolution we obtain
\[
\begin{split}
\int_\setR \abs{u-c} \rho_\eps(c-v) \d c
= \int_\setR \abs{u-c} \rho_\eps\bigl((u-v)-(u-c)\bigr) \d c
= \eta_\eps(u-v) + \int_\setR\abs{s}\rho_\eps(s)\ds,
\end{split}
\]
moreover
\[
\begin{split}
&\int_\setR \sign(u-c) \bigl(P^\alpha(t,x,u)-P^\alpha(t,x,c)\bigr)
	\frac12\eta_{\eps,v}''(c) \d c \\
&= \int_{-\infty}^u \bigl(P^\alpha(t,x,u)-P^\alpha(t,x,c)\bigr)
		\frac12\eta_{\eps,v}''(c) \d c
	-\int_u^\infty \bigl(P^\alpha(t,x,u)-P^\alpha(t,x,c)\bigr)
		\frac12\eta_{\eps,v}''(c) \d c \\
&= \int_v^u \bigl(P^\alpha(t,x,u)-P^\alpha(t,x,c)\bigr) \eta_{\eps,v}''(c) \d c
	\spliteq
	+\int_{-\infty}^v \bigl(P^\alpha(t,x,u)-P^\alpha(t,x,c)\bigr)
		\frac12\eta_{\eps,v}''(c) \d c
	-\int_v^\infty \bigl(P^\alpha(t,x,u)-P^\alpha(t,x,c)\bigr)
		\frac12\eta_{\eps,v}''(c) \d c \\
&= \left[\bigl(P^\alpha(t,x,u)-P^\alpha(t,x,c)\bigr) \eta_{\eps,v}'(c)\right]_v^u
	+\int_v^u \partial_3P^\alpha(t,x,c) \eta_{\eps,v}'(c) \d c
	\spliteq
	+\int_\setR \sign(v-c) \bigl(P^\alpha(t,x,u)-P^\alpha(t,x,c)\bigr)
		\frac12\eta_{\eps,v}''(c) \d c \\
&= P^\alpha_\eps(t,x,u,v)
	-\int_\setR \sign(v-c) P^\alpha(t,x,c) \frac12\eta_{\eps,v}''(c) \d c
\end{split}
\]
because in the last step the boundary terms vanish and
$\sign(v-c)P^\alpha(t,x,u)\eta_\eps''(c-v)$ is an odd function of $c-v$, and finally
\[
\begin{split}
&-\int_\setR \sign(u-c) \div_xP^\alpha(t,x,c) \frac12\eta_{\eps,v}''(c) \d c \\
&= -\int_{-\infty}^u \div_xP^\alpha(t,x,c) \frac12\eta_{\eps,v}''(c) \d c
	+\int_u^\infty \div_xP^\alpha(t,x,c) \frac12\eta_{\eps,v}''(c) \d c \\
&= -\int_v^u \div_xP^\alpha(t,x,c) \eta_{\eps,v}''(c) \d c
	\spliteq
	-\int_{-\infty}^v \div_xP^\alpha(t,x,c) \frac12\eta_{\eps,v}''(c) \d c
	+\int_v^\infty \div_xP^\alpha(t,x,c) \frac12\eta_{\eps,v}''(c) \d c \\
&= -\int_v^u \div_xP^\alpha(t,x,c) \eta_{\eps,v}''(c) \d c
	-\int_\setR \sign(v-c) \div_xP^\alpha(t,x,c) \frac12\eta_{\eps,v}''(c) \d c.
\end{split}
\]
Inserting these computations in the left hand side of \eqref{eq:entropy-ineq-reg-u} we get
\[
\begin{split}
& \int_\setR \Intd \bigl\{
	\abs{u-c}\partial_t\phi + \sign(u-c)\bigl[
		\bigl(P^\alpha(t,x,u)-P^\alpha(t,x,c)\bigr)\cdot\nabla_x\phi
		- \div_xP^\alpha(t,x,c)\phi \bigr] \bigr\} \frac12\eta_{\eps,v}''(c) \dx\dt \d c \\
&= \Intd \left\{
	\eta_{\eps,v}(u)\partial_t\phi + P^\alpha_\eps(t,x,u,v)\cdot\nabla_x\phi
	- \phi \int_v^u \div_xP^\alpha(t,x,s)\eta_{\eps,v}''(s)\d s
	\right\} \dx \dt
	\spliteq
	+ \Intd \partial_t\phi \int_\setR\abs{s}\rho_\eps(s)\ds \dx\dt
	- \int_\setR \Intd \sign(v-c) \div_x\bigl(P^\alpha(t,x,c)\phi\bigr)
		\frac12\eta_{\eps,v}''(c)\dx\dt \d c \\
&= \Intd \left\{
	\eta_{\eps,v}(u)\partial_t\phi + P^\alpha_\eps(t,x,u,v)\cdot\nabla_x\phi
	- \phi \int_v^u \div_xP^\alpha(t,x,s)\eta_{\eps,v}''(s)\d s
	\right\} \dx \dt
\end{split}
\]
because $\int_0^T \partial_t\phi \dt = 0$ and $\int_{\setR^n} \div_x\bigl(P^\alpha(t,x,c)\phi\bigr) \dx = 0$.

The right hand side of \eqref{eq:entropy-ineq-reg-u} remains unaltered because it does not depend on $c$ and the convolution kernel $\frac12\eta_{\eps,v}''(c)=\rho_\eps(c-v)$ is a probability measure.

Applying the same argument to \eqref{eq:entropy-ineq-reg-v}, the new pair of inequalities then becomes
\begin{subequations}\label{eq:entropy-ineq-eta}
\begin{equation}
\begin{split}
\Intd & \left\{
	\eta_{\eps,v}(u)\partial_t\phi + P^\alpha_\eps(t,x,u,v)\cdot\nabla_x\phi
	- \phi \int_v^u \div_xP^\alpha(t,x,s)\eta_{\eps,v}''(s)\d s
	\right\} \dx \dt \\
&\geq - \Intd \abs{\phi(t,x)} \d\mu_{0,t}(x)\dt
	- \Intd \abs{\nabla_x\phi(t,x)} \d\mu_{1,t}(x)\dt
	-\alpha e(P,\phi),
\end{split}
\end{equation}
\vspace{-.25cm}
\begin{equation}
\begin{split}
\Intd & \left\{
	\eta_{\eps,u}(v)\partial_t\phi + Q^\alpha_\eps(t,x,u,v)\cdot\nabla_x\phi
	- \phi \int_u^v \div_xQ^\alpha(t,x,s)\eta_{\eps,u}''(s)\d s
	\right\} \dx \dt \\
&\geq - \Intd \abs{\phi(t,x)} \d\nu_{0,t}(x)\dt
	- \Intd \abs{\nabla_x\phi(t,x)} \d\nu_{1,t}(x)\dt
	-\alpha e(Q,\phi).
\end{split}
\end{equation}
\end{subequations}

\paragraph{Doubling of variables.}
We perform the usual doubling of variables introduced by Kru\v{z}kov: given a test function $\tilde\phi(t,x,\tau,y)\in C^1_c(\setR\times\setR^n\times\setR\times\setR^n)$, we combine the entropy inequalities \eqref{eq:entropy-ineq-eta} for $u(t,x)$ and $v(\tau,y)$ and integrate w.r.t.\ the two additional variables. For conciseness, we omit the arguments of $u(t,x)$ and $v(\tau,y)$. We then obtain
\begin{subequations}\label{eq:entropy-ineq-doubl}
\begin{equation}\label{eq:entropy-ineq-doubl-u}
\begin{split}
\Intq & \left\{
	\eta_{\eps,v}(u)\partial_t\tilde\phi + P^\alpha_\eps(t,x,u,v)\cdot\nabla_x\tilde\phi
	- \tilde\phi \int_v^u \div_2P^\alpha(t,x,s)\eta_{\eps,v}''(s)\d s
		\right\} \dx \dy \dt \d\tau \\
&\geq -\Intd\left( \Intd \abs{\tilde\phi} \d\mu_{0,t}(x)\dt
	+ \Intd \abs{\nabla_x\tilde\phi} \d\mu_{1,t}(x)\dt \right) \dy\d\tau
	\spliteq
	-\alpha \Intd e\bigl(P,\tilde\phi(\plchldr,\plchldr,\tau,y)\bigr) \dy\d\tau,
\end{split}
\end{equation}
\vspace{-.75cm}
\begin{equation}\label{eq:entropy-ineq-doubl-v}
\begin{split}
\Intq & \left\{
	\eta_{\eps,u}(v)\partial_\tau\tilde\phi + Q^\alpha_\eps(\tau,y,u,v)\cdot\nabla_y\tilde\phi
	- \tilde\phi \int_u^v \div_2Q^\alpha(\tau,y,s)\eta_{\eps,u}''(s)\d s
		\right\} \dx \dy \dt \d\tau \\
&\geq -\Intd\left( \Intd \abs{\tilde\phi} \d\nu_{0,t}(y)\d\tau
	+ \Intd \abs{\nabla_x\tilde\phi} \d\nu_{1,t}(y)\d\tau \right) \dx\dt
	\spliteq
	-\alpha \Intd e\bigl(Q,\tilde\phi(t,x,\plchldr,\plchldr)\bigr) \dx\dt.
\end{split}
\end{equation}
\end{subequations}

We now consider a test function of the form
\[
\tilde\phi(t,x,\tau,y)
= \phi\oleft(\frac{t+\tau}2,\frac{x+y}2\right) \omega_\beta(x-y) \rho_\gamma(t-\tau),
\]
where $\phi\in C^1_c\bigl((0,T)\times\setR;[0,\infty)\bigr)$ is a space-time test function and $\beta,\gamma\in(0,1]$ are parameters. Once $\phi$ is fixed, for $\gamma$ small enough we have that $\tilde\phi(t,x,\tau,y)>0$ implies $t,\tau\in(0,T)$.
With this particular choice, we can bound the first integrals appearing in the right hand side of \eqref{eq:entropy-ineq-doubl-u} as
\[
\begin{split}
\Intq \abs{\tilde\phi} \d\mu_{0,t}(x)\dy\dt\d\tau
&= \Intq \abs*{\phi\oleft(\frac{t+\tau}2,\frac{x+y}2\right)} \omega_\beta(x-y)
	\rho_\gamma(t-\tau) \d\mu_{0,t}(x)\dy\dt\d\tau \\
&= \Intq \abs{\phi(s,z)} \omega_\beta\bigl(2(x-z)\bigr)
	\rho_\gamma\bigl(2(t-s)\bigr) \d\mu_{0,t}(x)2^n\dz\dt2\ds \\
&= \Intd \abs{\phi(s,z)} \Intd \omega_{\beta/2}(x-z)
	\rho_{\gamma/2}(t-s) \d\mu_{0,t}(x)\dt\dz\ds \\
&= \Intd \abs{\phi(s,z)} \mu_{0,\beta,\gamma}(s,z) \dz\ds,
\end{split}
\]
where $\mu_{0,\beta,\gamma}$ denotes the function
\[
\mu_{0,\beta,\gamma}(s,z)
= \Intd \omega_{\beta/2}(z-x) \rho_{\gamma/2}(s-t) \d\mu_{0,t}(x)\dt
= \bigl[(\rho_{\gamma/2}\omega_{\beta/2})*(\leb^1\otimes\mu_{0,t})\bigr](s,z),
\]
and similarly
\[
\begin{split}
& \Intq \abs{\nabla_x\tilde\phi} \d\mu_{1,t}(x)\dy\dt\d\tau \\
&\leq \Intq \frac12\abs*{\nabla_2\phi\oleft(\frac{t+\tau}2,\frac{x+y}2\right)}
	\omega_\beta(x-y) \rho_\gamma(t-\tau) \d\mu_{1,t}(x)\dy\dt\d\tau
	\spliteq
	+\Intq \abs*{\phi\oleft(\frac{t+\tau}2,\frac{x+y}2\right)}
	\abs{\nabla\omega_\beta(x-y)} \rho_\gamma(t-\tau) \d\mu_{1,t}(x)\dy\dt\d\tau \\
&= \frac12 \Intd \abs{\nabla_2\phi(s,z)} \mu_{1,\beta,\gamma}(s,z) \dz\ds
	+ \frac12 \Intd \abs{\phi(s,z)} \tilde\mu_{1,\beta,\gamma}(s,z) \dz\ds,
\end{split}
\]
where $\mu_{1,\beta,\gamma}$ and $\tilde\mu_{1,\beta,\gamma}$ denote the functions
\begin{align*}
\mu_{1,\beta,\gamma}(s,z)
&= \Intd \omega_{\beta/2}(z-x) \rho_{\gamma/2}(s-t) \d\mu_{1,t}(x)\dt
= \bigl[(\rho_{\gamma/2}\omega_{\beta/2})*(\leb^1\otimes\mu_{1,t})\bigr](s,z), \\
\tilde\mu_{1,\beta,\gamma}(s,z)
&= \Intd \abs{\nabla\omega_{\beta/2}(z-x)} \rho_{\gamma/2}(s-t) \d\mu_{1,t}(x)\dt
= \bigl[(\rho_{\gamma/2}\abs{\nabla\omega_{\beta/2}})*(\leb^1\otimes\mu_{1,t})\bigr](s,z).
\end{align*}
Of course analogous estimates hold with $\nu$ in place of $\mu$.
Informally speaking, observe that $\norm{\rho_{\gamma/2}\omega_{\beta/2}}_1=1$, therefore $\mu_{0,\beta,\gamma}$ and $\mu_{1,\beta,\gamma}$ are comparable to $\mu_0=\leb^1\otimes\mu_{0,t}$ and $\mu_1=\leb^1\otimes\mu_{1,t}$ respectively, whereas $\norm{\rho_{\gamma/2}\nabla\omega_{\beta/2}}_1\sim\beta^{-1}$, therefore $\tilde\mu_{1,\beta,\gamma}$ is comparable to $\beta^{-1}\mu_1$.
In order to estimate the errors $\int_0^T\!\!\int_{\setR^n} e\bigl(P,\tilde\phi(\plchldr,\plchldr,\tau,y)\bigr)\dy\d\tau$ and $\int_0^T\!\!\int_{\setR^n} e\bigl(Q,\tilde\phi(t,x,\plchldr,\plchldr)\bigr)\dx\dt$ we perform a similar computation as above replacing $\mu_{0,t}$ and $\mu_{1,t}$ with $\leb^n$.
We have that
\[
\begin{split}
\Intt \abs{\tilde\phi} \dx\dy\d\tau
&= \Intt \abs{\phi(s,z)} \omega_{\beta/2}(x-z)
	\rho_{\gamma/2}(t-s) \dx\dz\ds \\
&= \int_0^T \norm{\phi(s,\plchldr)}_{L^1(\setR^n)} \rho_{\gamma/2}(t-s) \ds
\leq \norm{\phi}_{L^\infty\bigl([0,T];L^1(\setR^n)\bigr)}
\end{split}
\]
and
\[
\begin{split}
\Intt \abs{\nabla_x\tilde\phi} \dx\dy\d\tau
&= \frac12 \Intt \abs{\nabla_2\phi(s,z)} \omega_{\beta/2}(x-z)
	\rho_{\gamma/2}(t-s) \dx\dz\ds
	\spliteq
	+ \frac12 \Intt \phi(s,z) \abs{\nabla\omega_{\beta/2}(x-z)}
	\rho_{\gamma/2}(t-s) \dx\dz\ds \\
&\leq \frac12\norm{\nabla_2\phi}_{L^\infty\bigl([0,T];L^1(\setR^n)\bigr)}
	+ \frac12
	\norm{\nabla\omega_{\beta/2}}_1 \norm{\phi}_{L^\infty\bigl([0,T];L^1(\setR^n)\bigr)},
\end{split}
\]
therefore, recalling the definition of $e(\plchldr,\plchldr)$ introduced in \eqref{eq:e-P-phi}, the error terms can be estimated as
\[
\begin{split}
&\Intd e\bigl(P,\tilde\phi(\plchldr,\plchldr,\tau,y)\bigr)\dy\d\tau
	+\Intd e\bigl(Q,\tilde\phi(t,x,\plchldr,\plchldr)\bigr)\dx\dt \\
&\leq
	\norm{\phi}_{L^\infty\bigl([0,T];L^1(\setR^n)\bigr)}
		\int_0^T [\Lip_3(\div_2P(t))+\Lip_3(\div_2Q(t))] \dt
	\spliteq
	+ \left( \norm{\nabla_2\phi}_{L^\infty\bigl([0,T];L^1(\setR^n)\bigr)}
	+ \norm{\nabla\omega_{\beta/2}}_1 \norm{\phi}_{L^\infty\bigl([0,T];L^1(\setR^n)\bigr)}
	\right) \int_0^T [\Lip_3(P(t))+\Lip_3(Q(t))] \dt \\
&\eqqcolon E(P,Q,\phi,\omega_{\beta/2}).
\end{split}
\]

Summing the two inequalities \eqref{eq:entropy-ineq-doubl} and using the previous estimates for the terms in the right hand side we get
\begin{equation}\label{eq:entropy-ineq-sum}
\begin{split}
\int_0^T \!\!\!\! & \int_0^T \!\!\!\! \int_{\setR^n} \!\! \int_{\setR^n} \biggl\{
	\eta_\eps(u-v)(\partial_t\tilde\phi+\partial_\tau\tilde\phi)
	+ \bigl[P^\alpha_\eps(t,x,u,v)\cdot\nabla_x\tilde\phi
	+ Q^\alpha_\eps(\tau,y,u,v)\cdot\nabla_y\tilde\phi\bigr]
	\\&\hspace{2cm}
	- \tilde\phi \left(\int_v^u \div_2P^\alpha(t,x,s)\eta_{\eps,v}''(s)\d s
	+ \int_u^v \div_2Q^\alpha(\tau,y,s)\eta_{\eps,u}''(s)\d s \right)
	\biggr\} \dx \dy \dt \d\tau \\
&\geq
	-\Intd \left(
	\abs{\phi} (\mu_{0,\beta,\gamma}+\nu_{0,\beta,\gamma})
	+\frac12 \abs{\nabla_2\phi} (\mu_{1,\beta,\gamma}+\nu_{1,\beta,\gamma})
	+\frac12 \abs{\phi} (\tilde\mu_{1,\beta,\gamma}+\tilde\nu_{1,\beta,\gamma})
	\right) \dz\ds
	\spliteq
	-\alpha E(P,Q,\phi,\omega_{\beta/2}) .
\end{split}
\end{equation}
Exploiting the identities
\begin{align*}
\partial_t\tilde\phi+\partial_\tau\tilde\phi &= \partial_1\phi\omega_\beta\rho_\gamma, &
\nabla_x\tilde\phi
	&= \frac12\nabla_2\phi\omega_\beta\rho_\gamma+\phi\nabla\omega_\beta\rho_\gamma, &
\nabla_y\tilde\phi
	&= \frac12\nabla_2\phi\omega_\beta\rho_\gamma-\phi\nabla\omega_\beta\rho_\gamma,
\end{align*}
we split the left hand side of \eqref{eq:entropy-ineq-sum} as $I+\II+\III$, where
\begin{align*}
I &= \Intq \eta_\eps(u-v) \partial_1\phi \omega_\beta \rho_\gamma \dx\dy\dt\d\tau, \\
\II
&= \Intq \biggl\{
	\frac12[P^\alpha_\eps(t,x,u,v) + Q^\alpha_\eps(\tau,y,u,v)]
		\cdot\nabla_2\phi \omega_\beta \rho_\gamma
	\spliteq \hspace{2cm}
	+ [P^\alpha_\eps(t,x,u,v) - Q^\alpha_\eps(\tau,y,u,v)]
		\cdot\nabla\omega_\beta \phi \rho_\gamma
	\biggr\} \dx \dy \dt \d\tau, \\
\III &= -\Intq \tilde\phi \left(\int_v^u \div_2P^\alpha(t,x,s)\eta_{\eps,v}''(s)\d s
	+ \int_u^v \div_2Q^\alpha(\tau,y,s)\eta_{\eps,u}''(s)\d s \right) \dx\dy\dt\d\tau.
\end{align*}

\paragraph{Dedoubling in time.}

We now perform the dedoubling of the time variables, which corresponds to passing to the limit $\gamma\to0$. All three terms $I$, $\II$ and $\III$ will be treated in a unified manner with the aid of \autoref{lem:dedoubling-time}. For each of the three we apply the lemma with the following choice of functions:
\begin{itemize}
\item dedoubling of $I$:
\begin{align*}
A(t,x,u,v) &= \eta_\eps(u-v), &
B &= \partial_1\phi, &
C &= \omega_\beta ;
\end{align*}
\item dedoubling of $\II$:
\begin{align*}
A(t,x,u,v) &= P^\alpha_\eps(t,x,u,v), &
B &= \nabla_2\phi, &
C &= \omega_\beta , \\
A(t,x,u,v) &= P^\alpha_\eps(t,x,u,v), &
B &= \phi, &
C &= \nabla\omega_\beta,
\end{align*}
and symmetrically for $Q^\alpha_\eps$;
\item dedoubling of $\III$:
\begin{align*}
A(t,x,u,v) &= \int_v^u \div_2P^\alpha(t,x,s)\eta''_\eps(s-v)\ds, &
B &= \phi, &
C &= \omega_\beta,
\end{align*}
and symmetrically for $Q^\alpha_\eps$.
\end{itemize}

As a result, in the limit $\gamma\to0$ we obtain the new terms
\begin{align*}
\widetilde I &= 
	\Intt \eta_\eps(u-v) \partial_1\phi \omega_\beta \dx\dy\dt , \\
\widetilde \II &= 
	\Intt \biggl\{
	\frac12[P^\alpha_\eps(t,x,u,v)
		+ Q^\alpha_\eps(t,y,u,v)]
		\cdot\nabla_2\phi \omega_\beta
	\spliteq\hspace{2cm}
	+ [P^\alpha_\eps(t,x,u,v)
		- Q^\alpha_\eps(t,y,u,v)]
		\cdot\nabla\omega_\beta \phi
	\biggr\} \dx \dy \dt , \\
\widetilde \III &= 
	-\Intt \left(\int_v^u \div_2P^\alpha(t,x,s)\eta_{\eps,v}''(s)\d s
		+ \int_u^v \div_2Q^\alpha(t,y,s)\eta_{\eps,u}''(s)\d s \right)
	\phi \omega_\beta \dx\dy\dt ,
\end{align*}
where we omitted the new implicit arguments
\[
u(t,x), \qquad
v(t,y), \qquad
\phi\oleft(t,\frac{x+y}2\right), \qquad
\omega_\beta(x-y).
\]

Let us now focus on the right hand side of \eqref{eq:entropy-ineq-sum}. We have that
\[
\begin{split}
\lim_{\gamma\to0} & \Intd \abs{\phi(s,z)}\mu_{0,\beta,\gamma}(s,z) \dz\ds \\
&= \lim_{\gamma\to0} \Intq \abs{\phi(s,z)}\rho_{\gamma/2}(s-t) \omega_{\beta/2}(z-x)
	\d\mu_{0,t}(x)\dz\dt\ds\\
&= \lim_{\gamma\to0} \Intd [\abs{\phi(\plchldr,z)}*\rho_{\gamma/2}](t)
	(\omega_{\beta/2}*\mu_{0,t})(z) \dz\dt \\
&= \Intd \abs{\phi(t,z)} (\omega_{\beta/2}*\mu_{0,t})(z) \dz\dt
\end{split}
\]
because the function $(t,z)\mapsto[\abs{\phi(\plchldr,z)}*\rho_{\gamma/2}](t)$ converges uniformly to $\abs{\phi(t,z)}$ and have bounded support, whereas $(\omega_{\beta/2}*\mu_{0,t})\dz\dt$ is a locally finite measure. With similar computation for the other terms, we get
\begin{equation}\label{eq:dedoubled-rhs}
\begin{split}
\lim_{\gamma\to0} & \Intd \left(
	\abs{\phi} (\mu_{0,\beta,\gamma}+\nu_{0,\beta,\gamma})
	+\frac12 \abs{\nabla_2\phi} (\mu_{1,\beta,\gamma}+\nu_{1,\beta,\gamma})
	+\frac12 \abs{\phi} (\tilde\mu_{1,\beta,\gamma}+\tilde\nu_{1,\beta,\gamma})
	\right) \dz\ds \\
&= \Intd \left(
	\abs{\phi} \omega_{\beta/2}*(\mu_{0,t}+\nu_{0,t})
	+\frac12 \abs{\nabla_2\phi} \omega_{\beta/2}*(\mu_{1,t}+\nu_{1,t})
	+\frac12 \abs{\phi} \abs{\nabla\omega_{\beta/2}}*(\mu_{1,t}+\nu_{1,t})
	\right) \dz\dt.
\end{split}
\end{equation}

\paragraph{Integration by parts and chain rule.}

Integrating by parts w.r.t.\ $x$ the term involving $\nabla\omega_\beta$, we can rewrite $\widetilde\II$ as
\[
\begin{split}
\widetilde\II
&= \Intt \biggl\{
	\frac12[P^\alpha_\eps(t,x,u,v) + Q^\alpha_\eps(t,y,u,v)] \cdot\nabla_2\phi \omega_\beta
	\spliteq \hspace{2cm}
	- \div_x[P^\alpha_\eps(t,x,u,v) - Q^\alpha_\eps(t,y,u,v)] \phi \omega_\beta
	\spliteq \hspace{2cm}
	- \frac12 [P^\alpha_\eps(t,x,u,v) - Q^\alpha_\eps(t,y,u,v)] \cdot\nabla_2\phi \omega_\beta
	\biggr\} \dx \dy \dt \\
&= \Intt \biggl\{
	Q^\alpha_\eps(t,y,u,v) \cdot\nabla_2\phi \omega_\beta
	- \div_x[P^\alpha_\eps(t,x,u,v) - Q^\alpha_\eps(t,y,u,v)] \phi \omega_\beta
	\biggr\} \dx \dy \dt \\
&\eqqcolon \widetilde\II_1 + \widetilde\II_2 ,
\end{split}
\]
where with an abuse of notation we denoted
\[
\begin{split}
\widetilde\II_2
&= -\Intt \div_x[P^\alpha_\eps(t,x,u,v) - Q^\alpha_\eps(t,y,u,v)] \phi \omega_\beta
	\dx \dy \dt
= -\Intt \phi \omega_\beta \d\sigma_{t,y}(x) \dy \dt,
\end{split}
\]
with $\sigma_{t,y}=\div_x\bigl[P^\alpha_\eps\bigl(t,x,u(t,x),v(t,y)\bigr) - Q^\alpha_\eps\bigl(t,y,u(t,x),v(t,y)\bigr)\bigr]$ being a measure in the variable $x$ parametrized by $(t,y)$.

Indeed, since the fluxes $P^\alpha_\eps$ and $Q^\alpha_\eps$ are $C^1$ and $u$ is $BV_\loc$, by the chain rule \cite[§ 13.2]{VolpertBV}\footnote{The chain rule is applied to the $C^1$ function $(x,u)\mapsto P^\alpha_\eps(t,x,u,v)$ with $(t,v)$ fixed and the $BV_\loc$ function $x\mapsto(x,u(x))$.} we get for every $t,y$ that the function $P^\alpha_\eps\bigl(t,x,u(t,x),v(t,y)\bigr) - Q^\alpha_\eps\bigl(t,y,u(t,x),v(t,y)\bigr)$ is $BV_\loc(\setR^n;\setR^n)$ and its divergence is the measure in the variable $x$ given by
\begin{equation*}
\begin{split}
\sigma_{t,y}
&= \div_2 P^\alpha_\eps\bigl(t,x,u(t,x),v(t,y)\bigr)\leb^n
	\spliteq
	+ \partial_3 \bigl[P^\alpha_\eps\bigl(t,x,u(t,x),v(t,y)\bigr)
	- Q^\alpha_\eps\bigl(t,y,u(t,x),v(t,y)\bigr)\bigr]
		\cdot \bigl(\D^a_x u(t)+\D^c_x u(t)\bigr)
	\spliteq
	+ \biggl\{\bigl[P^\alpha_\eps\bigl(t,x,u^+(t,x),v(t,y)\bigr)
	- Q^\alpha_\eps\bigl(t,y,u^+(t,x),v(t,y)\bigr)\bigr]
	\spliteq\qquad
	-\bigl[P^\alpha_\eps\bigl(t,x,u^-(t,x),v(t,y)\bigr)
	- Q^\alpha_\eps\bigl(t,y,u^-(t,x),v(t,y)\bigr)\bigr]\biggr\}
	\cdot\bm{n} \haus^{n-1}\rvert_{J_{u(t)}} ,
\end{split}
\end{equation*}
where $D^a_xu(t)$, $D^c_xu(t)$ and $\bigl(u^+(t,\plchldr)-u^-(t,\plchldr)\bigr)\bm{n}\haus^{n-1}\rvert_{J_{u(t)}}$ represent the absolutely continuous, the Cantor and the jump part of the derivative of the $BV_\loc$ function $u(t,\plchldr)$ respectively \cite{AFP}.

\paragraph{Estimates of $\widetilde I$, $\widetilde\II_1$, $\widetilde\II_2$, $\widetilde\III$.}


In this section of the proof we make use of several pointwise estimates, all of which are relevant only when $(t,x),(t,y)\in\supp\Theta$ and $u,v\in[0,R(T)]$, hence can be expressed in terms of \autoref{rmk:notation-norms}.

We fix the specific test function
\[
\phi(t,x) = \Theta(t,x) \theta_r(t-t_1) \theta_r(t_2-t),
\]
where $\theta\in C^\infty(\setR;[0,1])$ with $\theta(t)=0$ for $t\leq0$ and $\theta(t)=1$ for $\theta\geq1$, and $\theta_r(t)=\theta(t/r)$. For $r<(t_2-t_1)/2$ we have
\[
\begin{split}
\partial_1\phi\oleft(t,\frac{x+y}2\right)
&= \partial_1\Theta\oleft(t,\frac{x+y}2\right) \theta_r(t-t_1) \theta_r(t_2-t)
	+ \Theta\oleft(t,\frac{x+y}2\right) [\theta_r'(t-t_1) - \theta_r'(t_2-t)] .
\end{split}
\]
Exploiting the fact that $\eta'_{\eps,u}$ has always the same sign in the interval $[u\wedge v,u\vee v]$, we can estimate
\[
\abs{Q^\alpha_\eps(t,y,u,v)}
= \abs*{\int_u^v \partial_3Q^\alpha(t,y,s)\eta'_{\eps,u}(s)\d s}
\leq \Lip_3(Q^\alpha(t)) \int_u^v \abs*{\eta'_{\eps,u}(s)}\d s
= \Lip_3(Q^\alpha(t)) \eta_\eps(u-v).
\]
Using this Lipschitz estimate for $Q^\alpha_\eps$ and the hypothesis \eqref{eq:Theta-slope} we get
\begin{equation}\label{eq:time-survivors}
\begin{split}
\widetilde I + \widetilde\II_1
&= \Intt \eta_\eps(u-v) \partial_1\phi \omega_\beta \dx\dy\dt
	+ \Intt Q^\alpha_\eps(t,y,u,v) \cdot\nabla_2\phi \omega_\beta \dx\dy\dt \\
&= \Intt \bigl[ \eta_\eps(u-v)\partial_1\Theta + Q^\alpha_\eps(t,y,u,v)\cdot\nabla_2\Theta
	\bigr] \theta_r(t-t_1) \theta_r(t_2-t) \omega_\beta \dx\dy\dt
	\spliteq
	+ \Intt \eta_\eps(u-v) \Theta [\theta_r'(t-t_1) - \theta_r'(t_2-t)]
		\omega_\beta \dx\dy\dt \\
&\leq \Intt \bigl[\partial_1\Theta+\Lip_3(Q^\alpha(t))\abs{\nabla_2\Theta}\bigr]
	\eta_\eps(u-v) \theta_r(t-t_1) \theta_r(t_2-t) \omega_\beta \dx\dy\dt
	\spliteq
	+ \Intt \eta_\eps(u-v) \Theta [\theta_r'(t-t_1) - \theta_r'(t_2-t)]
		\omega_\beta \dx\dy\dt \\
&\!\!\overset{\eqref{eq:Theta-slope}}{\leq} \!\! \Intt \eta_\eps\bigl(u(t,x)-v(t,y)\bigr)
	\Theta\oleft(t,\frac{x+y}2\right)
	[\theta_r'(t-t_1) - \theta_r'(t_2-t)] \omega_\beta(x-y) \dx\dy\dt.
\end{split}
\end{equation}
Observe that for $\beta\leq1$ one has
\begin{equation}\label{eq:BV-estimate}
\begin{split}
&\int_{\setR^n} \!\! \int_{\setR^n} \abs*{u(t,x)-u\oleft(t,\frac{x+y}2\right)}	
	\Theta\oleft(t,\frac{x+y}2\right) \omega_\beta(x-y) \dx\dy \\
&\leq \frac\beta2 \norm{\Theta(t)}_\infty \TV_{\Omega_t}\bigl(u(t,\plchldr)\bigr)
\leq \frac\beta2 \norm{\Theta(t)}_\infty B(t).
\end{split}
\end{equation}
Using this observation and the fact that $\eta_\eps$ is $1$-Lipschitz,
we can continue the estimate \eqref{eq:time-survivors} as
\begin{equation}
\begin{split}\label{eq:time-survivors-estimate}
\widetilde I + \widetilde\II_1
&\leq \Intt \eta_\eps\oleft(u\oleft(t,\frac{x+y}2\right)
		-v\oleft(t,\frac{x+y}2\right)\right)
	\Theta\oleft(t,\frac{x+y}2\right)
	[\theta_r'(t-t_1) - \theta_r'(t_2-t)] \omega_\beta(x-y) \dx\dy\dt
	\spliteq
	+\Intt \abs*{u(t,x)-u\oleft(t,\frac{x+y}2\right)}
	\Theta\oleft(t,\frac{x+y}2\right)
	\abs{\theta_r'(t-t_1) - \theta_r'(t_2-t)} \omega_\beta(x-y) \dx\dy\dt
	\spliteq
	+\Intt \abs*{v\oleft(t,\frac{x+y}2\right)-v(t,y)}
	\Theta\oleft(t,\frac{x+y}2\right)
	\abs{\theta_r'(t-t_1) - \theta_r'(t_2-t)} \omega_\beta(x-y) \dx\dy\dt \\
&\leq \Intd \eta_\eps\bigl(u(t,x)-v(t,x)\bigr) \Theta(t,x)
	[\theta_r'(t-t_1) - \theta_r'(t_2-t)] \dx\dt
	\spliteq
	+\beta \int_0^T \norm{\Theta(t)}_\infty B(t)
	\abs{\theta_r'(t-t_1) - \theta_r'(t_2-t)} \dt
	\\
&\leq \Intd \eta_\eps\bigl(u(t,x)-v(t,x)\bigr) \Theta(t,x)
	[\theta_r'(t-t_1) - \theta_r'(t_2-t)] \dx\dt
	+2\beta \norm{\Theta}_\infty B(t_2),
\end{split}
\end{equation}
where we used that $B(t)\leq B(t_2)$, $\norm{\theta'_r}_1=1$, and the integrand is supported in $[t_1,t_2]$.

Let us now deal with the third term.
\begin{equation}\label{eq:estimate-III}
\begin{split}
\widetilde\III
&= -\Intt \left\{
	\int_v^u \div_2P^\alpha(t,x,s)\eta_{\eps,v}''(s)\d s
	+ \int_u^v \div_2Q^\alpha(t,y,s)\eta_{\eps,u}''(s)\d s
	\right\} \phi \omega_\beta \dx\dy\dt \\
&= -\Intt \biggl\{
	\div_2P^\alpha(t,x,u)\eta_\eps'(u-v)
		- \int_v^u \div_2\partial_3P^\alpha(t,x,s)\eta_{\eps,v}'(s)\d s
	\spliteq\qquad
	- \div_2Q^\alpha(t,y,v)\eta_\eps'(u-v)
		- \int_u^v \div_2\partial_3Q^\alpha(t,y,s)\eta_{\eps,u}'(s)\d s
	\biggr\} \phi \omega_\beta \dx\dy\dt \\
&= -\Intt \biggl\{
	\Bigl[ \div_2P^\alpha(t,x,u) - \div_2P^\alpha(t,x,v) \Bigr]
		\spliteq\hspace{2cm}
	+ \Bigl[ \div_2P^\alpha(t,x,v) - \div_2P^\alpha(t,y,v) \Bigr]
		\spliteq\hspace{2cm}
	+ \Bigl[ \div_2P^\alpha(t,y,v) - \div_2Q^\alpha(t,y,v) \Bigr]
	\biggr\} \eta'_\eps(u-v) \phi \omega_\beta \dx\dy\dt
	\spliteq
	+\Intt \biggl\{
	\int_v^u \div_2\partial_3P^\alpha(t,x,s)\eta_{\eps,v}'(s)\d s
	+\int_u^v \div_2\partial_3Q^\alpha(t,y,s)\eta_{\eps,u}'(s)\d s
	\biggr\} \phi \omega_\beta \dx\dy\dt ,
\end{split}
\end{equation}
where as usual $u=u(t,x)$ and $v=v(t,y)$.
We now estimate the differences of divergences inside the square brackets:
\begin{align*}
\abs*{\div_2P^\alpha(t,x,u) - \div_2P^\alpha(t,x,v)}
&\leq \Lip_3(\div_2P^\alpha(t)) \abs{u-v} , \\
\abs*{\div_2P^\alpha(t,x,v) - \div_2P^\alpha(t,y,v)}
&\leq \Lip_2(\div_2P^\alpha(t)) \abs{x-y} , \\
\abs*{\div_2P^\alpha(t,y,v) - \div_2Q^\alpha(t,y,v)}
&\leq \norm{\div_2(P^\alpha-Q^\alpha)(t,y,\plchldr)}_{L^\infty([0,R(T)])}.
\end{align*}
On the other hand, using the fact that $\abs{\eta_\eps'}\leq1$, we have
\begin{equation}\label{eq:div2-P-eps-u-v}
\abs*{\int_v^u \div_2\partial_3P^\alpha(t,x,s)\eta_{\eps,v}'(s)\d s}
\leq \Lip_3(\div_2 P^\alpha(t)) \abs{u-v},
\end{equation}
and similarly for the integral involving $Q^\alpha$.
Therefore, continuing \eqref{eq:estimate-III}, we get
\begin{equation}\label{eq:estimate-III-intermediate}
\begin{split}
\widetilde\III
\leq \intt \Bigl\{&
	\Bigl[2\Lip_3(\div_2P^\alpha(t))+\Lip_3(\div_2Q^\alpha(t)) \Bigr] \abs{u-v}
	+ \Lip_2(\div_2P^\alpha(t)) \abs{x-y}
	\spliteq
	+ \norm{\div_2(P^\alpha-Q^\alpha)(t,y,\plchldr)}_{L^\infty([0,R(T)])}
		\Bigr\} \Theta\omega_\beta \dx\dy\dt.
\end{split}
\end{equation}
Using \eqref{eq:BV-estimate} we can estimate
\begin{equation}\label{eq:ux-vy}
\begin{split}
\intrr & \abs{u(t,x)-v(t,y)} \Theta\oleft(t,\frac{x+y}2\right) \omega_\beta(x-y) \dx\dy\\
&\leq \intrr \abs*{u\oleft(t,\frac{x+y}2\right)-v\oleft(t,\frac{x+y}2\right)}
	\Theta\oleft(t,\frac{x+y}2\right) \omega_\beta(x-y) \dx\dy
	\spliteq
	+ \intrr \abs*{u(t,x)-u\oleft(t,\frac{x+y}2\right)}	
	\Theta\oleft(t,\frac{x+y}2\right) \omega_\beta(x-y) \dx\dy
	\spliteq
	+ \intrr \abs*{v\oleft(t,\frac{x+y}2\right)-v(t,y)}	
	\Theta\oleft(t,\frac{x+y}2\right) \omega_\beta(x-y) \dx\dy \\
&\leq \int_{\setR^n} \abs{u(t,x)-v(t,x)}\Theta(t,x) \dx
	+ \beta \norm{\Theta(t)}_\infty B(t),
\end{split}
\end{equation}
\[
\begin{split}
\intrr \abs{x-y} \Theta\oleft(t,\frac{x+y}2\right) \omega_\beta(x-y) \dx\dy
&= \intrr \abs{2x-2z} \Theta(t,z) \omega_\beta(2x-2z) 2^n\dx\dz \\
&= \int_{\setR^n} \abs{x}\omega_\beta(x)\dx \norm{\Theta(t)}_1
\leq \beta \norm{\Theta(t)}_1 ,
\end{split}
\]
and
\[
\begin{split}
&\intt \norm{\div_2(P^\alpha-Q^\alpha)(t,y,\plchldr)}_{L^\infty([0,R(T)])}
	\Theta\oleft(t,\frac{x+y}2\right) \omega_\beta(x-y) \dx\dy\dt \\
&= \intd \norm{\div_2(P^\alpha-Q^\alpha)(t,y,\plchldr)}_{L^\infty([0,R(T)])}
		\Theta(t,y) \dy\dt
	\spliteq
	+ \intt \norm{\div_2(P^\alpha-Q^\alpha)(t,y,\plchldr)}_{L^\infty([0,R(T)])}
	\left[\Theta\oleft(t,\frac{x+y}2\right)-\Theta(t,y)\right] \omega_\beta(x-y) \dx\dy\dt\\
&\leq \intd \norm{\div_2(P^\alpha-Q^\alpha)(t,y,\plchldr)}_{L^\infty([0,R(t)])}
	\Theta(t,y) \dy\dt
	\spliteq
	+ \frac\beta2 \int_{t_1}^{t_2} \Lip_2(\Theta(t))
	\int_{\Omega_t} \norm{\div_2(P^\alpha-Q^\alpha)(t,y,\plchldr)}_{L^\infty([0,R(t)])}
	\dy \dt,
\end{split}
\]
which combined and inserted in \eqref{eq:estimate-III-intermediate} lead to
\begin{equation}\label{eq:estimate-III-final}
\begin{split}
\widetilde\III
&\leq \intd [2\Lip_3(\div_2P^\alpha(t))+\Lip_3(\div_2Q^\alpha(t))]
	\abs{u(t,x)-v(t,x)}\Theta(t,x) \dx\dt
	\spliteq
	+ \beta \norm{\Theta}_\infty B(t_2) \int_{t_1}^{t_2}
		[2\Lip_3(\div_2P^\alpha(t))+\Lip_3(\div_2Q^\alpha(t))] \dt
	\spliteq
	+ \beta \norm{\Theta}_{L^\infty_t L^1_x}
		\int_{t_1}^{t_2} \Lip_2(\div_2P^\alpha(t)) \dt
	\spliteq
	+ \intd \norm{\div_2(P^\alpha-Q^\alpha)(t,y,\plchldr)}_{L^\infty([0,R(t)])}
		\Theta(t,y) \dy\dt
	\spliteq
	+ \frac\beta2 \norm{\Theta}_{L^\infty_t\Lip_x}
	\int_{t_1}^{t_2} \!\!\!\! \int_{\Omega_t}
		\norm{\div_2(P^\alpha-Q^\alpha)(t,y,\plchldr)}_{L^\infty([0,R(t)])} \dy\dt.
\end{split}
\end{equation}

Finally, we turn to the estimate for $\widetilde\II_2$, which will be obtained by considering separately the three integrals which constitute it. First of all, by \eqref{eq:div2-P-eps-u-v} and \eqref{eq:ux-vy} we have
\begin{equation}\label{eq:II2-div2}
\begin{split}
&\abs*{\Intt \div_2 P^\alpha_\eps\bigl(t,x,u(t,x),v(t,y)\bigr) \phi\omega_\beta \dx\dy\dt}\\
&\leq \int_{t_1}^{t_2} \Lip_3\bigl(\div_2P^\alpha(t)\bigr) \left(
	\int_{\setR^n} \abs{u(t,x)-v(t,x)}\Theta(t,x) \dx
	+ \beta \norm{\Theta(t)}_\infty B(t) \right) \dt \\
&\leq \intd \Lip_3\bigl(\div_2P^\alpha(t)\bigr) \abs{u(t,x)-v(t,x)}\Theta(t,x) \dx \dt
	+ \beta \norm{\Theta}_\infty B(t_2) \int_{t_1}^{t_2} \Lip_3(\div_2P^\alpha(t)) \dt.
\end{split}
\end{equation}
Let us now estimate
\[
\begin{split}
&\abs*{\partial_3(P^\alpha_\eps-Q^\alpha_\eps)(t,x,u,v)} \\
&= \abs*{\partial_3\left(
	\int_v^u \partial_3P^\alpha(t,x,s)\eta'_\eps(s-v)\d s
	-\int_u^v \partial_3Q^\alpha(t,x,s)\eta'_\eps(s-u)\d s
	\right)} \\
&= \abs*{\partial_3\left(
	\int_v^u \partial_3P^\alpha(t,x,s)\eta'_\eps(s-v)\d s
	-\int_v^u \partial_3Q^\alpha(t,x,s)\eta'_\eps(u-s)\d s
	\right)} \\
&= \abs*{\partial_3P^\alpha(t,x,u)\eta'_\eps(u-v)
	- \int_v^u \partial_3Q^\alpha(t,x,s)\eta''_\eps(u-s)\d s}.
\end{split}
\]
If $\abs{u-v}\leq\eps$ then this is less than
\[
\begin{split}
&\abs*{\partial_3P^\alpha(t,x,u)\eta'_\eps(u-v)
	+ \int_v^u \partial_3\partial_3Q^\alpha(t,x,s)\eta'_\eps(s-u)\d s
	+ \partial_3Q^\alpha(t,x,v)\eta'_\eps(v-u)} \\
&\leq \abs*{\partial_3(P^\alpha-Q^\alpha)(t,x,u)}
	+ \abs*{\partial_3Q^\alpha(t,x,u)-\partial_3Q^\alpha(t,x,v)}
	+ \int_v^u \abs{\partial_3\partial_3Q^\alpha(t,x,s)}\d s \\
&\leq \Lip_3\bigl((P^\alpha-Q^\alpha)(t,x)\bigr)
	+ 2\Lip_3\bigl(\partial_3Q^\alpha(t)\bigr)\abs{u-v} \\
&\leq \Lip_3\bigl((P^\alpha-Q^\alpha)(t,x)\bigr)
	+ 2\Lip_3\bigl(\partial_3Q^\alpha(t)\bigr)\eps.
\end{split}
\]
Otherwise, if $\abs{u-v}>\eps$ then it is less than
\[
\begin{split}
& \abs*{\partial_3(P^\alpha-Q^\alpha)(t,x,u)}
	+ \abs*{\partial_3Q^\alpha(t,x,u)\eta'_\eps(u-v)
	-\int_v^u \partial_3Q^\alpha(t,x,s)\eta''_\eps(u-s)\d s} \\
&\leq \Lip_3\bigl((P^\alpha-Q^\alpha)(t,x)\bigr)
	+ \abs*{\int_v^u [\partial_3Q^\alpha(t,x,u)-\partial_3Q^\alpha(t,x,s)]\eta''_\eps(u-s)\d s} \\
&\leq \Lip_3\bigl((P^\alpha-Q^\alpha)(t,x)\bigr)
	+ \Lip_3\bigl(\partial_3Q^\alpha(t)\bigr)\eps
\end{split}
\]
since $\supp(\eta''_\eps) \subset[-\eps,\eps]$ and $\int_v^u \eta''_\eps(u-s)\ds=\eta'_\eps(u-v)$. Therefore in both cases we have
\begin{equation}\label{eq:partial3-P-Q}
\abs*{\partial_3(P^\alpha_\eps-Q^\alpha_\eps)(t,x,u,v)}
\leq \Lip_3\bigl((P^\alpha-Q^\alpha)(t,x)\bigr)
	+ 2\Lip_3\bigl(\partial_3Q^\alpha(t)\bigr)\eps.
\end{equation}
Moreover,
\begin{equation}\label{eq:lip2-partial3-Q}
\begin{split}
\abs*{\partial_3[Q^\alpha_\eps(t,x,u,v) - Q^\alpha_\eps(t,y,u,v)]}
&= \abs*{
	\int_u^v [\partial_3Q^\alpha(t,x,s)-\partial_3Q^\alpha(t,y,s)]\eta''_\eps(s-u)\d s} \\
&\leq \Lip_2\bigl(\partial_3Q^\alpha(t)\bigr)\abs{x-y}
\end{split}
\end{equation}
because $\eta''_\eps(s-u)$ has mass less than $1$ in $[u\wedge v,u\vee v]$.
With \eqref{eq:partial3-P-Q} and \eqref{eq:lip2-partial3-Q} we can estimate
\begin{equation}\label{eq:II2-Da-Dc}
\begin{split}
& \abs*{\Intt \phi\omega_\beta \partial_3 \bigl[P^\alpha_\eps\bigl(t,x,u(t,x),v(t,y)\bigr)
	- Q^\alpha_\eps\bigl(t,y,u(t,x),v(t,y)\bigr)\bigr]
		\cdot \d\bigl(\D^a_x u(t)+\D^c_x u(t)\bigr)(x) \dy\dt } \\
&\leq \abs*{ \Intt \phi\omega_\beta
	\partial_3 (P^\alpha_\eps-Q^\alpha_\eps)\bigl(t,x,u(t,x),v(t,y)\bigr)
		\cdot \d\bigl(\D^a_x u(t)+\D^c_x u(t)\bigr)(x) \dy\dt }
	\spliteq
	+\abs*{ \Intt \phi\omega_\beta
	\partial_3 \bigl[Q^\alpha_\eps\bigl(t,x,u(t,x),v(t,y)\bigr)
	- Q^\alpha_\eps\bigl(t,y,u(t,x),v(t,y)\bigr)\bigr]
		\cdot \d\bigl(\D^a_x u(t)+\D^c_x u(t)\bigr)(x) \dy\dt } \\
&\leq \Intt \phi\omega_\beta
	[\Lip_3\bigl((P^\alpha-Q^\alpha)(t,x)\bigr)
		+ 2\Lip_3\bigl(\partial_3Q^\alpha(t)\bigr)\eps]
	\d\abs{\D^a_x u(t)+\D^c_x u(t)}(x) \dy\dt
	\spliteq
	+ \Intt \Lip_2\bigl(\partial_3Q^\alpha(t)\bigr)
		\phi\omega_\beta \abs{x-y} \d\abs{\D^a_x u(t)+\D^c_x u(t)}(x) \dy\dt \\
&\leq \intt \Lip_3\bigl((P^\alpha-Q^\alpha)(t,x)\bigr) \Theta(t,x) \omega_\beta(x-y)
		\d\abs{\D^a_x u(t)+\D^c_x u(t)}(x) \dy\dt
	\spliteq
	+\intt \Lip_3\bigl((P^\alpha-Q^\alpha)(t,x)\bigr)
	\abs*{\Theta\oleft(t,\frac{x+y}2\right) - \Theta(t,x)}\omega_\beta(x-y)
		\d\abs{\D^a_x u(t)+\D^c_x u(t)}(x) \dy\dt
	\spliteq
	+\Intt 2\Lip_3\bigl(\partial_3Q^\alpha(t)\bigr)\eps \phi\omega_\beta
		\d\abs{\D^a_x u(t)+\D^c_x u(t)}(x) \dy\dt
	\spliteq
	+ \Intt \Lip_2\bigl(\partial_3Q^\alpha(t)\bigr)
		\phi\omega_\beta \abs{x-y} \d\abs{\D^a_x u(t)+\D^c_x u(t)}(x) \dy\dt \\
&\leq B(t_2) \int_{t_1}^{t_2}
	\norm{\Lip_3\bigl((P^\alpha-Q^\alpha)(t,\plchldr)\bigr) \Theta(t,\plchldr)}_\infty \dt
	+\frac\beta2 \norm{\Theta}_{L^\infty_t\Lip_x} B(t_2)
		\int_{t_1}^{t_2} \Lip_3\bigl((P^\alpha-Q^\alpha)(t)\bigr) \dt
	\spliteq
	+2\eps \norm{\Theta}_\infty B(t_2)
		\int_{t_1}^{t_2} \Lip_3\bigl(\partial_3Q^\alpha(t)\bigr) \dt
	+\beta \norm{\Theta}_\infty B(t_2) \int_{t_1}^{t_2}
		\Lip_2\bigl(\partial_3Q^\alpha(t)\bigr) \dt .
\end{split}
\end{equation}
In a similar fashion we have
\[
\begin{split}
&\Intt \biggl\{\bigl[P^\alpha_\eps\bigl(t,x,u^+(t,x),v(t,y)\bigr)
	- Q^\alpha_\eps\bigl(t,y,u^+(t,x),v(t,y)\bigr)\bigr]
	\spliteq\hspace{2cm}
	-\bigl[P^\alpha_\eps\bigl(t,x,u^-(t,x),v(t,y)\bigr)
	- Q^\alpha_\eps\bigl(t,y,u^-(t,x),v(t,y)\bigr)\bigr]\biggr\}
	\cdot\bm{n} \phi\omega_\beta \d\haus^{n-1}\rvert_{J_{u(t)}}(x) \dy\dt \\
&= \Intt \biggl\{
	(P^\alpha_\eps-Q^\alpha_\eps)\bigl(t,x,u^+(t,x),v(t,y)\bigr)
	-(P^\alpha_\eps-Q^\alpha_\eps)\bigl(t,x,u^-(t,x),v(t,y)\bigr)\biggr\}
	\cdot\bm{n} \phi\omega_\beta \d\haus^{n-1}\rvert_{J_{u(t)}}(x) \dy\dt
	\spliteq
	+\Intt \biggl\{
	\bigl[Q^\alpha_\eps\bigl(t,x,u^+(t,x),v(t,y)\bigr)
		-Q^\alpha_\eps\bigl(t,x,u^-(t,x),v(t,y)\bigr)\bigr]
	\spliteq\qquad\qquad
	-\bigl[Q^\alpha_\eps\bigl(t,y,u^+(t,x),v(t,y)\bigr)
		-Q^\alpha_\eps\bigl(t,y,u^-(t,x),v(t,y)\bigr)\bigr]\biggr\}
	\cdot\bm{n} \phi\omega_\beta \d\haus^{n-1}\rvert_{J_{u(t)}}(x) \dy\dt \\
&= \Intt \int_{u^-(t,x)}^{u^+(t,x)}
	\partial_3(P^\alpha_\eps-Q^\alpha_\eps)\bigl(t,x,s,v(t,y)\bigr)
	\cdot\bm{n} \phi\omega_\beta \ds \d\haus^{n-1}\rvert_{J_{u(t)}}(x) \dy\dt
	\spliteq
	+\Intt \int_{u^-(t,x)}^{u^+(t,x)}
		[\partial_3Q^\alpha_\eps\bigl(t,x,s,v(t,y)\bigr)
		-\partial_3Q^\alpha_\eps\bigl(t,y,s,v(t,y)\bigr)]
	\cdot\bm{n} \phi\omega_\beta \ds \d\haus^{n-1}\rvert_{J_{u(t)}}(x) \dy\dt.
\end{split}
\]
By \eqref{eq:partial3-P-Q} and \eqref{eq:lip2-partial3-Q} we have
\[
\begin{split}
&\abs*{\int_{u^-(t,x)}^{u^+(t,x)} \bigl\{
	\partial_3(P^\alpha_\eps-Q^\alpha_\eps)\bigl(t,x,s,v(t,y)\bigr)
	+[\partial_3Q^\alpha_\eps\bigl(t,x,s,v(t,y)\bigr)
		-\partial_3Q^\alpha_\eps\bigl(t,y,s,v(t,y)\bigr)]
	\bigr\}\ds } \\
&\leq \bigl[\Lip_3\bigl((P^\alpha-Q^\alpha)(t,x)\bigr)
	+ 2\Lip_3\bigl(\partial_3Q^\alpha(t)\bigr)\eps
	+\Lip_2\bigl(\partial_3Q^\alpha(t)\bigr)\abs{x-y}\bigr] \cdot \abs{u^+(t,x)-u^-(t,x)},
\end{split}
\]
therefore the previous integral can be estimated in absolute value as
\begin{equation}\label{eq:II2-jump}
\begin{split}
&\biggl\lvert \Intt \biggl\{\bigl[P^\alpha_\eps\bigl(t,x,u^+(t,x),v(t,y)\bigr)
	- Q^\alpha_\eps\bigl(t,y,u^+(t,x),v(t,y)\bigr)\bigr]
	\spliteq\hspace{2cm}
	-\bigl[P^\alpha_\eps\bigl(t,x,u^-(t,x),v(t,y)\bigr)
	- Q^\alpha_\eps\bigl(t,y,u^-(t,x),v(t,y)\bigr)\bigr]\biggr\}
	\cdot\bm{n} \phi\omega_\beta \d\haus^{n-1}\rvert_{J_{u(t)}}(x) \dy\dt \biggr\rvert \\
&\leq \Intt
	\bigl[\Lip_3\bigl((P^\alpha-Q^\alpha)(t,x)\bigr)
	+ 2\Lip_3\bigl(\partial_3Q^\alpha(t)\bigr)\eps\bigr] \phi \omega_\beta
	\abs{u^+-u^-} \d\haus^{n-1}\rvert_{J_{u(t)}}(x) \dy\dt
	\spliteq
	+ \Intt \Lip_2\bigl(\partial_3Q^\alpha(t)\bigr)\abs{x-y} \phi \omega_\beta
	\abs{u^+-u^-} \d\haus^{n-1}\rvert_{J_{u(t)}}(x) \dy\dt \\
&\leq B(t_2) \int_{t_1}^{t_2}
	\norm{\Lip_3\bigl((P^\alpha-Q^\alpha)(t,\plchldr)\bigr) \Theta(t,\plchldr)}_\infty \dt
	+2\eps \norm{\Theta}_\infty B(t_2)
		\int_{t_1}^{t_2} \Lip_3\bigl(\partial_3Q^\alpha(t)\bigr) \dt
	\spliteq
	+\frac\beta2 \norm{\Theta}_{L^\infty_t\Lip_x} B(t_2)
		\int_{t_1}^{t_2} \Lip_3\bigl((P^\alpha-Q^\alpha)(t)\bigr) \dt
	+\beta \norm{\Theta}_\infty B(t_2) \int_{t_1}^{t_2}
		\Lip_2\bigl(\partial_3Q^\alpha(t)\bigr) \dt .
\end{split}
\end{equation}
In conclusion, from \eqref{eq:II2-div2}, \eqref{eq:II2-Da-Dc} and \eqref{eq:II2-jump} we deduce
\begin{equation}\label{eq:estimate-II2}
\begin{split}
\widetilde\II_2
&\leq \intd \Lip_3\bigl(\div_2P^\alpha(t)\bigr) \abs{u(t,x)-v(t,x)}\Theta(t,x) \dx \dt
	\spliteq
	+2B(t_2) \int_{t_1}^{t_2}
	\norm{\Lip_3\bigl((P^\alpha-Q^\alpha)(t,\plchldr)\bigr) \Theta(t,\plchldr)}_\infty \dt
	+4\eps \norm{\Theta}_\infty B(t_2)
		\int_{t_1}^{t_2} \Lip_3\bigl(\partial_3Q^\alpha(t)\bigr) \dt
	\spliteq
	+\beta \norm{\Theta}_{L^\infty_t\Lip_x} B(t_2)
		\int_{t_1}^{t_2} \Lip_3\bigl((P^\alpha-Q^\alpha)(t)\bigr) \dt
	+2\beta \norm{\Theta}_\infty B(t_2) \int_{t_1}^{t_2}
		\Lip_2\bigl(\partial_3Q^\alpha(t)\bigr) \dt
	\spliteq
	+\beta \norm{\Theta}_\infty B(t_2) \int_{t_1}^{t_2} \Lip_3(\div_2P^\alpha(t)) \dt.
\end{split}
\end{equation}

\paragraph{Conclusion.}

After the dedoubling in time of \eqref{eq:entropy-ineq-sum}, using the estimates \eqref{eq:dedoubled-rhs}, \eqref{eq:time-survivors-estimate}, \eqref{eq:estimate-III-final} and \eqref{eq:estimate-II2} we deduce that
\begin{equation}\label{eq:conclusion-inequality}
\begin{split}
& \Intd \eta_\eps\bigl(u(t,x)-v(t,x)\bigr) \Theta(t,x)
		[\theta_r'(t-t_1) - \theta_r'(t_2-t)] \dx\dt
	\spliteq
	+\intd [3\Lip_3(\div_2P^\alpha(t))+\Lip_3(\div_2Q^\alpha(t))]
		\abs{u(t,x)-v(t,x)}\Theta(t,x) \dx\dt
	\spliteq
	+ \intd \norm{\div_2(P^\alpha-Q^\alpha)(t,y,\plchldr)}_{L^\infty([0,R(t)])}
		\Theta(t,y) \dy\dt
	\spliteq
	+2B(t_2) \int_{t_1}^{t_2}
	\norm{\Lip_3\bigl((P^\alpha-Q^\alpha)(t,\plchldr)\bigr) \Theta(t,\plchldr)}_\infty \dt
	\spliteq
	+ \beta \norm{\Theta}_\infty B(t_2) \left( 2 + \int_{t_1}^{t_2}
		[3\Lip_3(\div_2P^\alpha(t))+\Lip_3(\div_2Q^\alpha(t))
			+2\Lip_2(\partial_3Q^\alpha(t))] \dt \right)
	\spliteq
	+\beta \norm{\Theta}_{L^\infty_t L^1_x} \int_{t_1}^{t_2} \Lip_2(\div_2P^\alpha(t)) \dt
	+\beta \norm{\Theta}_{L^\infty_t\Lip_x} B(t_2)
		\int_{t_1}^{t_2} \Lip_3\bigl((P^\alpha-Q^\alpha)(t)\bigr) \dt
	\spliteq
	+ \frac\beta2 \norm{\Theta}_{L^\infty_t\Lip_x}
	\int_{t_1}^{t_2} \!\!\!\! \int_{\Omega_t}
	\norm{\div_2(P^\alpha-Q^\alpha)(t,y,\plchldr)}_{L^\infty([0,R(t)])} \dy\dt
	\spliteq
	+4\eps \norm{\Theta}_\infty B(t_2)
		\int_{t_1}^{t_2} \Lip_3\bigl(\partial_3Q^\alpha(t)\bigr) \dt \\
&\geq - \Intd \left(
	\abs{\phi} \omega_{\beta/2}*(\mu_{0,t}+\nu_{0,t})
	+\frac12 \abs{\nabla_2\phi} \omega_{\beta/2}*(\mu_{1,t}+\nu_{1,t})
	+\frac12 \abs{\phi} \abs{\nabla\omega_{\beta/2}}*(\mu_{1,t}+\nu_{1,t})
	\right) \dz\dt
	\spliteq
	-\alpha E(P,Q,\phi,\omega_{\beta/2}).
\end{split}
\end{equation}
Letting $\eps\to0$ and $\alpha\to0$ in this order and estimating the integral in the right hand side we get\footnote{Notice that for $\alpha\to0$ we have $\Lip(P^\alpha-Q^\alpha) \to \Lip(P-Q)$ and similarly for the other Lipschitz norms.}
\[
\begin{split}
& \Intd \abs{u(t,x)-v(t,x)} \Theta(t,x) [\theta_r'(t-t_1) - \theta_r'(t_2-t)] \dx\dt
	\spliteq
	+\intd [3\Lip_3(\div_2P(t))+\Lip_3(\div_2Q(t))]
		\abs{u(t,x)-v(t,x)}\Theta(t,x) \dx\dt
	\spliteq
	+\intd \norm{\div_2(P-Q)(t,y,\plchldr)}_{L^\infty([0,R(t)])}
		\Theta(t,y) \dy\dt
	\spliteq
	+2B(t_2) \int_{t_1}^{t_2}
		\norm{\Lip_3\bigl((P-Q)(t,\plchldr)\bigr) \Theta(t,\plchldr)}_\infty \dt
	\spliteq
	+ \beta \norm{\Theta}_\infty B(t_2) \left( 2 + \int_{t_1}^{t_2}
		[3\Lip_3(\div_2P(t))+\Lip_3(\div_2Q(t))
			+2\Lip_2(\partial_3Q(t))] \dt \right)
	\spliteq
	+\beta \norm{\Theta}_{L^\infty_t L^1_x} \int_{t_1}^{t_2} \Lip_2(\div_2P(t)) \dt
	+\beta \norm{\Theta}_{L^\infty_t\Lip_x} B(t_2)
		\int_{t_1}^{t_2} \Lip_3\bigl((P-Q)(t)\bigr) \dt
	\spliteq
	+\frac\beta2 \norm{\Theta}_{L^\infty_t\Lip_x}
	\int_{t_1}^{t_2} \!\!\!\! \int_{\Omega_t}
	\norm{\div_2(P-Q)(t,y,\plchldr)}_{L^\infty([0,R(t)])} \dy\dt \\
&\geq -\norm{\Theta}_\infty \int_{t_1}^{t_2}
	(\mu_{0,t}+\nu_{0,t})\bigl(\supp\Theta(t,\plchldr)\bigr)_\beta \dt
	\spliteq
	-\left(\frac12\norm{\nabla_2\Theta}_\infty
		+\beta^{-1} \norm{\Theta}_\infty \norm{\nabla\omega}_1\right)
		\int_{t_1}^{t_2} (\mu_{1,t}+\nu_{1,t})\bigl(\supp\Theta(t,\plchldr)\bigr)_\beta \dt.
\end{split}
\]
For $\beta\leq1$, letting $r\to0$, rearranging the terms and enlarging the set $\bigl(\supp\Theta(t,\plchldr)\bigr)_\beta$ to $\Omega_t=\bigl(\supp\Theta(t,\plchldr)\bigr)_1$ we deduce
\[
\begin{split}
&\left[\int_{\setR^n} \abs{u(t,x)-v(t,x)} \Theta(t,x) \dx \right]_{t_1}^{t_2} \\
&\leq \intd [3\Lip_3(\div_2P(t))+\Lip_3(\div_2Q(t))] \abs{u(t,x)-v(t,x)}\Theta(t,x) \dx\dt
	\spliteq
	+\intd \norm{\div_2(P-Q)(t,y,\plchldr)}_{L^\infty([0,R(t)])}
		\Theta(t,y) \dy\dt
	\spliteq
	+2B(t_2) \int_{t_1}^{t_2}
		\norm{\Lip_3\bigl((P-Q)(t,\plchldr)\bigr) \Theta(t,\plchldr)}_\infty \dt
	\spliteq
	+\norm{\Theta}_\infty \int_{t_1}^{t_2}
		(\mu_{0,t}+\nu_{0,t})(\Omega_t) \dt
	+\frac12\norm{\nabla_2\Theta}_\infty \int_{t_1}^{t_2}
		(\mu_{1,t}+\nu_{1,t})(\Omega_t) \dt
	\spliteq
	+ \beta \norm{\Theta}_\infty B(t_2) \left( 2 + \int_{t_1}^{t_2}
		[3\Lip_3(\div_2P(t))+\Lip_3(\div_2Q(t))
			+2\Lip_2(\partial_3Q(t))] \dt \right)
	\spliteq
	+\beta \norm{\Theta}_{L^\infty_t L^1_x} \int_{t_1}^{t_2} \Lip_2(\div_2P(t)) \dt
	+\beta \norm{\Theta}_{L^\infty_t\Lip_x} B(t_2)
		\int_{t_1}^{t_2} \Lip_3\bigl((P-Q)(t)\bigr) \dt
	\spliteq
	+\frac\beta2 \norm{\Theta}_{L^\infty_t\Lip_x}
		\int_{t_1}^{t_2} \!\!\!\! \int_{\Omega_t}
		\norm{\div_2(P-Q)(t,y,\plchldr)}_{L^\infty([0,R(t)])} \dy\dt
	\spliteq
	+c_n \beta^{-1} \norm{\Theta}_\infty \int_{t_1}^{t_2}
		(\mu_{1,t}+\nu_{1,t})(\Omega_t) \dt,
\end{split}
\]
where in the last line we estimated $\norm{\nabla\omega}_1\leq c_n$.

The last remaining step is to (almost) optimize in $\beta$. Recalling the definitions \eqref{eq:M-C} of $M$ and $C$ and picking
\[
\beta = \min\oleft\{ M(t_1,t_2)^{1/2}, 1 \right\}
\]
we have that the last four lines of the right hand side, which are of the form $C(t_1,t_2)\beta+c_n\norm{\Theta}_\infty\beta^{-1}M(t_1,t_2)$, can be bounded by
\[
C(t_1,t_2) \min\oleft\{ M(t_1,t_2)^{1/2}, 1 \right\}
	+c_n \norm{\Theta}_\infty
		\max\oleft\{ M(t_1,t_2)^{1/2}, M(t_1,t_2) \right\},
\]
which leads to the inequality \eqref{eq:MRS} claimed in the statement of the theorem.
\end{proof}

\begin{remark}\label{rmk:entropic-statement-proof}
We provide here details on how to recover the analogous result for the entropic case stated in \autoref{rmk:entropic-statement}. Heuristically speaking, this could be done by performing the dedoubling in space too, which amounts to setting $\beta=0$ in the last part of the proof. In view of recycling as much as possible of the given proof, one can proceed as follows instead.

All the steps of the proof are retraced unaltered except the last one. At the beginning of the step marked \textbf{Conclusion}, in equation \eqref{eq:conclusion-inequality}, we get rid of the integral in the right hand side which depends on the measures $\mu_i,\nu_i$ since they are all zero by assumption. Next we take the limit as $\beta\to0$ and as a result eliminate all the terms with $\beta$, appearing in the rows 5 through 7 of the left hand side. Then the limits $\eps\to0$ and $\alpha\to0$ are taken in the same described way and we immediately reach the conclusion because the last step (optimization in $\beta$) is now unnecessary.

Notice in particular that we can get rid of the assumptions regarding $\Lip_2(\div_2 P(t))$ and $\Lip_2(\partial_3 Q(t))$ because they do not appear in the final statement. Indeed, in \eqref{eq:conclusion-inequality} there are $\Lip_2(\div_2 P^\alpha(t))$ and $\Lip_2(\partial_3 Q^\alpha(t))$, which are finite because the fluxes are regularized with $\alpha>0$, but these terms disappear by taking $\beta\to0$ before $\alpha\to0$.
\end{remark}

\begin{proof}[Proof of \autoref{prop:gronwall}]
By assumption we can directly apply \autoref{thm:stability} and deduce the validity of \eqref{eq:MRS}. Using the assumptions \eqref{eq:gronwall-assumptions} we can estimate the terms
\[
\begin{split}
\intd & \norm{\div_2(P-Q)(t,x,\plchldr)}_{L^\infty([0,R(t)])}
		\Theta(t,x) \dx\dt
	+2B(t_2) \int_{t_1}^{t_2}
		\norm{\Lip_3\bigl((P-Q)(t,\plchldr)\bigr) \Theta(t,\plchldr)}_\infty \dt \\
&\leq [1+2B(t_2)] \int_{t_1}^{t_2} h(t)
	\int_{\setR^n} \abs{u(t,x)-v(t,x)} \Theta(t,x) \dx\dt.
\end{split}
\]
Letting
\[
w(t) = \int_{\setR^n} \abs{u(t,x)-v(t,x)} \Theta(t,x) \dx,
\]
the above estimate combined with \eqref{eq:MRS} leads to an inequality of the form
\[
w(t_2) \leq w(t_1) + \Phi(t_1,t_2) + \int_{t_1}^{t_2} f(t)w(t)\dt,
\]
where $f(t) = [3\Lip_3(\div_2P(t))+\Lip_3(\div_2Q(t))] + [1+2B(T)]h(t)$ and $\Phi(t_1,t_2)$ encompasses all the remaining terms.
Applying Grönwall's inequality leads to
\[
w(t_2) \leq [w(t_1) + \Phi(t_1,t_2)] \exp\oleft(\int_{t_1}^{t_2} f(t)\dt\right).\qedhere
\]
\end{proof}


%
%
%
%


\begin{lemma}[Dedoubling in time]\label{lem:dedoubling-time}
Let $A,\partial_3A,\partial_4A \in L^1_\loc\bigl([0,T);L^\infty_\loc(\setR^n\times\setR\times\setR)\bigr)$, $B\in C_c((0,T)\times\setR^n)$, $C\in C^1_c(\setR^n)$, and $u,v\in C\bigl([0,T);L^1_\loc(\setR^n)\bigr)\cap L^\infty_\loc([0,T)\times\setR^n)$.
Then
\[
\begin{split}
\lim_{\gamma\to0} &
	\Intq A\bigl(t,x,u(t,x),v(\tau,y)\bigr)
	B\oleft(\frac{t+\tau}2,\frac{x+y}2\right)
	C(x-y) \rho_\gamma(t-\tau) \dx\dy\dt\d\tau \\
&= \Intt A\bigl(t,x,u(t,x),v(t,y)\bigr)
	B\oleft(t,\frac{x+y}2\right) C(x-y) \dx\dy\dt.
\end{split}
\]
\end{lemma}

\begin{proof}
Let $T'\in(0,T)$ and $r>0$ be such that $\supp B\subseteq[2\gamma,T'-2\gamma]\times B_r(0)$ and $\supp C\subseteq B_r(0)$ for $\gamma>0$ sufficiently small and let $R:[0,T)\to[0,\infty)$ be an increasing function such that $u(t,x),v(t,x)\in[0,R(t)]$ for every $(t,x)\in[0,T)\times B_r(0)$.
Let
\[
H(\gamma)
= \sup \set{\norm{v(t_1,\plchldr)-v(t_2,\plchldr)}_{L^1(B_{2r}(0))}}
	{t_1,t_2\in[0,T'],\ \abs{t_1-t_2}\leq\gamma}
\]
denote the modulus of continuity of the map $[0,T']\to L^1\bigl(B_{2r}(0)\bigr):t\mapsto v(t,\plchldr)$ and let
\[
K(\gamma)
= \sup \set{\norm{B(t_1,\plchldr)-B(t_2,\plchldr)}_{L^\infty(\setR^n)}}
	{t_1,t_2\in[0,T'],\ \abs{t_1-t_2}\leq\gamma}
\]
denote the modulus of continuity of the map $[0,T']\to L^\infty(\setR^n):t\mapsto B(t,\plchldr)$.

Thanks to the assumed continuities we have both $H(\gamma),K(\gamma)\to0$ for $\gamma\to0$. Moreover
\[
\int_0^{T'} \!\!\!\! \int_{B_{2r}} \! \int_{B_{2r}}
	\abs*{A\bigl(t,x,u(t,x),v(t,y)\bigr)} \dx\dy\dt
\leq \int_0^{T'} \norm{A(t,\plchldr,\plchldr,\plchldr)}_{L^\infty(B_{2r},[0,R(t)],[0,R(t)])} \dt
< \infty.
\]

Then, using the notation $\Lip_4\bigl(A(t)\bigr) = \norm{\partial_4A(t,\plchldr,\plchldr,\plchldr)}_{L^\infty(B_{2r},[0,R(t)],[0,R(t)])}$, we have
\[
\begin{split}
& \biggl\lvert \Intq A\bigl(t,x,u(t,x),v(\tau,y)\bigr)
	B\oleft(\frac{t+\tau}2,\frac{x+y}2\right)
	C(x-y) \rho_\gamma(t-\tau) \dx\dy\dt\d\tau
	\spliteq
	- \Intt A\bigl(t,x,u(t,x),v(t,y)\bigr)
	B\oleft(t,\frac{x+y}2\right) C(x-y) \dx\dy\dt \biggr\rvert \\
&\leq \Intq \abs*{A\bigl(t,x,u(t,x),v(\tau,y)\bigr) - A\bigl(t,x,u(t,x),v(t,y)\bigr)}
	\abs*{B\oleft(\frac{t+\tau}2,\frac{x+y}2\right)}
	\abs{C(x-y)} \rho_\gamma(t-\tau) \dx\dy\dt\d\tau
	\spliteq
	+ \Intq \abs*{A\bigl(t,x,u(t,x),v(t,y)\bigr)}
	\abs*{B\oleft(\frac{t+\tau}2,\frac{x+y}2\right)-B\oleft(t,\frac{x+y}2\right)}
	\abs{C(x-y)} \rho_\gamma(t-\tau) \dx\dy\dt\d\tau \\
&\leq \Intq \Lip_4\bigl(A(t)\bigr) \abs{v(\tau,y)-v(t,y)}
	\abs*{B\oleft(\frac{t+\tau}2,\frac{x+y}2\right)}
	\abs{C(x-y)} \rho_\gamma(t-\tau) \dx\dy\dt\d\tau
	\spliteq
	+ \int_0^{T'} \!\!\!\! \int_{\setR^n} \!\! \int_{\setR^n}
	\abs*{A\bigl(t,x,u(t,x),v(t,y)\bigr)}
	K(\gamma) \bm{1}_{B_r}\oleft(\frac{x+y}2\right)
	\norm{C}_\infty \bm{1}_{B_r}(x-y) \dx\dy\dt \\
&\leq \int_0^{T'} \Lip_4\bigl(A(t)\bigr) \norm{B}_\infty \norm{C}_1
	\int_0^{T'-\gamma} \norm{v(\tau,\plchldr)-v(t,\plchldr)}_{L^1(B_{2r})}
		\rho_\gamma(t-\tau) \d\tau \d t
	\spliteq
	+ K(\gamma) \norm{C}_\infty \int_0^{T'} \!\!\!\! \int_{B_{2r}} \! \int_{B_{2r}}
	\abs*{A\bigl(t,x,u(t,x),v(t,y)\bigr)} \dx\dy\dt \\
&\leq \norm{B}_\infty \norm{C}_1 H(\gamma) \norm{\Lip_4\bigl(A(t)\bigr)}_{L^1([0,T'])}
	+ K(\gamma) \norm{C}_\infty
	\int_0^{T'} \norm{A(t,\plchldr,\plchldr,\plchldr)}_{L^\infty(B_{2r},[0,R(t)],[0,R(t)])} \dt.
\end{split}
\]

The right hand side of the estimate converges to $0$ when $\gamma\to0$ and this proves the limit in the statement.
\end{proof}

\section{Applications}

In this section we demonstrate how to employ \autoref{thm:stability} and \autoref{prop:gronwall} to both theoretical and numerical applications.

\subsection{Conditional existence via Cauchy sequences in \texorpdfstring{$L^1_\loc$}{L¹-loc}}\label{sec:existence}

Our stability result \autoref{thm:stability} leads to a conditional existence theorem of entropy solutions for the problems $\partial_tu+\div_x\bigl(P(t,x,u)\bigr)=0$ and $\partial_tu+\div_x\bigl(P[u](t,x,u)\bigr)=0$ whenever one can provide a sequence of quasi-entropy solutions with vanishing errors.

\begin{theorem}[Conditional existence, fixed flux]\label{thm:existence-fixed}
Let $P:[0,T)\times\setR^n\times[0,\infty)\to\setR^n$ be a flux satisfying \autoref{as:flux} with the property that $\Lip_3\bigl(P(t)\bigr)\in L^\infty_\loc\bigl([0,T)\bigr)$ and let $u_0\in L^\infty(\setR^n)\cap BV(\setR^n)$ be a non-negative initial datum.
Let $(u_k)_{k\in\setN}$ be a sequence of $(\mu_{k,0},\mu_{k,1})$-quasi-entropy solutions according to \autoref{def:mu-quasi-entropy-solution} for the problem
\begin{equation}\label{eq:cons-law-existence-fixed}
\partial_t u + \div_x\bigl(P(t,x,u)\bigr) = 0
\end{equation}
with non-negative initial datum $u_k(0)\in L^\infty(\setR^n)\cap BV(\setR^n)$ converging to $u_0$ in $L^1_\loc(\setR^n)$,
satisfying the uniform bounds
\[
u_k(t,\plchldr) \leq R(t), \qquad
\norm{u_k(t,\plchldr)}_{BV} \leq B(t), \qquad \forall t\in[0,T), k\in\setN,
\]
for some increasing functions $R,B:[0,T)\to[0,\infty)$.
Assume in addition that the measures $\mu_{k,i,t}\in L^1\bigl([0,T);\Meas_+(\setR^n)\bigr)$ vanish in the sense that for every $r>1/T$ we have
\begin{equation}\label{eq:vanishing-mu-k}
\lim_{k\to\infty} \int_0^{T-1/r} \abs{\mu_{k,i,t}}(B_r) \dt = 0.
\end{equation}

Then for every $t\in[0,T)$ the sequence $u_k(t,\plchldr)$ converges in $L^1_\loc(\setR^n)$ to a function $u(t,x)$ which is an entropy solution of \eqref{eq:cons-law-existence-fixed} according to \autoref{def:entropy-solution} with initial datum $u_0$.
\end{theorem}

\begin{proof}
Let $\theta\in C^\infty\bigl(\setR;[0,1]\bigr)$ be such that $\theta(r)=0$ for $r\leq0$, $\theta(r)=1$ for $r\geq1$, and $\theta'(r)\leq2$ for every $r$.
Given $t_2\in(0,T)$ and a radius $r>0$, let $c=\norm{\Lip_3\bigl(P(t)\bigr)}_{L^\infty([0,t_2])}$ and define $\Theta\in C^\infty_c\bigl([0,t_2]\times\setR^n;[0,1]\bigr)$ as
\[
\Theta(t,x) = \theta\bigl(c(t_2-t) - \abs{x} + 1+r \bigr).
\]
Then we have
\[
\partial_t\Theta(t,x) = -c\theta'\bigl(c(t_2-t) - \abs{x} + 1+r \bigr)
= -c\abs{\nabla_x\Theta(t,x)},
\]
so that \eqref{eq:Theta-slope} is satisfied up to $t=t_2$.

For every $k,l\in\setN$ we can then apply \autoref{thm:stability} to the pair $(u_k,u_l)$ with $t_1=0$, from which we deduce that for every $t\in(0,t_2)$ we have
\[
\begin{split}
\int_{\setR^n} \abs{u_k(t,x)-u_l(t,x)} \Theta(t,x) \dx
&\leq 4\int_0^t \Lip_3(\div_2P(s))
		\int_{\setR^n} \abs{u_k(s,x)-u_l(s,x)}\Theta(s,x) \dx\ds
	\spliteq
	+ \int_{\setR^n} \abs{u_k(0,x)-u_l(0,x)} \Theta(0,x) \dx
	+ \eps_{k,l}
\end{split}
\]
where
\[
\begin{split}
\eps_{k,l}
&= \int_0^t
		(\mu_{k,0,s}+\mu_{l,0,s})(\Omega_s) \ds
	+ \int_0^t
		(\mu_{k,1,s}+\mu_{l,1,s})(\Omega_s) \ds
	\spliteq
	+C(0,t) \min\oleft\{ M_{k,l}(0,t)^{1/2}, 1 \right\}
		+ c_n \max\oleft\{ M_{k,l}(0,t)^{1/2}, M_{k,l}(0,t) \right\}
\end{split}
\]
is a quantity that goes to $0$ as $k,l\to0$ thanks to \eqref{eq:vanishing-mu-k}.
By Grönwall's inequality we get
\[
\int_{\setR^n} \abs{u_k(t,x)-u_l(t,x)} \Theta(t,x) \dx
\leq \left[
	\int_{\setR^n} \abs{u_k(0,x)-u_l(0,x)} \Theta(0,x) \dx
	+\eps_{k,l}\right]
	\exp\oleft(4\int_0^t \Lip_3(\div_2P(s))\ds\right).
\]
Considering that
\[
\norm{u_k(t,\plchldr)-u_l(t,\plchldr)}_{L^1(B_r)}
\leq \int_{\setR^n} \abs{u_k(t,x)-u_l(t,x)} \Theta(t,x) \dx
\]
for $t\in[0,t_2]$ and using the fact that $u_k(0,\plchldr)$ is a Cauchy sequence in $L^1_\loc(\setR^n)$, we deduce that $u_k(t,\plchldr)$ is a Cauchy sequence in $L^1(B_r)$.

From the arbitrariness of $t_2$ and $r$ we obtain that there exists a limit function $u(t,\plchldr)$ to which the sequence $u_k(t,\plchldr)$ converges in $L^1_\loc(\setR^n)$ for every $t\in[0,T)$.

We need to show that $u$ satisfies the entropy inequality \eqref{eq:def-entropy-ineq-hyp}. Since there are some technical difficulties passing to the limit the entropy inequality with Kruzkov type entropies $\abs{u-c}$, this will be done by regularizing the absolute value mimicking the second step in the main proof of the stability theorem and passing to the limit this regularized inequality instead. The desired entropy inequality of Kruzkov type for the limit function $u$ is then recovered by approximation.

By assumption, $u_k$ satisfies
\[
\begin{split}
\Intd & \bigl\{
	\abs{u_k-c'}\partial_t\phi + \sign(u_k-c')\bigl[
		\bigl(P(t,x,u_k)-P(t,x,c')\bigr)\cdot\nabla_x\phi - \div_xP(t,x,c')\phi
	\bigr] \bigr\} \dx \dt \\
&\geq -	\Intd \abs{\phi(t,x)} \d\mu_{k,0,t}(x)\dt
	- \Intd \abs{\nabla_x\phi(t,x)} \d\mu_{k,1,t}(x)\dt
\end{split}
\]
for every constant $c'\in[0,\infty)$ and non-negative test function $\phi\in C^\infty_c\bigl((0,T)\times\setR^n;[0,\infty)\bigr)$.
Letting
\[
\eta_\eps(u) = \int_\setR \abs{s} \rho_\eps(u-s) \d s
	-\int_\setR \abs{s} \rho_\eps(s) \d s
\]
and its translation $\eta_{\eps,c}(u)=\eta_\eps(u-c)$, with a similar computation as in the step \emph{Regularization of the absolute value} on page~\pageref{par:reg-abs}, convolving the previous quasi-entropy inequality with $\frac12\eta''_{\eps,c}(c')$ we deduce that
\begin{equation}\label{eq:exist-reg-entropy-ineq-k}
\begin{split}
\Intd & \left\{
	\eta_{\eps,c}(u_k)\partial_t\phi + P_\eps(t,x,u_k,c)\cdot\nabla_x\phi
	- \phi \int_c^{u_k} \div_xP(t,x,s)\eta_{\eps,c}''(s)\d s
	\right\} \dx \dt \\
&\geq - \Intd \abs{\phi(t,x)} \d\mu_{k,0,t}(x)\dt
	- \Intd \abs{\nabla_x\phi(t,x)} \d\mu_{k,1,t}(x)\dt,
\end{split}
\end{equation}
where
\[
P_\eps(t,x,u,c) = \int_c^u \partial_3P(t,x,s)\eta_{\eps,c}'(s) \d s.
\]

We now want to pass to the limit \eqref{eq:exist-reg-entropy-ineq-k} as $k\to\infty$.
For the first two terms in the left hand side, we use the fact that the dependence on $u_k$ is Lipschitz, indeed
\begin{align*}
\abs{\eta_{\eps,c}(u_k)-\eta_{\eps,c}(u)} &\leq \abs{u_k-u}, \\
\abs{P_\eps(t,x,u_k,c)-P_\eps(t,x,u,c)} &\leq \Lip_3(P(t)) \abs{u_k-u}.
\end{align*}
For the third term we use the fact that
\[
\lim_{k\to\infty} \Intd \int_{u(t,x)}^{u_k(t,x)} \phi(t,x)\div_xP(t,x,s)\eta_{\eps,c}''(s)\ds\dx\dt = 0
\]
because the integrand is a fixed $L^1$ function which is integrated over the set
\[
\set{(t,x,s)\in [0,T)\times\setR^n\times[0,R(T)]}
	{(t,x)\in\supp\phi,\ s\in[u_k(t,x)\wedge u(t,x),u_k(t,x)\vee u(t,x)]}
\]
whose measure is
\[
\iint_{\supp\phi} \abs{u_k(t,x)-u(t,x)}\dx\dt \to 0, \qquad \text{for $k\to\infty$}.
\]
For the right hand side we use the assumption \eqref{eq:vanishing-mu-k} to deduce that both terms are vanishing.
As a consequence we obtain that $u$ satisfies the regularized entropy inequality
\[
\Intd \left\{
	\eta_{\eps,c}(u)\partial_t\phi + P_\eps(t,x,u,c)\cdot\nabla_x\phi
	- \phi \int_c^{u} \div_xP(t,x,s)\eta_{\eps,c}''(s)\d s
	\right\} \dx \dt
\geq 0,
\]
for every constant $c\in[0,\infty)$ and non-negative test function $\phi\in C^\infty_c\bigl((0,T)\times\setR^n;[0,\infty)\bigr)$.

Letting $\eps\to0$ gives us the desired entropy inequality for $u$ with Kruzkov entropies:
\[
\Intd \bigl\{
	\abs{u-c}\partial_t\phi + \sign(u-c)\bigl[
		\bigl(P(t,x,u)-P(t,x,c)\bigr)\cdot\nabla_x\phi - \div_xP(t,x,c)\phi
	\bigr] \bigr\} \dx \dt
\geq 0.
\]

By a standard argument, setting $c=0$ and $c=R(T)$ in the previous inequality we deduce that $u$ solves the conservation law $\partial_t u + \div_x\bigl(P(t,x,u)\bigr)=0$ too.

We now prove the continuity $u\in C\bigl([0,T);L^1_\loc(\setR^n)\bigr)$. From the fact that $u$ solves the conservation law we deduce that for every $\phi\in C^\infty_c(\setR^n)$ the map
\[
[0,T)\to\setR: t\mapsto \int_{\setR^n} u(t,x)\phi(x)\dx
\]
is differentiable with finite derivative equal to
\[
-\int_{\setR^n} P\bigl(t,x,u(t,x)\bigr) \nabla_x\phi(x)\dx.
\]
In particular, the map is continuous, hence $u(t_n,\plchldr)\weakto u(t,\plchldr)$ whenever $t_n\to t$. This limit holds also in the strong topology $L^1_\loc(\setR^n)$ because the the equi-boundedness in $BV(\setR^n)$ of the sequence implies its relatively compactness in $L^1_\loc(\setR^n)$.
\end{proof}

%
%
%
%
%
%
%

\begin{theorem}[Conditional existence, non-local flux]\label{thm:existence-nonlocal}
Let
\[
\begin{split}
P:\mathscr{S}\bigl([0,T),\setR^n\bigr)
	&\to \{p:[0,T)\times\setR^n\times[0,\infty)\to\setR^n\} \\
u &\mapsto P[u]
\end{split}
\]
be a map associating to every function $u$ a flux $P[u]$ satisfying the \autoref{as:flux},
and let $u_0\in L^\infty(\setR^n)\cap BV(\setR^n)$ be a non-negative initial datum.
Let $(u_k)_{k\in\setN}$ be a sequence of $(\mu_{k,0},\mu_{k,1})$-quasi-entropy solutions according to \autoref{def:mu-quasi-entropy-solution} for the respective non-local problems
\begin{equation}\label{eq:cons-law-existence-nonlocal-k}
\partial_t u_k + \div_x\bigl(P[u_k](t,x,u_k)\bigr) = 0
\end{equation}
with non-negative initial datum $u_k(0)\in L^\infty(\setR^n)\cap BV(\setR^n)$ converging to $u_0$ in $L^1_\loc(\setR^n)$,
satisfying the uniform bounds
\[
u_k(t,\plchldr) \leq R(t), \qquad
\norm{u_k(t,\plchldr)}_{BV} \leq B(t), \qquad \forall t\in[0,T), k\in\setN,
\]
for some increasing functions $R,B:[0,T)\to[0,\infty)$.
Assume in addition that the measures $\mu_{k,i,t}\in L^1\bigl([0,T);\Meas_+(\setR^n)\bigr)$ vanish in the sense that for every $r>1/T$ we have
\begin{equation}\label{eq:vanishing-mu-k-non-local}
\lim_{k\to\infty} \int_0^{T-1/r} \abs{\mu_{k,i,t}}(B_r) \dt = 0.
\end{equation}

Assume that if $u_k$ converges in $L^\infty_\loc\bigl([0,T);L^1_\loc(\setR^n)\bigr)$ to some function $u$ then
\begin{align*}
P[u_k](t,x,s) &\to P[u](t,x,s) && \text{for a.e.\ $(t,x,s)$}, \\
\div_2P[u_k](t,x,s) &\to \div_2P[u](t,x,s) && \text{for a.e.\ $(t,x,s)$}, \\
\partial_3P[u_k](t,x,s) &\to \partial_3P[u](t,x,s) && \text{for a.e.\ $(t,x,s)$}.
\end{align*}

Let $\mathscr{F}\subseteq C^\infty\bigl([0,T)\times\setR^n;[0,\infty)\bigr)$ be a class of test functions compactly supported in space with the properties
\begin{align*}
\forall k\in\setN\ \forall(t,x)\in[0,T)\times\setR^n\ \forall\Theta\in\mathscr{F} \ &:\ 
	\partial_t\Theta(t,x)\leq-\Lip_3\bigl(P[u_k](t)\bigr)\abs{\nabla_x\Theta(t,x)}, \\
\forall(t,x)\in[0,T)\times\setR^n\ \exists\Theta\in\mathscr{F} \ &:\ 
	\Theta(t,x)>0.
\end{align*}

Assume that for every $\Theta\in\mathscr{F}$ the quantities\footnote{Recall the notation introduced in \autoref{rmk:notation-norms} depends implicitly on $\Theta$.}
\begin{gather*}
\Lip_3\bigl(P[u_k](t)\bigr), \quad
\Lip_3\bigl(\div_2 P[u_k](t)\bigr), \quad
\Lip_2\bigl(\partial_3 P[u_k](t)\bigr), \quad
\Lip_2\bigl(\div_2 P[u_k](t)\bigr), \\
\int_{\Omega_t}
	\norm{\div_2 P[u_k](t,x,\plchldr)}_{L^\infty([0,R(t)])} \dx
\end{gather*}
are equi-bounded in $L^1_\loc\bigl([0,T)\bigr)$ uniformly in $k$ and assume that there is a non-negative function $h\in L^1([0,T))$ such that
\begin{subequations}\label{eq:gronwall-non-local-h}
\begin{equation}
\int_{\setR^n} \norm{\div_2(P[u_k]-P[u_l])(t,x,\plchldr)}_{L^\infty([0,R(t)])} \Theta(t,x) \dx
\leq h(t)\int_{\setR^n} \abs{u_k(t,x)-u_l(t,x)} \Theta(t,x) \dx,
\end{equation}
\begin{equation}
\norm{\Lip_3\bigl((P[u_k]-P[u_l])(t,\plchldr)\bigr) \Theta(t,\plchldr)}_{L^\infty(\setR^n)}
\leq h(t) \int_{\setR^n} \abs{u_k(t,x)-u_l(t,x)} \Theta(t,x) \dx,
\end{equation}
\end{subequations}
hold for a.e.\ $t\in[0,T)$ and for every $k,l\in\setN$.

Then for every $t\in[0,T)$ the sequence $u_k(t,\plchldr)$ converges in $L^1_\loc(\setR^n)$ to a function $u(t,x)$ which is an entropy solution of
\begin{equation}\label{eq:cons-law-existence-nonlocal}
\partial_t u + \div_x\bigl(P[u](t,x,u)\bigr) = 0
\end{equation}
according to \autoref{def:entropy-solution} with initial datum $u_0$.
\end{theorem}

\begin{proof}[Proof of \autoref{thm:existence-fixed} and \autoref{thm:existence-nonlocal}]
For a fixed $\Theta\in\mathscr{F}$, we can apply \autoref{prop:gronwall} to the pair of quasi-solutions $(u_k,u_l)$ with $t_1=0$ and $t_2=t$, from which we deduce that 
\[
\begin{split}
\int_{\setR^n} & \abs{u_k(t,x)-u_l(t,x)} \Theta(t,x) \dx \\
&\leq \left( \int_{\setR^n} \abs{u_k(0,x)-u_l(0,x)} \Theta(0,x) \dx
	+ \Phi_{k,l}(0,t) \right) \exp\oleft(\int_0^t f_{k,l}(s)\ds\right),
\end{split}
\]
where
\begin{align*}
f_{k,l}(t) &= [3\Lip_3(\div_2P[u_k](t))+\Lip_3(\div_2P[u_l](t))] + [1+2B(T)]h(t), \\
\Phi_{k,l}(0,t) &= \norm{\Theta}_\infty \int_0^t
		(\mu_{k,0,s}+\mu_{l,0,s})\bigl(\supp\Theta(s,\plchldr)\bigr)_1 \ds
	+\frac12\norm{\nabla_x\Theta}_\infty \int_0^t
		(\mu_{k,1,s}+\mu_{l,1,s})\bigl(\supp\Theta(s,\plchldr)\bigr)_1 \ds
	\spliteq
	+C_{k,l}(0,t) \min\left\{ M_{k,l}(0,t)^{1/2}, 1 \right\}
		+c_n \norm{\Theta}_\infty
		\max\left\{ M_{k,l}(0,t)^{1/2}, M_{k,l}(0,t) \right\}, \\
C_{k,l}(0,t) &=
	\norm{\Theta}_\infty B(t) \left( 2 + \int_0^t
		[3\Lip_3(\div_2P[u_k](s))+\Lip_3(\div_2P[u_l](s))
			+2\Lip_2(\partial_3P[u_l](s))] \ds \right)
	\spliteq
	+\norm{\Theta}_{L^\infty_t L^1_x} \int_0^t \Lip_2(\div_2P[u_k](s)) \ds
	+\norm{\Theta}_{L^\infty_t\Lip_x} B(t)
		\int_0^t \Lip_3\bigl((P[u_k]-P[u_l])(s)\bigr) \ds
	\spliteq
	+\frac12\norm{\Theta}_{L^\infty_t\Lip_x}
		\int_0^t \!\!\! \int_{\bigl(\supp \Theta(s,\plchldr)\bigr)_1}
			\norm{\div_2(P[u_k]-P[u_l])(s,x,\plchldr)}_{L^\infty([0,R(s)])} \dx\ds.
\end{align*}
From the equi-boundedness in $k$ of the various norms, we get that $f_{k,l}(t)$ and $C_{k,l}(0,t)$ are uniformly bounded for every $k,l\in\setN$. Together with \eqref{eq:vanishing-mu-k-non-local}, this implies that $\Phi_{k,l}(0,t)\to0$ as $k,l\to\infty$.

Thanks to the properties of the family $\mathscr{F}$, for every radius $r>0$ we can find a finite subfamily $\mathscr{E}\subseteq\mathscr{F}$ and coefficients $c:\mathscr{E}\to\setR_+$ such that
\[
\sum_{\Theta\in\mathscr{E}} c_\Theta \Theta(t,\plchldr) \geq \bm1_{B_r}(\plchldr).
\]
Summing the previous Grönwall inequalities over $\Theta\in\mathscr{E}$ and using that $u_k(0,\plchldr)$ is a Cauchy sequence in $L^1_\loc(\setR^n)$, we deduce that $u_k(t,\plchldr)$ is a Cauchy sequence in $L^1(B_r)$, hence there exists a limit function $u(t,\plchldr)$ to which the sequence $u_k(t,\plchldr)$ converges in $L^1_\loc(\setR^n)$ for every $t\in[0,T)$. Since the Grönwall estimate is locally uniform in $t$, we have that $u_k\to u$ in $L^\infty_\loc\bigl([0,T);L^1_\loc(\setR^n)\bigr)$.

We need to show that $u$ satisfies the entropy inequality associated to the problem \eqref{eq:cons-law-existence-nonlocal}. We follow the same strategy as for the previous theorem.
By assumption, $u_k$ satisfy the $(\mu_{k,0},\mu_{k,1})$-quasi-entropy inequality for the problem \eqref{eq:cons-law-existence-nonlocal-k}. Regularizing the absolute value we deduce
\[
\begin{split}
\Intd & \left\{
	\eta_{\eps,c}(u_k)\partial_t\phi + P[u_k]_\eps(t,x,u_k,c)\cdot\nabla_x\phi
	- \phi \int_c^{u_k} \div_xP[u_k](t,x,s)\eta_{\eps,c}''(s)\d s
	\right\} \dx \dt \\
&\geq - \Intd \abs{\phi(t,x)} \d\mu_{k,0,t}(x)\dt
	- \Intd \abs{\nabla_x\phi(t,x)} \d\mu_{k,1,t}(x)\dt,
\end{split}
\]
where
\[
P[u_k]_\eps(t,x,u,c) = \int_c^u \partial_3P[u_k](t,x,s)\eta_{\eps,c}'(s) \d s.
\]
We now have to pass to the limit this inequality for $k\to\infty$.
The first term of the left hand side and the full right hand side are standard, as in the previous proof.
Moreover, we have that $P[u_k]_\eps(t,x,u_k,c) \to P[u]_\eps(t,x,u,c)$: indeed
\[
\begin{split}
&\abs{P[u_k]_\eps(t,x,u_k,c) - P[u]_\eps(t,x,u,c)} \\
&\leq \abs{P[u_k]_\eps(t,x,u_k,c) - P[u_k]_\eps(t,x,u,c)}
	+ \abs{P[u_k]_\eps(t,x,u,c) - P[u]_\eps(t,x,u,c)} \\
&\leq \Lip_3\bigl(P[u_k](t)\bigr)\abs{u_k-u}
	+ \abs*{\int_c^u (\partial_3P[u_k]-\partial_3P[u])(t,x,s)\eta_{\eps,c}'(s) \d s} \\
&\leq \Lip_3\bigl(P[u_k](t)\bigr)\abs{u_k-u}
	+ \int_0^{R(t)} \abs{\partial_3P[u_k]-\partial_3P[u]}(t,x,s) \d s,
\end{split}
\]
the first term converges to $0$ in $L^1$ and the second converges to $0$ a.e.\ and enjoys the uniform bound
\[
\abs*{\int_0^{R(t)} \abs{\partial_3P[u_k]-\partial_3P[u]}(t,x,s) \ds}
\leq \bigl[\Lip_3(P[u_k](t))+\Lip_3(P[u](t))\bigr] R(t).
\]
Finally, for the third term we decompose it as
\[
\begin{split}
&\abs*{\int_c^{u_k} \div_xP[u_k](t,x,s)\eta_{\eps,c}''(s)\ds
	- \int_c^u \div_xP[u](t,x,s)\eta_{\eps,c}''(s)\ds} \\
&\leq \abs*{\int_u^{u_k} \div_xP[u](t,x,s)\eta_{\eps,c}''(s)\ds}
	+\abs*{\int_c^{u_k} (\div_xP[u_k]-\div_xP[u])(t,x,s)\eta_{\eps,c}''(s)\ds}.
\end{split}
\]
The first integral vanishes in the limit
\[
\lim_{k\to\infty} \Intd \int_{u(t,x)}^{u_k(t,x)} \phi(t,x)\div_xP[u](t,x,s)\eta_{\eps,c}''(s)\ds\dx\dt = 0
\]
because the integrand is a fixed $L^1$ function which is integrated over the set
\[
\set{(t,x,s)\in [0,T)\times\setR^n\times[0,R(T)]}
	{(t,x)\in\supp\phi,\ s\in[u_k(t,x)\wedge u(t,x),u_k(t,x)\vee u(t,x)]}
\]
whose measure is
\[
\iint_{\supp\phi} \abs{u_k(t,x)-u(t,x)}\dx\dt \to 0, \qquad \text{for $k\to\infty$}.
\]
The second integral can be estimated as
\[
\abs*{\int_c^{u_k} (\div_xP[u_k]-\div_xP[u])(t,x,s)\eta_{\eps,c}''(s)\ds}
\leq \int_0^{R(t)} \abs{\div_xP[u_k]-\div_xP[u]}(t,x,s)\eta_{\eps,c}''(s)\ds,
\]
which goes to zero thanks to the assumptions on $\div_xP$.

After this limit $k\to0$ we get that $u$ satisfies the regularized entropy inequality
\[
\Intd \left\{
	\eta_{\eps,c}(u)\partial_t\phi + P[u]_\eps(t,x,u,c)\cdot\nabla_x\phi
	- \phi \int_c^u \div_xP[u](t,x,s)\eta_{\eps,c}''(s)\d s
	\right\} \dx \dt
\geq 0.
\]

Finally, letting $\eps\to0$ we deduce that $u$ satisfies the entropy inequality with Kruzkov entropies, and then that $u$ solves the conservation law and is continuous in time, exactly as for the previous theorem.
\end{proof}

\subsection{Uniqueness}\label{sec:uniqueness}

The uniqueness theorem for a fixed flux we would obtain from our \autoref{thm:stability} is weaker than \cite{Kruzkov} because we require the solutions to be $BV$ in addition to $L^\infty$.

On the other hand, \autoref{thm:stability} gives an interesting uniqueness result for the problem 
\[
\partial_t u(t,x) + \div_x\bigl(P[u](t,x,u(t,x))\bigr) = 0
\]
where the flux $P[u](t,x,s)$ depends non-locally on the whole solution itself, for instance $P[u](t,x,s) = v(s)(W*u)(t,x)$. The uniqueness in this setting seems to require the $BV$ regularity.

We provide two statements of the uniqueness, one dealing with solutions compactly supported in space and one with more general solutions. The reason for this distinction is that in the compact case it is sufficient to verify the assumptions with a single weight function $\Theta$ which does not need to satisfy the slope condition \eqref{eq:Theta-slope}, whereas in the general case we need to work with a suitable family of weight functions.

%
%
%

\begin{theorem}[Uniqueness with non-local fluxes, compact solutions]\label{thm:uniqueness-nonlocal-compact}
Let
\[
\begin{split}
P:\mathscr{S}\bigl([0,T),\setR^n\bigr) &\to \{p:[0,T)\times\setR^n\times[0,\infty)\to\setR\} \\
u &\mapsto P[u]
\end{split}
\]
be such that $P[u]$ is a flux satisfying (A1)\nobreakdash--(A4) of \autoref{as:flux}.

Suppose that $u,v\in\mathscr{S}\bigl([0,T),\setR^n\bigr)$ are compactly supported entropy solutions of the Cauchy problems
\[
\partial_t u + \div_x\bigl(P[u](t,x,u)\bigr) = 0, \qquad
\partial_t v + \div_x\bigl(P[v](t,x,v)\bigr) = 0
\]
with the same compactly supported initial datum $u_0\in L^\infty_\loc(\setR^n)\cap BV_\loc(\setR^n)$.

Assume that there is $\Theta\in C^1\bigl([0,T)\times\setR^n;[0,\infty)\bigr)$ compactly supported in space with $\Theta=1$ on a neighborhood of $\supp(u)\cup\supp(v)$ and a non-negative function $h\in L^1_\loc\bigl([0,T)\bigr)$ such that
\begin{subequations}
\begin{equation}
\int_{\setR^n} \norm{\div_2(P[u]-P[v])(t,x,\plchldr)}_{L^\infty([0,R(t)])} \Theta(t,x) \dx
\leq h(t)\int_{\setR^n} \abs{u(t,x)-v(t,x)} \Theta(t,x) \dx,
\end{equation}
\begin{equation}
\norm{\Lip_3\bigl((P[u]-P[v])(t,\plchldr)\bigr) \Theta(t,\plchldr)}_{L^\infty(\setR^n)}
\leq h(t) \int_{\setR^n} \abs{u(t,x)-v(t,x)} \Theta(t,x) \dx
\end{equation}
\end{subequations}
hold for a.e.\ $t\in[0,T)$.

Then $u=v$.
\end{theorem}

\begin{proof}
Recalling \autoref{rmk:entropic-statement} which allows us to omit (A5), we can apply \autoref{prop:gronwall} with $t_1=0$, $\Phi=0$ and $f(t) = [3\Lip_3(\div_2P[u](t))+\Lip_3(\div_2P[v](t))] + [1+2B(T)]h(t)$, thus obtaining
\[
\begin{split}
\int_{\setR^n} \abs{u(t_2,x)-v(t_2,x)} \dx
&= \int_{\setR^n} \abs{u(t_2,x)-v(t_2,x)} \Theta(t_2,x) \dx \\
&\leq \left( \int_{\setR^n} \abs{u(0,x)-v(0,x)} \Theta(0,x) \dx \right)
	\exp\oleft(\int_0^{t_2} f(t)\dt\right) = 0. \qedhere
\end{split}
\]
\end{proof}

\begin{theorem}[Uniqueness with non-local fluxes, non-compact solutions]\label{thm:uniqueness-nonlocal}
Let
\[
\begin{split}
P:\mathscr{S}\bigl([0,T),\setR^n\bigr) &\to \{p:[0,T)\times\setR^n\times[0,\infty)\to\setR\} \\
u &\mapsto P[u]
\end{split}
\]
be such that $P[u]$ is a flux satisfying (A1)\nobreakdash--(A4) of \autoref{as:flux}.


Suppose that $u,v\in\mathscr{S}\bigl([0,T),\setR^n\bigr)$ are entropy solutions of the Cauchy problems
\[
\partial_t u + \div_x\bigl(P[u](t,x,u)\bigr) = 0, \qquad
\partial_t v + \div_x\bigl(P[v](t,x,v)\bigr) = 0
\]
with the same initial datum $u_0\in L^\infty_\loc(\setR^n)\cap BV_\loc(\setR^n)$.


Let $\mathscr{F}\subseteq C^1\bigl([0,T)\times\setR^n;[0,\infty)\bigr)$ be a class of test functions compactly supported in space with the properties
\begin{align*}
\forall(t,x)\in[0,T)\times\setR^n\ \forall\Theta\in\mathscr{F} \ &:\ 
	\partial_t\Theta(t,x)\leq-\Lip_3\bigl(P[v](t)\bigr)\abs{\nabla_x\Theta(t,x)}, \\
\forall(t,x)\in[0,T)\times\setR^n\ \exists\Theta\in\mathscr{F} \ &:\ 
	\Theta(t,x)>0.
\end{align*}

Assume that for every $\Theta\in\mathscr{F}$ there is a non-negative function $h\in L^1_\loc\bigl([0,T)\bigr)$ such that
\begin{subequations}\label{eq:gronwall-non-compact-h}
\begin{equation}
\int_{\setR^n} \norm{\div_2(P[u]-P[v])(t,x,\plchldr)}_{L^\infty([0,R(t)])} \Theta(t,x) \dx
\leq h(t)\int_{\setR^n} \abs{u(t,x)-v(t,x)} \Theta(t,x) \dx,
\end{equation}
\begin{equation}
\norm{\Lip_3\bigl((P[u]-P[v])(t,\plchldr)\bigr) \Theta(t,\plchldr)}_{L^\infty(\setR^n)}
\leq h(t) \int_{\setR^n} \abs{u(t,x)-v(t,x)} \Theta(t,x) \dx,
\end{equation}
\end{subequations}
hold for a.e.\ $t\in[0,T)$.

Then $u=v$.
\end{theorem}

\begin{proof}
Given $\Theta\in\mathscr{F}$ and the corresponding function $h$ satisfying \eqref{eq:gronwall-non-compact-h}, recalling \autoref{rmk:entropic-statement} which allows us to omit (A5), we can apply \autoref{prop:gronwall} with $t_1=0$, $\Phi=0$ and $f(t) = [3\Lip_3(\div_2P[u](t))+\Lip_3(\div_2P[v](t))] + [1+2B(T)]h(t)$, thus obtaining
\[
\int_{\setR^n} \abs{u(t_2,x)-v(t_2,x)} \Theta(t_2,x) \dx 
\leq \left( \int_{\setR^n} \abs{u(0,x)-v(0,x)} \Theta(0,x) \dx \right)
	\exp\oleft(\int_0^{t_2} f(t)\dt\right)
= 0,
\]
hence $u(t_2,x)=v(t_2,x)$ where $\Theta(t_2,x)>0$. Letting $\Theta$ vary in $\mathscr{F}$ we get the desired conclusion.
\end{proof}

\subsection{Rate of convergence of various approximating schemes}\label{sec:rates}


%

In this subsection we study the convergence properties of some numerical schemes for the solution of scalar conservation laws. We consider schemes producing approximate solutions which are quasi-entropic and continuous in time, so that \autoref{def:mu-quasi-entropy-solution} is satisfied.

Among them, we devote more details to a recent particle method \cite{RadiciStra} applied to a non-local scalar conservation law used to model for instance traffic with congestion. In particular, we are able to derive an explicit rate of convergence (the convergence in the cited article was by compactness) and show that it is optimal (\autoref{rmk:particles-rate-sharp}). Furthermore, the Cauchy property shown in \autoref{thm:particles-rate-of-convergence} provides an independent way to prove the existence which bypasses the compactness argument used in \cite{RadiciStra}.

In addition to this, we treat the classical vanishing viscosity and front tracking methods and recover the well known convergence rates. Since our stability theorem holds also for fluxes depending on $(t,x)$, the convergence rate of the vanishing viscosity method is derived in this more general setting.

We avoid the discussion of finite difference/volume/elements such as Godunov or higher order ones since they are described by processes discrete in time. They could be studied within our framework once one constructs an interpolation in time which produces an error which is $L^1$ in time, as the right hand side of \eqref{eq:def-quasi-entropy-ineq-hyp}.
Another alternative is to generalize our stability result \autoref{thm:stability} in order to treat errors which are bounded by arbitrary measures in time instead of $\leb^1$, this however causes drastic changes to both the statement and its proof because the quasi-solutions are not necessarily continuous in time and one can no longer de-double the time variables. We leave this research direction for future work.

\subsubsection{Particle method}\label{sec:particles}

Following the series of articles \cite{DiFrancescoRosini,DiFrancescoFagioliRadici,DiFrancescoStivaletta} and concurrently with \cite{FagioliTse}, in \cite{RadiciStra} the authors study the scalar conservation law
\begin{equation}\label{eq:particles-conservation-law}
\partial_t\rho(t,x) + \div_x\bigl[\rho(t,x) v\bigl(\rho(t,x)\bigr)\bigl(V(t,x)-(\partial_xW*\rho)(t,x)\bigr)\bigr] = 0,
\end{equation}
where the convolution $\partial_xW*\rho$ is in space only,
and construct a particle based numerical scheme that produces piecewise constant approximations of the solution.

For $N\in\setN$ fixed, the piecewise constant approximation $\bar\rho^N$ is defined as
\begin{align*}
\bar\rho^N(t,x) &= \sum_{i=1}^N \rho_i^N(t) \bm{1}_{(x_{i-1}(t),x_i(t))}(x), &
\rho_i^N(t) &= \frac{1}{N(x_i(t)-x_{i-1}(t))},
\end{align*}
where the particles $X=(x_0,\dots,x_N)$ solve the ODE
\[
\tag{ODE$_I$}
\left\{\begin{aligned}
x_i'(t) &= v_i(t) \bar U_i(t), \\
\bar U_i(t) &= V\bigl(t,x_i(t)\bigr) - (\partial_xW*\bar\rho)\bigl(t,x_i(t)\bigr) \\
	&= V\bigl(t,x_i(t)\bigr) - \sum_{j=0}^N(\rho_{j+1}(t) - \rho_j(t)) W\bigl(t,x_i(t)-x_j(t)\bigr), \\
v_i(t) &= \begin{cases}
	v\bigl(\rho_i(t)\bigr), & \text{if } \bar U_i(t) < 0, \\
	v\bigl(\rho_{i+1}(t)\bigr), & \text{if } \bar U_i(t) \geq 0. \end{cases}
\end{aligned}\right.
\]

In \cite[Theorem 1.3]{RadiciStra} it is shown that the $\bar\rho^N$ converge in $L^1_\loc\bigl([0,\infty)\times\setR\bigr)$ to the entropy solution $\rho$ of \eqref{eq:particles-conservation-law}. The proof relies on a compactness argument that does not allow to establish the rate of convergence. The stability result \autoref{thm:stability} gives an alternative way to deduce the convergence $\bar\rho^N\to\rho$ together with an explicit rate.

According to \cite[Proposition~2.6, Proposition~2.10, Corollary~2.12]{RadiciStra} these functions belong to $C\bigl([0,\infty);L^1(\setR)\bigr) \cap L^\infty_\loc\bigl([0,\infty);L^\infty(\setR)\bigr) \cap L^\infty_\loc\bigl([0,\infty);BV(\setR)\bigr)$ uniformly in $N$.
In \cite[Proposition~2.13]{RadiciStra} it is shown that the piecewise constant densities satisfy an approximate entropy inequality. However, the notion adopted there differs from \autoref{def:mu-quasi-entropy-solution} because the error term is not written in integral form.
Therefore we cannot directly use \cite[(2.14)]{RadiciStra} in order to apply our stability result and must instead slightly modify the way we estimate the error terms in its proof.
The statement we can prove is the following.

\newcommand{\citeprop}{\cite[Proposition~2.13]{RadiciStra}}
\begin{proposition}[Modified version of \citeprop]\label{prop:particles-quasi-entropy}
For $N\in\setN$ let $\bar\rho^N$ be the piecewise constant density associated to the particles $X^N=(x_i)_{i=0}^N$ solving $(\mathrm{ODE}_I)$.

Then $\bar\rho^N$ is a $(\mu_0^N,\mu_1^N)$-quasi-entropy solution of \eqref{eq:particles-conservation-law} in the sense of \autoref{def:mu-quasi-entropy-solution} with $\mu_{0,t}^N = 0$ and
\[
\begin{split}
\mu_{1,t}^N
&= \sum_{i=1}^N m\bigl(\rho_i^N(t)\bigr) \abs*{\bar U(t,x)-\bar U\bigl(t,x_{i-1}(t)\bigr)}
	\bm{1}_{[x_{i-1}(t),x_i(t)]}(x) \leb^1
	\spliteq
	+\sum_{i=1}^N \rho_i^N(t) \left[
	\abs*{x'_i(t)-x'_{i-1}(t)} + \abs*{x'_{i-1}(t)
		- v\bigl(\rho_i^N(t)\bigr) \bar U\bigl(t,x_{i-1}(t)\bigr)}
	\right] \bm{1}_{[x_{i-1}(t),x_i(t)]}(x) \leb^1
\end{split}
\]
enjoying the estimate of the total mass
\[
\mu_{1,t}^N(\setR) \leq \frac1N H(t)
\]
for some increasing function $H:[0,\infty)\to[0,\infty)$ independent of $N$.
\end{proposition}

\begin{proof}
In the proof of \cite[Proposition~2.13]{RadiciStra} it is shown that
\[
\int_0^\infty \int_\setR \left\{ \abs{\bar\rho^N-c}\partial_t\phi
	+\sign(\bar\rho^N-c) \bigl[\bigl(m(\bar\rho^N)-m(c)\bigr)\bar U^N\partial_x\phi
	-m(c)\partial_x\bar U^N\phi\bigr] \right\} \d x \d t
\]
can be written as $I+\II+\III$ where $I,\II\geq0$ and\footnote{For simplicity of notation we write $\rho_i$ instead of $\rho_i^N$.}
\begin{multline*}
\III = \sum_{i=1}^N \int_0^T \sign(\rho_i-c) \biggl\{
	m(\rho_i) \biggl[ \phi(x_i) \bigl(\bar U(x_i) - \bar U(x_{i-1})\bigr)
		-\int_{x_{i-1}}^{x_i} \partial_x\bar U \phi \d x \biggr]
	\\
	+ \rho_i \Bigl[ (x'_i-x'_{i-1}) (\phi(\bar x_i) - \phi(x_i))
	- \bigl(\phi(x_i) - \phi(x_{i-1})\bigr)\bigl(x'_{i-1} - v(\rho_i)\bar U(x_{i-1})\bigr) \Bigr] \biggr\} \d t,
\end{multline*}
with $\bar x_i\in(x_{i-1},x_i)$ such that $\phi(\bar x_i)=\dashint_{x_{i-1}}^{x_i} \phi(x) \dx$.

The first line of $\III$ is the integral in time of
\[
\begin{split}
&\sum_{i=1}^N \sign(\rho_i-c)
	m(\rho_i) \biggl[ \phi(x_i) \bigl(\bar U(x_i) - \bar U(x_{i-1})\bigr)
		-\int_{x_{i-1}}^{x_i} \partial_x\bar U \phi \d x \biggr] \\
&= -\sum_{i=1}^N \sign(\rho_i-c) m(\rho_i)
	\int_{x_{i-1}}^{x_i} \partial_x\bar U(x) \bigl(\phi(x)-\phi(x_i)\bigr) \d x \\
&= \sum_{i=1}^N \sign(\rho_i-c) m(\rho_i)
	\int_{x_{i-1}}^{x_i} \bar U(x) \partial_x\phi(x) \d x
	+\sum_{i=1}^N \sign(\rho_i-c) m(\rho_i)
	\bar U(x_{i-1}) \bigl(\phi(x_{i-1})-\phi(x_i)\bigr) \\
&= \sum_{i=1}^N \sign(\rho_i-c) m(\rho_i)
	\int_{x_{i-1}}^{x_i} [\bar U(x)-\bar U(x_{i-1})] \partial_x\phi(x) \d x \\
&\geq -\sum_{i=1}^N m(\rho_i)
	\int_{x_{i-1}}^{x_i} \abs{\bar U(x)-\bar U(x_{i-1})} \abs{\partial_x\phi(x)} \d x;
\end{split}
\]
whereas the second line of $\III$ is the integral in time of
\[
\begin{split}
&\sum_{i=1}^N \sign(\rho_i-c)
	\rho_i \left[ (x'_i-x'_{i-1}) \left(\dashint_{x_{i-1}}^{x_i} \phi(x) \dx - \phi(x_i)\right)
	- \bigl(\phi(x_i) - \phi(x_{i-1})\bigr)\bigl(x'_{i-1} - v(\rho_i)\bar U(x_{i-1})\bigr) \right] \\
&=-\sum_{i=1}^N \sign(\rho_i-c)
	\rho_i \left[ \frac{x'_i-x'_{i-1}}{x_i-x_{i-1}}
	\int_{x_{i-1}}^{x_i} \int_x^{x_i} \partial_y\phi(y) \dy \dx
	+ \bigl(x'_{i-1} - v(\rho_i)\bar U(x_{i-1})\bigr)
	\int_{x_{i-1}}^{x_i} \partial_y\phi(y) \dy \right] \\
&=-\sum_{i=1}^N \sign(\rho_i-c)
	\rho_i \left[ (x'_i-x'_{i-1})
	\int_{x_{i-1}}^{x_i} \frac{y-x_{i-1}}{x_i-x_{i-1}} \partial_y\phi(y) \dy
	+ \bigl(x'_{i-1} - v(\rho_i)\bar U(x_{i-1})\bigr)
	\int_{x_{i-1}}^{x_i} \partial_y\phi(y) \dy \right] \\
&\geq -\sum_{i=1}^N \rho_i \left[
	\abs*{x'_i-x'_{i-1}} + \abs*{x'_{i-1} - v(\rho_i)\bar U(x_{i-1})}
	\right] \int_{x_{i-1}}^{x_i} \abs{\partial_y\phi(y)} \dy
\end{split}
\]
Combining these last two computations we deduce that
\[
\III \geq -\int_0^T \!\!\!\! \int_\setR \abs{\partial_x\phi} \d\mu_{1,t}^N(x) \dt
\]
where $\mu_{1,t}^N$ is the one claimed in the statement of the proposition.

Making use of some of the computations performed in \cite{RadiciStra} to bound $\III$, the total mass of $\mu_{1,t}^N$ can be estimated by
\[
\begin{split}
\mu_{1,t}^N(\setR)
&\leq \sum_{i=1}^N \rho_i v(\rho_i) (x_i-x_{i-1}) \int_{x_{i-1}}^{x_i} \abs{\partial_x\bar U} \dx
	+\sum_{i=1}^N \rho_i \left[
	\abs*{x'_i-x'_{i-1}} + \abs{x'_{i-1} - v(\rho_i) \bar U(x_{i-1})}
	\right] (x_i-x_{i-1}) \\
&\leq \frac1N \norm{v}_\infty \int_{-S(t)}^{S(t)} \abs{\partial_x\bar U} \dx
	+ \frac1N \sum_{i=1}^N \left(
	\abs*{x'_i-x'_{i-1}} + \abs{v_{i-1} - v(\rho_i)} \abs{\bar U(x_{i-1})}
	\right) \\
&\leq \frac1N \norm{v}_\infty F(t)\bigl[1+G\bigl(S(t)\bigr)+G\bigl(2S(t)\bigr)\bigr]
	\spliteq
	+ \frac1N \biggl[
	2 \norm{v}_\infty F(t)G\bigl(2S(t)\bigr)[1+R(t)]S(t)
	+8 F(t) G\bigl(R(t)\bigr) G\bigl(2S(t)\bigr) B(t)
	\biggr],
\end{split}
\]
in terms of the functions $F,G,B,R,S$ defined in \cite{RadiciStra}, which are independent of $N$.
\end{proof}

%
%

Once we have established that the piecewise constant approximations $\bar\rho^N$ are quasi-entropy solutions, we can proceed to show the converge of $\bar\rho^N$ to the exact entropy solution of \eqref{eq:particles-conservation-law}. With respect to \cite[Assumptions~1.4]{RadiciStra}, we need to require better regularity of the velocity fields $V$ and $W$ so that the resulting flux satisfies (A5) of \autoref{as:flux}.

\begin{theorem}[Cauchy property and rate of convergence]\label{thm:particles-rate-of-convergence}
Let $v,V,W,\rho_0$ be as in \cite[Theorem~1.3]{RadiciStra}. Assume in addition that $V,W \in L^1_\loc\bigl([0,\infty);W^{2,\infty}_\loc(\setR)\bigr)$. For $N\in\setN$ let $\bar\rho^N$ be the piecewise constant density associated to the particles solving $(\mathrm{ODE}_I)$, with initial datum $\bar\rho_0^N$ as in \cite[Theorem~1.3]{RadiciStra}.

Then $(\bar\rho^N)_{N\in\setN}$ is a Cauchy sequence in $L^\infty\bigl([0,T];L^1(\setR)\bigr)$ for every $T>0$. More precisely, for $M,N\in\setN$ large enough so that
$\left(\frac1M+\frac1N\right) H(T) < 1$, we have
\begin{equation}\label{eq:particles-cauchy}
\norm{\bar\rho^N(t)-\bar\rho^M(t)}_{L^1(\setR)}
\leq K(t) \norm{\bar\rho_0^N-\bar\rho_0^M}_{L^1(\setR)}
	+ L(t) \left(\frac1M+\frac1N\right)^{1/2}
\end{equation}
for some increasing functions $K,L:[0,\infty)\to[0,\infty)$ independent of $M,N$.

Moreover, the rate of convergence of $\bar\rho^N$ to the unique entropy solution $\rho$ is
\begin{equation}\label{eq:particles-rate-of-convergence}
\norm{\bar\rho^N(t)-\rho(t)}_{L^1(\setR)}
\leq K(t) \norm{\bar\rho_0^N-\rho_0}_{L^1(\setR)} + \frac{L(t)}{\sqrt N}.
\end{equation}
\end{theorem}

\begin{proof}
Under the assumptions \cite[Assumptions~1.4]{RadiciStra}, we have that the flux
\[
J^N(t,x,s) = m(s) \bigl(V(t,x)-(\partial_xW*\bar\rho^N)(t,x)\bigr)
\]
satisfies (A1)\nobreakdash--(A4) of \autoref{as:flux}. Under the additional assumptions on $V$ and $W$ we have that $J^N$ satisfies (A5) as well because
\[
\partial_x^2 J^N(t,x,s)
= m(s)[ \partial_x^2V(t,x) + (\partial_x^2W*D\bar\rho^N)(t,x)]
\in L^1_\loc\bigl([0,\infty);L^\infty_\loc(\setR)\bigr)
\]
and
\[
\partial_x\partial_s J^N(t,x,s)
= m'(s)[ \partial_xV(t,x) + (\partial_x^2W*\bar\rho^N)(t,x)]
\in L^1_\loc\bigl([0,\infty);L^\infty_\loc(\setR)\bigr).
\]

As in \cite[Proposition~2.5]{RadiciStra}, let $S(t)$ be such that $\supp\bigl(\bar\rho^N(t)\bigr)\subseteq[-S(t),S(t)]$.
Taking $\Theta\leq1$ as in \autoref{rmk:compact-support} with $\supp \Theta(t,\plchldr)\subset[-2S(t),2S(t)]$, we have
\[
\begin{split}
\int_\setR & \norm{\div_2(J^N-J^M)(t,x,\plchldr)}_{L^\infty([0,R(t)])} \Theta(t,x) \dx \\
&= \int_\setR \norm*{\partial_x\bigl[
	m(\plchldr) \bigl(\partial_xW*(\bar\rho^N-\bar\rho^M)\bigr)(t,x)
	\bigr]}_{L^\infty([0,R(t)])} \Theta(t,x) \dx \\
&\leq \norm{m}_{L^\infty([0,R(t)])} \int_{-2S(t)}^{2S(t)}
	\abs*{\bigl(D_x\partial_xW*(\bar\rho^N-\bar\rho^M)\bigr)(t,x)} \dx \\
&\leq \norm{m}_{L^\infty([0,R(t)])}
	\int_{-3S(t)}^{3S(t)} \d\abs{D_x\partial_xW(t)}
	\int_{-S(t)}^{S(t)} \abs{\bar\rho^N(t,x)-\bar\rho^M(t,x)}\dx \\
&\leq \norm{m}_{L^\infty([0,R(t)])} F(t)[1+6G(3S(t))S(t)]
	\int_\setR \abs{\bar\rho^N(t,x)-\bar\rho^M(t,x)} \Theta(t,x) \dx \dt
\end{split}
\]
and
\[
\begin{split}
&\norm{\Lip_3\bigl((J^N-J^M)(t,\plchldr)\bigr) \Theta(t,\plchldr)}_{L^\infty(\setR)} \\
&\leq \Lip(m) \norm{\bigl(\partial_xW*(\bar\rho^N-\bar\rho^M)\bigr)(t,\plchldr)\Theta(t,\plchldr)}_{L^\infty([-2S(t),2S(t)])} \\
&\leq \Lip(m) \norm{\bigl(\partial_xW*(\bar\rho^N-\bar\rho^M)\bigr)(t,\plchldr)}_{L^\infty([-2S(t),2S(t)])} \\
&\leq \Lip(m) \norm{\partial_xW(t,\plchldr)}_{L^\infty([-3S(t),3S(t)])}
	\norm{\bar\rho^N(t,\plchldr)-\bar\rho^M(t,\plchldr)}_{L^1([-S(t),S(t)])} \\
&\leq \Lip(m) F(t)G(3S(t)) \int_\setR \abs{\bar\rho^N(t,x)-\bar\rho^M(t,x)} \Theta(t,x) \dx.
\end{split}
\]
Therefore $P=J^M,Q=J^N$, $u=\bar\rho^M,v=\bar\rho^N$ satisfy \eqref{eq:gronwall-assumptions}
with
\[
h(t) = \norm{m}_{L^\infty([0,R(t)])} F(t)[1+6G(3S(t))S(t)]
	+\Lip(m) F(t)G(3S(t)).
\]
Therefore, taking into account \cite[Proposition~2.6, Proposition~2.10]{RadiciStra} and \autoref{prop:particles-quasi-entropy}, we can apply \autoref{prop:gronwall} with $t_1=0$ and $t_2=t$. As a consequence we obtain the estimate
\begin{equation}\label{eq:particles-gronwall-estimate}
\int_\setR \abs{\bar\rho^N(t,x)-\bar\rho^M(t,x)}\dx
\leq \left(\norm{\bar\rho_0^N-\bar\rho_0^M}_{L^1(\setR)}
	+ \Phi(0,t)\right) \exp\oleft(\int_0^t f(s)\ds \right).
\end{equation}
Noticing that under the assumption on $M,N$ we have
\[
\min\left\{M(0,t)^{1/2},1\right\}
= \max\left\{M(0,t)^{1/2},M(0,t)\right\}
= M(0,t)^{1/2}
\leq H(t)^{1/2} \left(\frac1M+\frac1N\right)^{1/2},
\]
the right hand side of \eqref{eq:particles-gronwall-estimate} can be bounded by
\[
\exp\oleft(\int_0^t f(s)\ds \right) \left\{
	\norm{\bar\rho_0^N-\bar\rho_0^M}_{L^1(\setR)}
	+ \left[\frac12 \norm{\nabla_2\Theta}_\infty + C(0,t) + c_n\right]
	H(t)^{1/2} \left(\frac1M+\frac1N\right)^{1/2} \right\}.
\]
The proof is concluded by estimating $f(s)$ and $C(0,t)$ in terms of the functions $F,G,R,S,B,m$ from \cite{RadiciStra} independently of $M,N$, in a similar way as we already did for $H$ and $h$.

The second part of the claim, the rate of convergence of $\bar\rho^N(t)\to\rho(t)$ in $L^1(\setR)$, is obtained either by sending $M\to\infty$ or by replicating the above Grönwall argument to the pair $(\bar\rho^N,\rho)$.
\end{proof}

\begin{remark}
Notice that the existence of $\rho$, apart from being already ensured by \cite[Theorem~1.3]{RadiciStra}, can now be deduced as a consequence of the Cauchy estimate \eqref{eq:particles-cauchy}, which implies the compactness in $C\bigl([0,T];L^1(\setR)\bigr)$ for every $T>0$.
\end{remark}

\begin{remark}\label{rmk:particles-rate-sharp}
Notice that any initial datum $\rho_0$ as in \cite[Theorem~1.3]{RadiciStra} can be approximated by $\bar\rho_0^N$ as in \cite[Lemma~1.2]{RadiciStra} with the additional property that $\norm{\bar\rho_0^N-\rho_0}_{L^1(\setR)}\leq C N^{-1/2}$, therefore \eqref{eq:particles-rate-of-convergence} as a whole is of the order $N^{-1/2}$.
Moreover, in general the initial datum $\rho_0$ cannot be approximated in $L^1(\setR)$ by some $\bar\rho_0^N$ better than $N^{-1/2}$, thus this rate is sharp.
A precise formulation of these claims is given in the next lemma.
%
\end{remark}

\begin{lemma}
Given a fixed $\rho_0\in\Prob(\setR)\cap L^\infty(\setR)\cap BV(\setR)$ with the bounds
\[
\rho_0\leq R_0, \qquad
\supp(\rho_0)\subseteq[-S_0,S_0], \qquad
\TV(\rho_0)\leq B_0,
\]
for every $N\in\setN$ there is a family of sorted particles $X^N=(x_0,\dots,x_N)$ such that the corresponding piecewise constant $\bar\rho^N$ satisfies
\[
\bar\rho^N \leq R_0, \qquad
\supp(\bar\rho^N) \subseteq [-S_0,S_0], \qquad
\TV(\bar\rho^N) \leq B_0, \qquad
\norm{\bar\rho^N - \rho_0}_{L^1(\setR)} \leq \frac{2S_0 + B_0}{\sqrt{2N}}.
\]


On the other hand, we have the following counterexample. For every positive sequence $(a_N)_{N\in\setN}$ such that $a_N=o(1/\sqrt{N})$ there exists $\rho_0\in\Prob(\setR)\cap L^\infty(\setR)\cap BV(\setR)$ with
\[
\rho_0 \leq 1, \qquad
\supp(\rho_0) \subseteq [-1,1], \qquad
\TV(\rho_0) \leq 3, \qquad
\]
such that for every sequence of sorted particles $X^N=(x_0,\dots,x_N)$ the corresponding piecewise constant $\bar\rho^N$ satisfy
\[
\limsup_{N\to\infty} \frac{\norm{\bar\rho^N - \rho_0}_{L^1(\setR)}}{a_N} = \infty.
\]
\end{lemma}

\begin{proof}
Define $x_0$ and $x_N$ be such that $[x_0,x_N]$ is the smallest interval containing $\supp(\rho_0)$, i.e.\ it is its convex hull, and then consider intermediate particles $x_i$ for $i=1,\dots,N-1$ with the property that $\rho_0([x_{i-1},x_i])=1/N$. The estimates on the $L^\infty$ norm, the support and the total variation of $\bar\rho^N$ are established in \cite[Lemma~1.2]{RadiciStra}.
Letting $l_i=x_i-x_{i-1}$ and $v_i=\TV_{(x_{i-1},x_i)}(\rho_0)$, we have
\begin{gather*}
\norm{\bar\rho^N-\rho_0}_{L^1([x_{i-1},x_i])}
	\leq l_i \left(\esssup_{(x_{i-1},x_i)}\rho_0 - \essinf_{(x_{i-1},x_i)}\rho_0 \right)
	\leq l_i v_i, \\
\norm{\bar\rho^N-\rho_0}_{L^1([x_{i-1},x_i])}
	\leq \norm{\bar\rho^N}_{L^1([x_{i-1},x_i])} + \norm{\rho_0}_{L^1([x_{i-1},x_i])}
	\leq \frac2N.
\end{gather*}
Therefore
\[
\begin{split}
\norm{\bar\rho^N-\rho_0}_{L^1(\setR)}
&\leq \sum_{i=1}^N \norm{\bar\rho^N-\rho_0}_{L^1([x_{i-1},x_i])}
\leq \sum_{i=1}^N \min\oleft\{l_iv_i,\frac2N\right\} \\
&\leq \sum_{i=1}^N \sqrt{\frac2N}\sqrt{l_iv_i}
\leq \sqrt{\frac2N} \sum_{i=1}^N \frac{l_i+v_i}2
\leq \frac{2S_0 + B_0}{\sqrt{2N}}.
\end{split}
\]

Let us now move on to the counterexample.
Given $N,C_N\in\setN_+$ define the building block
\[
\beta_N(x) = \sum_{k=0}^{C_N} \frac1{\sqrt{N}} \bm1_{[2k/\sqrt{N},(2k+1)/\sqrt{N}]}(x).
\]
We have
\[
\beta_N \leq 1, \qquad
\supp(\beta_N) \subseteq \left[0, \frac{2C_N+1}{\sqrt{N}}\right], \qquad
\TV(\beta_N) \leq 2\frac{C_N+1}{\sqrt{N}}, \qquad
\norm{\beta_N}_{L^1(\setR)} = \frac{C_N+1}N.
\]
Let $X^N=(x_0,\dots,x_N)$ be an arbitrary family of sorted particles and let $\bar\rho^N$ be its associated piecewise constant density.
Given $k=0,\dots,C_N-1$, we distinguish two cases:
\begin{itemize}
\item $\bar\rho^N(x) < 1/(2\sqrt{N})$ for some $x\in[(2k+1)/\sqrt{N},(2k+2)/\sqrt{N}]\setminus X^N$: letting $i$ be the index such that $x\in(x_{i-1},x_i)$, we must have $x_i-x_{i-1}>2/\sqrt{N}$, therefore the set
\[
\set*{y}{\bar\rho^N(y) < 1/(2\sqrt{N})} \cap \left(
	\left[\frac{2k+1/2}{\sqrt{N}},\frac{2k+1}{\sqrt{N}}\right] \cup
	\left[\frac{2k+2}{\sqrt{N}},\frac{2k+5/2}{\sqrt{N}}\right] \right)
\]
has measure at least $1/(2\sqrt{N})$, thus
\[
\norm{\bar\rho^N-\beta_N}_{L^1([(2k+1/2)/\sqrt{N},(2k+5/2)/\sqrt{N}])} \geq 1/(2\sqrt{N})^2 = 1/(4N);
\]
\item $\bar\rho^N(x) \geq 1/(2\sqrt{N})$ for every $x\in[(2k+1)/\sqrt{N},(2k+2)/\sqrt{N}]$: since $\beta_N=0$ in the interval we have
\[
\norm{\bar\rho^N-\beta_N}_{L^1([(2k+1)/\sqrt{N},(2k+2)/\sqrt{N}])}\geq 1/(2\sqrt{N}) \cdot 1/\sqrt{N} = 1/(2N).
\]
\end{itemize}
Summing over $k$ and noting that the intervals over which we have the $L^1$ estimate are essentially disjoint we get
\[
\norm{\bar\rho^N-\beta_N}_{L^1([0,(2C_N+1)/\sqrt{N}])}
\geq \sum_{k=0}^{C_N-1} \frac1{4N}
= \frac{C_N}{4N}.
\]

By assumption, $a_N N \ll \sqrt{N}$. Choose $C_N$ such that $a_NN \ll C_N \ll \sqrt{N}$, for instance $C_N = a_N^{1/2}N^{3/4}$.
The counterexample will be given by the probability
\[
\rho_0(x) = \sum_{j=1}^\infty \beta_{N_j}(x-b_{N_j})
	+ \left(1-\sum_{j=1}^\infty \frac{C_{N_j}+1}{N_j}\right)\bm1_{[-1,0]}(x),
\]
where
\[
b_{N_j} = \sum_{h=1}^{j-1} \frac{2C_{N_h}+1}{\sqrt{N_h}}
\]
is chosen so that the supports of the building blocks are disjoint and $N_j$ is a sequence growing sufficiently fast so that
\[
\sum_{j=1}^\infty \frac{2C_{N_j}+1}{\sqrt{N_j}} \leq 1, \qquad
\sum_{j=1}^{\infty} 2\frac{C_{N_j}+1}{\sqrt{N_j}} \leq 1, \qquad
\sum_{j=1}^{\infty} \frac{C_{N_j}+1}{N_j} \leq 1.
\]
Indeed $\rho_0\leq1$,
\begin{gather*}
\supp(\rho_0) \subseteq \left[-1, \sum_{j=1}^\infty \frac{2C_{N_j}+1}{\sqrt{N_j}}\right]
	\subseteq [-1, 1], \\
\TV(\rho_0) \leq 2 + \sum_{j=1}^{\infty} \TV(\beta_{N_j})
\leq 2 + \sum_{j=1}^{\infty} 2\frac{C_{N_j}+1}{\sqrt{N_j}} \leq 3,
\end{gather*}
and for any family of sorted particles $X^{N_j}=(x_0,\dots,x_{N_j})$ we have
\[
\frac{\norm{\bar\rho^{N_j}-\rho_0}_{L^1(\setR)}}{a_{N_j}}
\geq \frac{\norm{\bar\rho^{N_j} - \beta_{N_j}(\plchldr-b_{N_j})}
	_{L^1([b_{N_j},b_{N_j}+(2C_{N_j}+1)/\sqrt{N_j}])}}{a_{N_j}}
\geq \frac{C_{N_j}}{4N_j a_{N_j}} \to \infty. \qedhere
\]
%
%
%
%
%
%
%
%
%
%
\end{proof}



\subsubsection{Vanishing viscosity}\label{sec:viscosity}

The vanishing viscosity method is a way to construct solutions to the scalar conservation law
\begin{equation}\label{eq:cons-law}
\partial_tu+\div_x\bigl(P(t,x,u)\bigr)=0
\end{equation}
as the limit for $\eps\to0$ of functions $u_\eps$ that solve
\begin{equation}\label{eq:cons-law-epsDelta}
\partial_tu_\eps+\div_x\bigl(P(t,x,u_\eps)\bigr) = \eps\Delta u_\eps.
\end{equation}

In this section we demonstrate how to apply \autoref{thm:stability} to get the rate of convergence of solutions $u_\eps$ to $u$. \autoref{thm:vanishing-rate} does not address the issue of the existence of the solutions $u_\eps$ for the equation \eqref{eq:cons-law-epsDelta}. In this regard, see for instance \cite{Ladyzenskaja} which provides assumptions under which it is guaranteed.

In \cite[Theorem~5]{Kruzkov} the author uses the vanishing viscosity method to show the existence of solutions for the problem \eqref{eq:cons-law}; however, he does not provide a rate of convergence of the approximating sequence $u_\eps$ because he relies on a compactness argument.
In \cite[Theorem~3]{Kuznetsov} the author shows the rate of convergence $\sim \eps^{1/2}$ of the vanishing viscosity method for $BV$ solutions in the particular case of a flux $P(t,x,u)=P(u)$ independent of space and time. In this section we extend its result to the general case.

The following proposition shows how the approximating viscous solutions of \eqref{eq:cons-law-epsDelta} fall in this framework and can be interpreted as quasi-entropy solutions of \eqref{eq:cons-law}. This will be used in \autoref{thm:vanishing-rate} to deduce the rate of convergence.

\begin{proposition}\label{prop:vanishing-quasi-entropy}
Given $P$ satisfying \autoref{as:flux}, let
\[
u_\eps\in C\bigl([0,T);L^1(\setR^n)\bigr) \cap L^\infty_\loc\bigl([0,T);BV(\setR^n)\bigr) \cap L^\infty_\loc\bigl((0,T);H^1(\setR^n)\bigr)
\]
be a solution of \eqref{eq:cons-law-epsDelta}. Then $u_\eps$ is a $(\mu_0^\eps,\mu_1^\eps)$-quasi-entropy solution of \eqref{eq:cons-law} in the sense of \autoref{def:mu-quasi-entropy-solution} with
\begin{align*}
\mu_{0,t}^\eps &= 0, &
\mu_{1,t}^\eps &= \eps \abs{\nabla u_\eps(t,\plchldr)} \leb^n.
\end{align*}
\end{proposition}

\begin{proof}
Fix $\eta\in C^2(\setR)$ convex.
Letting
\[
P_\eta(t,x,u) = \int_0^u \partial_3P(t,x,s) \eta'(s) \ds,
\]
by the chain rule \cite[Theorem~1.5]{DeCiccoLeoni} we have the identity
\[
\begin{split}
& \div_x\bigl(P(t,x,u_\eps(t,x))\bigr) \eta'(u_\eps(t,x)) \\
&= \div_2P(t,x,u_\eps(t,x)) \eta'(u_\eps(t,x)) + \partial_3P(t,x,u_\eps(t,x))\cdot\nabla_x u_\eps(t,x) \eta'(u_\eps(t,x)) \\
&= \div_2P(t,x,u_\eps(t,x)) \eta'(u_\eps(t,x)) + \partial_3P_\eta(t,x,u_\eps(t,x))\cdot\nabla_x u_\eps(t,x) \\
&= \div_x\bigl(P_\eta(t,x,u_\eps(t,x))\bigr) + [\div_2P(t,x,u_\eps(t,x)) \eta'(u_\eps(t,x))-\div_2P_\eta(t,x,u_\eps(t,x))].
\end{split}
\]
Testing the left hand side of the equation \eqref{eq:cons-law-epsDelta} satisfied by $u_\eps$ with $\eta'(u_\eps)\phi$ and using the previous identity we get
\[
\begin{split}
&\Intd [\partial_t u_\eps + \div_x\bigl(P(t,x,u_\eps)\bigr)]\eta'(u_\eps)\phi \dx\dt \\
&= -\Intd \left\{\eta(u_\eps)\partial_t\phi + P_\eta(t,x,u_\eps)\cdot\nabla\phi
	-[\div_2P(t,x,u_\eps) \eta'(u_\eps)-\div_2P_\eta(t,x,u_\eps)] \phi \right\} \dx\dt,
\end{split}
\]
whereas for the right hand side we can estimate
\[
\begin{split}
\eps\Intd \Delta u_\eps \eta'(u_\eps)\phi\dx\dt
&= -\eps\Intd \abs{\nabla u_\eps}^2 \eta''(u_\eps)\phi \dx\dt
	-\eps\Intd \nabla u_\eps\cdot\nabla\phi \eta'(u_\eps) \dx\dt \\
&\leq \eps \norm{\eta'}_\infty \Intd \abs{\nabla\phi} \abs{\nabla u_\eps} \dx\dt.
\end{split}
\]
Given a fixed constant $c\in\setR$, in the limit as $\eta(u)$ approximates $\eta_c(u)=\abs{u-c}$, we have that $\eta'(u)$ approximates $\eta_c'(u)=\sign(u-c)$, the flux $P_\eta(t,x,u)$ approximates $P_{\eta_c}(t,x,u)=\sign(u-c)[P(t,x,u)-P(t,x,c)]$ and $[\div_2P(t,x,u) \eta'(u)-\div_2P_\eta(t,x,u)]$ approximates
\[
\begin{split}
&\div_2P(t,x,u) \sign(u-c)-\div_2P_{\eta_c}(t,x,u) \\
&= \div_2P(t,x,u) \sign(u-c)-\sign(u-c)\div_2[P(t,x,u)-P(t,x,c)]
= \sign(u-c)P(t,x,c).
\end{split}
\]
Combining the previous computations we deduce
\[
\begin{split}
\Intd & \left\{\abs{u_\eps-c}\partial_t\phi
	+ \sign(u_\eps-c)[P(t,x,u_\eps)-P(t,x,c)]\cdot\nabla\phi
	-\div_2P(t,x,c) \phi \right\} \dx\dt \\
&\geq -\eps \Intd \abs{\nabla\phi} \abs{\nabla u_\eps} \dx\dt,
\end{split}
\]
which is the definition of quasi-entropy solution with $\mu_0=0$ and $\mu_{1,t} = \eps\abs{\nabla u_\eps(t,\plchldr)}\leb^n$.
\end{proof}

\begin{theorem}[Rate of convergence of vanishing viscosity method]\label{thm:vanishing-rate}
Given $P$ satisfying \autoref{as:flux}, for every $\eps>0$ let
\[
u_\eps\in C\bigl([0,T);L^1(\setR^n)\bigr) \cap L^\infty_\loc\bigl([0,T);L^\infty(\setR^n)\bigr) \cap L^\infty_\loc\bigl([0,T);BV(\setR^n)\bigr) \cap L^\infty_\loc\bigl((0,T);H^1(\setR^n)\bigr)
\]
be a solution of \eqref{eq:cons-law-epsDelta} with initial datum $u_0\in L^1(\setR^n)\cap L^\infty(\setR^n)\cap BV(\setR^n)$.
Assume that $u_\eps$ are equi-continuous in $C\bigl([0,T);L^1(\setR^n)\bigr)$ and equi-bounded in $L^\infty_\loc\bigl([0,T);L^\infty(\setR^n)\bigr)$ and $L^\infty_\loc\bigl([0,T);BV(\setR^n)\bigr)$.

Then for $\eps\to0$ we have that $u_\eps$ converges in $C\bigl([0,T);L^1(\setR^n)\bigr)$ to
\[
u\in C\bigl([0,T);L^1(\setR^n)\bigr) \cap L^\infty_\loc\bigl([0,T);L^\infty(\setR^n)\bigr) \cap L^\infty_\loc\bigl([0,T);BV(\setR^n)\bigr)
\]
which is the unique entropy solution of \eqref{eq:cons-law} in the sense of \autoref{def:entropy-solution} with initial datum $u_0$.

Given a weight function $\Theta\in C^1\bigl([0,T)\times\setR^n;[0,\infty)\bigr)$ compactly supported in space for every time and satisfying the property\footnote{Recall that the notation $\Lip_3\bigl(P(t)\bigr)$ depends implicitly on both $\Omega_t=\bigl(\supp\Theta(t,\plchldr)\bigr)_1$ and the $L^\infty$ norm of $u$.}
\[
\partial_t\Theta(t,x)
\leq-\Lip_3\bigl(P(t)\bigr)\abs{\nabla_x\Theta(t,x)},
\qquad \forall (t,x)\in[0,T)\times\setR^n,
\]
we have the convergence rate (for $\eps<1$)
\[
\int_{\setR^n} \abs{u_\eps(t,x)-u(t,x)}\Theta(t,x)\dx
\leq \eps^{1/2} L(t), \qquad \forall t\in[0,T),
\]
for some increasing function $L:[0,T)\to[0,\infty)$ independent of $\eps$.
\end{theorem}

\begin{proof}
The claimed convergence $u_\eps\to u$ is ensured by the assumptions on $u_\eps$ and Ascoli-Arzelà Theorem. Passing to the limit the quasi-entropy inequality \eqref{eq:def-quasi-entropy-ineq-hyp} obtains \eqref{eq:def-entropy-ineq-hyp}, which says that $u$ is an entropy solution of \eqref{eq:cons-law}.

Let $B:[0,T)\to[0,\infty)$ be a non-decreasing function providing the bounds $\TV\bigl(u_\eps(t,\plchldr)\bigr),\TV\bigl(u(t,\plchldr)\bigr)\leq B(t)$.
Recalling \autoref{prop:vanishing-quasi-entropy} we have
\[
M(t_1,t_2)
= \eps \int_{t_1}^{t_2} \int_{\Omega_t}
	\abs{\nabla u_\eps(t,x)} \dx\dt
\leq \eps \int_{t_1}^{t_2} B(t)\dt
\leq \eps B(t_2)t_2.
\]
Applying \autoref{thm:stability} with $P=Q$, $\nu_0=\nu_1=0$ we get
\[
\begin{split}
&\left[\int_{\setR^n} \abs{u_\eps(t,x)-u(t,x)}\Theta(t,x)\dx\right]_{t_1}^{t_2} \\
&\leq \int_{t_1}^{t_2} 4\Lip_3(\div_2P(t))
		\int_{\setR^n} \abs{u(t,x)-v(t,x)}\Theta(t,x) \dx\dt
		+ \frac12\norm{\nabla_2\Theta}_\infty M(t_1,t_2)
	\spliteq
	+C(t_1,t_2) \min\oleft\{ M(t_1,t_2)^{1/2}, 1 \right\}
		+ c_n \norm{\Theta}_\infty
		\max\oleft\{ M(t_1,t_2)^{1/2}, M(t_1,t_2) \right\},
\end{split}
\]
where
\[
\begin{split}
C(t_1,t_2) &=
	\norm{\Theta}_\infty B(t_2) \left( 2 + \int_{t_1}^{t_2}
		[4\Lip_3(\div_2P(t))+2\Lip_2(\partial_3P(t))] \dt \right)
	\spliteq
	+\norm{\Theta}_{L^\infty_t L^1_x} \int_{t_1}^{t_2} \Lip_2(\div_2P(t)) \dt.
\end{split}
\]
Letting $t_1=0$, $t_2=t$, calling $w(t)=\int_{\setR^n} \abs{u_\eps(t,x)-u(t,x)}\Theta(t,x)\dx$ and using the fact that $u_\eps$ and $u$ have the same initial datum $u_0$ we deduce
\[
\begin{split}
w(t)
&\leq \int_0^t 4\Lip_3(\div_2P(s)) w(s) \ds
	+ \eps B(t)t \frac12\norm{\nabla_2\Theta}_\infty
	\spliteq
	+ C(0,t) \min\bigl\{\eps^{1/2}B(t)^{1/2}t^{1/2},1\bigr\}
	+ c_n \norm{\Theta}_\infty \max\bigl\{\eps^{1/2}B(t)^{1/2}t^{1/2},\eps B(t)t\bigr\}.
\end{split}
\]
For $\eps<1$ Grönwall Theorem implies then
\[
w(t)
\leq \eps^{1/2}
	\left(\frac12\norm{\nabla_2\Theta}_\infty + C(0,t) + c_n \norm{\Theta}_\infty\right)
	\bigl[B(t)t+B(t)^{1/2}t^{1/2}\bigr]
	\exp\oleft(\int_0^t 4\Lip_3(\div_2P(s)) \ds\right).
	\qedhere
\]
\end{proof}

Similarly to what we did in \autoref{thm:particles-rate-of-convergence}, a variation of our presented argument where we apply the stability theorem to $u_{\eps_1}$ and $u_{\eps_2}$ directly implies that $(u_{\eps_n})_{n\in\setN}$ is a Cauchy sequence in $L^\infty_\loc\bigl([0,T),L^1_\loc(\setR^n)\bigr)$ whenever $\eps_n\to0$. This is an alternative way to deduce the existence of a solution for the limiting problem \eqref{eq:cons-law} which does not rely on the equi-continuity and the compactness arguments.

%
%

\subsubsection{Front tracking}\label{sec:front-tracking}

Given $f\in C^2(\setR)$, consider the flux $P(t,x,u)=f(u)$. The front tracking is a scheme introduced by \cite{DafermosPolygonal,Holden} to solve the conservation law
\begin{equation}\label{eq:pde-front}
\partial_t u + \partial_x f(u) = 0.
\end{equation}
For $\nu\in\setN$, define the piecewise linear function $f_\nu:\setR\to\setR$ which interpolates $f$ on the grid $2^{-\nu}\setZ$. The approximating functions $u_\nu$ considered by the front tracking are the entropy solutions of the conservation law associated to the modified flux
\begin{equation}\label{eq:pde-front-piecewise}
\partial_t u_\nu + \partial_x f_\nu(u_\nu) = 0.
\end{equation}
The initial datum $u_{\nu,0}$ used in combination with \eqref{eq:pde-front-piecewise} is a discretization of the initial datum $u_0$ of \eqref{eq:pde-front} taking values only in $2^{-\nu}\setZ$.

Specializing \autoref{thm:stability} to the case of a problem independent of space and time, we are able to recover the well-known convergence rate \cite{Holden}.

Instead of considering $u_\nu$ as a $(\mu_0,\mu_1)$-quasi-entropy solution of the same problem \eqref{eq:pde-front} to which $u$ is an exact entropy solution with $\mu_0=0$ and $\mu_1=2\norm{f-f_\nu}_{L^\infty(\setR)}\leb$, a better stability estimate is obtained exploiting the fact that $u_\nu$ is an exact entropy solution of the nearby problem \eqref{eq:pde-front-piecewise}. The net result is that the estimate becomes independent of the measure of the support of the weight function $\Theta$.

\begin{theorem}[Rate of convergence of front tracking method]
Let $f\in C^2(\setR)$ and $u_0\in L^\infty(\setR)\cap BV(\setR)$ with $\norm{u_0}_{L^\infty}\leq R_0$ and $\TV(u_0)\leq B_0$. For $\nu\in\setN$, let $f_\nu:\setR\to\setR$ be the piecewise linear function which interpolates $f$ on the grid $2^{-\nu}\setZ$.
Let $u_\nu$ be the entropy solution of the problem \eqref{eq:pde-front-piecewise} with an initial datum such that $\mathop{\mathrm{Im}}(u_{\nu,0})\subseteq 2^{-\nu}\setZ$, $\norm{u_{\nu,0}}_{L^\infty}\leq R_0$, $\TV(u_{\nu,0})\leq B_0$ and $u_{\nu,0} \to u_0$ in $L^1(\setR)$. Then
\[
\begin{split}
\norm{u_\nu(t)-u(t)}_{L^1(\setR)} \leq \norm{u_{\nu,0}-u_0}_{L^1(\setR)}
	+ 2^{-\nu} B_0 \norm{f}_{C^2([-R_0,R_0])} t,
\end{split}
\]
where $u$ is the unique entropy solution of \eqref{eq:pde-front} with initial datum $u_0$.
\end{theorem}

\begin{proof}
By construction we have $\Lip_{[-R_0,R_0]}(f_\nu-f) \leq \norm{f}_{C^2([-R_0,R_0])} 2^{-\nu-1}$.

In \cite[Section~14.1]{Dafermos} it is shown that $\norm{u_\nu(t)}_{L^\infty} \leq \norm{u_{\nu,0}}_{L^\infty}\leq R_0$ and $\TV(u_\nu(t))\leq \TV(u_{\nu,0})\leq B_0$ for every $t\in[0,\infty)$.

Taking $\Theta\in C^1\bigl([0,\infty)\times\setR;[0,1]\bigr)$ compactly supported in space and satisfying $\partial_t\Theta(t,x) \leq -\Lip(f) \abs{\partial_x\Theta(t,x)}$, we can apply \autoref{thm:stability} with \autoref{rmk:entropic-statement}. Observing that the terms $\div_2P,\div_2Q=0$ because the fluxes are independent of space, we deduce
\[
\begin{split}
\left[\int_{\setR} \abs{u_\nu(t,x)-u(t,x)} \Theta(t,x) \dx \right]_{t_1}^{t_2}
&\leq 2B_0 \int_{t_1}^{t_2} \Lip_{[-R_0,R_0]}(f_\nu-f) \norm{\Theta(t,\plchldr)}_\infty \dt \\
&\leq 2^{-\nu} B_0 \norm{f}_{C^2([-R_0,R_0])} (t_2-t_1).
\end{split}
\]
Taking the limit as $\Theta$ approximates the constant function identically equal to $1$ we get the thesis.
\end{proof}

%
%


\phantomsection
\addcontentsline{toc}{section}{\refname}
\printbibliography

@book{AFP,
	Author = {Ambrosio, Luigi and Fusco, Nicola and Pallara, Diego},
	Title = {Functions of bounded variation and free discontinuity problems},
	isbn={9780198502456},
	lccn={99046602},
	series={Oxford Science Publications},
	year={2000},
	publisher={Clarendon Press}
}

@article{BouchutPerthame1998,
  title={Kruzkov's estimates for scalar conservation laws revisited},
  author={François Bouchut and Beno{\^i}t Perthame},
  journal={Transactions of the American Mathematical Society},
  year={1998},
  volume={350},
  pages={2847-2870}
}

@article{DaneriRadiciRunaLinear,
title = {Deterministic particle approximation of aggregation-diffusion equations on unbounded domains},
journal = {Journal of Differential Equations},
volume = {312},
pages = {474-517},
year = {2022},
issn = {0022-0396},
doi = {https://doi.org/10.1016/j.jde.2021.12.019},
url = {https://www.sciencedirect.com/science/article/pii/S0022039621007890},
author = {Sara Daneri and Emanuela Radici and Eris Runa},
}

@online{DaneriRadiciRunaNonlinear,
  author = {Daneri, Sara and Radici, Emanuela and Runa, Eris},
  title = {Deterministic particle approximation of aggregation diffusion equations with nonlinear mobility},
  year = {2022},
  eprinttype = {arXiv},
  eprint = {2209.10884},
  primaryClass = {math.AP},
  url = {https://arxiv.org/abs/2209.10884},
}

@article{FagioliTse,
title = {On gradient flow and entropy solutions for nonlocal transport equations with nonlinear mobility},
journal = {Nonlinear Analysis},
volume = {221},
pages = {112904},
year = {2022},
issn = {0362-546X},
doi = {https://doi.org/10.1016/j.na.2022.112904},
url = {https://www.sciencedirect.com/science/article/pii/S0362546X22000785},
author = {Simone Fagioli and Oliver Tse},
}

@book{Dafermos,
	title = {Hyperbolic Conservation Laws in Continuum Physics},
	author = {Constantine M. Dafermos},
	publisher = {Springer Berlin, Heidelberg},
	series = {Grundlehren der mathematischen Wissenschaften},
	year = {2016},
	doi = {https://doi.org/10.1007/978-3-662-49451-6},
	isbn = {978-3-662-49451-6}
}

@article{DafermosPolygonal,
title = {Polygonal approximations of solutions of the initial value problem for a conservation law},
journal = {Journal of Mathematical Analysis and Applications},
volume = {38},
number = {1},
pages = {33-41},
year = {1972},
issn = {0022-247X},
doi = {https://doi.org/10.1016/0022-247X(72)90114-X},
url = {https://www.sciencedirect.com/science/article/pii/0022247X7290114X},
author = {Constantine M Dafermos}
}

@article{DeCiccoLeoni,
author="De Cicco, Virginia
and Leoni, Giovanni",
title="A chain rule in $L^1(\mathrm{div};\Omega)$ and its applications to lower semicontinuity",
journal="Calculus of Variations and Partial Differential Equations",
year="2003",
day="01",
volume="19",
number="1",
pages="23--51",
issn="1432-0835",
doi="10.1007/s00526-003-0192-2",
url="https://doi.org/10.1007/s00526-003-0192-2"
}

@article{DiFrancescoStivaletta,
	title = {Convergence of the follow-the-leader scheme for scalar conservation laws with space dependent flux},
	volume = {40},
	url = {https://www.aimsciences.org/article/doi/10.3934/dcds.2020010},
	doi = {10.3934/dcds.2020010},
	pages = {233-266},
	number = {1},
	journaltitle = {Discrete \& Continuous Dynamical Systems},
	author = {Di Francesco, Marco and Stivaletta, Graziano},
	year = {2020}
}

@article{DiFrancescoRosini,
	title = {Rigorous Derivation of Nonlinear Scalar Conservation Laws from Follow-the-Leader Type Models via Many Particle Limit},
	volume = {217},
	issn = {1432-0673},
	url = {https://doi.org/10.1007/s00205-015-0843-4},
	doi = {10.1007/s00205-015-0843-4},
	pages = {831--871},
	number = {3},
	journaltitle = {Archive for Rational Mechanics and Analysis},
	author = {Di Francesco, M. and Rosini, M.D.},
	year = {2015}
}

@article{DiFrancescoFagioliRadici,
	title = {Deterministic particle approximation for nonlocal transport equations with nonlinear mobility},
	volume = {266},
	issn = {0022-0396},
	url = {https://www.sciencedirect.com/science/article/pii/S0022039618305102},
	doi = {10.1016/j.jde.2018.08.047},
	pages = {2830--2868},
	number = {5},
	journaltitle = {Journal of Differential Equations},
	author = {Di Francesco, Marco and Fagioli, Simone and Radici, Emanuela},
	year = {2019},
}

@article{FagioliRadiciDiffusion,
	title = {Solutions to aggregation–diffusion equations with nonlinear mobility constructed via a deterministic particle approximation},
	volume = {28},
	issn = {0218-2025},
	url = {https://www.worldscientific.com/doi/abs/10.1142/S0218202518400067},
	doi = {10.1142/S0218202518400067},
	pages = {1801--1829},
	number = {9},
	journaltitle = {Mathematical Models and Methods in Applied Sciences},
	author = {Fagioli, Simone and Radici, Emanuela},
	year = {2018}
}

@article{Holden,
title = {A numerical method for first order nonlinear scalar conservation laws in one-dimension},
journal = {Computers \& Mathematics with Applications},
volume = {15},
number = {6},
pages = {595-602},
year = {1988},
issn = {0898-1221},
doi = {https://doi.org/10.1016/0898-1221(88)90282-9},
url = {https://www.sciencedirect.com/science/article/pii/0898122188902829},
author = {H. Holden and L. Holden and R. Høegh-Krohn}
}

@article{KarlsenRisebro,
	Author = {Karlsen, Kenneth Hvistendahl and Risebro, Nils Henrik},
	Title = {On the uniqueness and stability of entropy solutions of nonlinear degenerate parabolic equations with rough coefficients},
	Journal = {Discrete and Continuous Dynamical Systems},
	Volume = {9},
	Number = {5},
	Year = {2003},
	Month = {09},
	Pages = {1081-1104},
	Doi = {10.3934/dcds.2003.9.1081}
}

@article{Kruzkov,
	Author = {Kru\v{z}kov, S. N.},
	Title = {First order quasilinear equations with several independent variables},
	Journal = {Matematicheskiĭ Sbornik},
	Volume = {81},
	Number = {123},
	Year = {1970},
	Pages = {228-255},
	MRCLASS = {35.37},
	MRNUMBER = {0267257},
}

@article{Kuznetsov,
title = {Accuracy of some approximate methods for computing the weak solutions of a first-order quasi-linear equation},
journal = {USSR Computational Mathematics and Mathematical Physics},
volume = {16},
number = {6},
pages = {105-119},
year = {1976},
issn = {0041-5553},
doi = {https://doi.org/10.1016/0041-5553(76)90046-X},
url = {https://www.sciencedirect.com/science/article/pii/004155537690046X},
author = {N.N. Kuznetsov}
}

@book{Ladyzenskaja,
	author = {Ladyzhenskaya, Olga A. and  V. A. Solonnikov and Nina Nikolaevna Uraltseva},
	year = {1968},
	title = {Linear and quasilinear equations of parabolic type},
	publisher = {American Mathematical Society, Providence, R.I.},
	series = {Transl. Math. Monographs},
	volume = {23}
}

@online{RadiciStra,
  author = {Radici, Emanuela and Stra, Federico},
  title = {Entropy solutions of non-local scalar conservation laws with congestion via deterministic particle method},
  eprinttype = {arxiv},
	eprint = {2012.01966},
	primaryClass={math.AP},
  url = {https://arxiv.org/abs/2107.10760},
  year = {2021},
  copyright = {Creative Commons Attribution Non Commercial No Derivatives 4.0 International},
  note = {Accepted for publication on SIAM SIMA.}
}

@article{VolpertBV,
	title = {The spaces $BV$ and quasilinear equations},
	volume = {2},
	issn = {0025-5734},
	url = {https://iopscience.iop.org/article/10.1070/SM1967v002n02ABEH002340/meta},
	doi = {10.1070/SM1967v002n02ABEH002340},
	pages = {225--267},
	number = {2},
	journal = {Mathematics of the {USSR}-Sbornik},
	publisher = {{IOP} Publishing},
	author = {Vol'pert, A. I.},
	year = {1967},
	langid = {english},
}

@article{VolpertParabolic,
	doi = {10.1070/sm1969v007n03abeh001095},
	url = {https://doi.org/10.1070/sm1969v007n03abeh001095},
	year = 1969,
	publisher = {{IOP} Publishing},
	volume = {7},
	number = {3},
	pages = {365--387},
	author = {A. I. Vol'pert and S. I. Hudjaev},
	title = {Cauchy's problem for degenerate second order quasilinear parabolic equations},
	journal = {Mathematics of the {USSR}-Sbornik}
}

\end{document}